\documentclass[a4paper]{article}
\usepackage{lmodern}

\usepackage{stmaryrd,amsthm,amsmath,amssymb,mathrsfs,hyperref}

\newtheorem{prop}{Proposition}
\newtheorem{thm}[prop]{Theorem}
\newtheorem{cor}[prop]{Corollary}
\newtheorem{lem}[prop]{Lemma}

\theoremstyle{remark}
\newtheorem{rem}[prop]{Remark}

\newcommand{\gen}[1]{\langle#1\rangle}

\newcommand{\normal}{\trianglelefteq}

\newcommand{\PSL}{\mathrm{PSL}}
\newcommand{\SL}{\mathrm{SL}}
\newcommand{\PGL}{\mathrm{PGL}}
\newcommand{\GL}{\mathrm{GL}}
\newcommand{\PSU}{\mathrm{PSU}}
\newcommand{\PGU}{\mathrm{PGU}}
\newcommand{\PSp}{\mathrm{PSp}}
\newcommand{\SU}{\mathrm{SU}}
\newcommand{\Sp}{\mathrm{Sp}}
\newcommand{\POmega}{\mathrm{P}\Omega}
\newcommand{\Spin}{\mathrm{Spin}}

\newcommand{\Alt}{\mathrm{Alt}}
\newcommand{\Sym}{\mathrm{Sym}}

\newcommand{\C}{\mathbb{C}}
\newcommand{\F}{\mathbb{F}}
\newcommand{\Q}{\mathbb{Q}}
\newcommand{\I}{\mathrm{i}}

\newcommand{\cf}{\mathrm{cf}}
\newcommand{\mbG}{\mathbf{G}}
\newcommand{\mb}[1]{\mathbf{#1}}
\newcommand{\ms}[1]{\mathscr{#1}}

\newcommand{\mc}[1]{\mathcal{#1}}

\newcommand{\Aut}{\mathrm{Aut}}
\newcommand{\Out}{\mathrm{Out}}
\newcommand{\Hom}{\mathrm{Hom}}
\newcommand{\Ext}{\mathrm{Ext}}
\newcommand{\Outdiag}{\mathrm{Outdiag}}

\numberwithin{prop}{section}
\numberwithin{equation}{section}
\title{The Maximal Subgroups of the Exceptional Groups $F_4(q)$, $E_6(q)$ and ${}^2\!E_6(q)$ and Related Almost Simple Groups}
\author{David A.\ Craven, University of Birmingham}
\date{\today}

\begin{document}
\maketitle
\begin{abstract}\textbf{This version is an updated version of the paper that appeared in \emph{Invent.\ Math.}. It incorporates all known corrections to the tables, which are listed at the end of this paper.}

\medskip

This article produces a complete list of all maximal subgroups of the finite simple groups of type $F_4$, $E_6$ and twisted $E_6$ over all finite fields. Along the way, we determine the collection of Lie primitive almost simple subgroups of the corresponding algebraic groups. We give the stabilizers under the actions of outer automorphisms, from which one can obtain complete information about the maximal subgroups of all almost simple groups with socle one of these groups. We also provide a new maximal subgroup of ${}^2\!F_4(8)$, correcting the maximal subgroups for that group from the list of Malle. This provides the first new exceptional groups of Lie type to have their maximal subgroups enumerated for three decades. The techniques are a mixture of algebraic groups, representation theory, computational algebra, and use of the trilinear form on the $27$-dimensional minimal module for $E_6$. We provide a collection of supplementary Magma files that prove the author's computational claims, yielding existence and the number of conjugacy classes of all maximal subgroups mentioned in the text.

\let\thefootnote\relax\footnotetext{\noindent Funding: Supported by a Royal Society University Research Fellowship.\\Code availability: all Magma programs are available on the arXiv, on the author's website, and submitted to the journal.}
\end{abstract}

\setcounter{tocdepth}{2}
\tableofcontents

\section{Introduction}
In order to use the classification of the finite simple groups, having names for the simple groups is usually not enough. Frequently, detailed information about their structure is required. One such piece of structure that appears often---for example in applications to permutation groups, regular graphs, and so on---is knowledge of the maximal subgroups of a finite simple group.

On the other hand, knowledge of the maximal subgroups of \emph{all} finite groups is equivalent to knowledge of the maximal subgroups and $1$-cohomology groups of all (almost) simple groups \cite{aschbacherscott1985}. Thus the maximal subgroups of (almost) simple groups take centre stage.

The maximal subgroups of the sporadic simple groups are all known except for a few cases for the Monster. For the alternating groups, the O'Nan--Scott theorem lists the maximal subgroups (up to conjugacy), with several structural classes and one collection of `leftover' cases, often labelled $\mc S$, which consist of almost simple primitive subgroups. The same occurs with the classical groups, where Aschbacher's theorem \cite{aschbacher1984} (see also \cite{bhrd,kleidmanliebeck}) outlines several geometric classes, and then has a final class of almost simple groups acting absolutely irreducibly, about which much is known but still much remains unknown.

For exceptional groups in characteristic $p$, a similar theorem classifying maximal subgroups holds \cite{borovik1989,liebeckseitz1990}: there are several families that arise from positive-dimensional subgroups of the corresponding exceptional algebraic group; a set of `exotic $r$-local subgroups', normalizers of certain $r$-subgroups for $r\neq p$; the Borovik subgroup, a particular subgroup of $E_8(q)$; a set of `subfield subgroups', which are the fixed points of outer automorphisms of the finite group (for $F_4(q)$ they are $F_4(q_0)$ if $q_0=q^r$ for $r$ a prime and ${}^2\!F_4(q)$ if $q$ is an odd power of $2$, for $E_6(q)$ they are $E_6(q_0)$ and also ${}^2\!E_6(\sqrt q)$ if $q$ is a square, and for ${}^2\!E_6(q)$ they are ${}^2\!E_6(q_0)$ if $q_0=q^r$ for $r$ an odd prime); and a set $\mc S$ of almost simple maximal subgroups that do not lie in any of the above sets. (See later in this introduction for more details about the history of the study of maximal subgroups of these groups.)

The exact set $\mc S$ has been determined for small-rank exceptional groups: the Suzuki groups ${}^2\!B_2(q)$, $G_2(q)$ and ${}^2\!G_2(q)$ \cite{kleidman1988,cooperstein1981}, ${}^3\!D_4(q)$ \cite{kleidman1988b}, and ${}^2\!F_4(q)$ \cite{malle1991}. Unpublished work of Magaard \cite{magaardphd} and Aschbacher \cite{aschbacherE6Vun} made considerable progress on determining $\mc S$ for $F_4(q)$ and $E_6(q)$ respectively, but there were numerous remaining candidates. (The papers \cite{aschbacherE6Vun,magaardphd} do produce a list of subgroups such that every element of $\mc S$ appears on the list, although whether each element yields maximal subgroups, and how many conjugacy classes, was left undetermined for some of them.)

This paper makes the first complete determination of $\mc S$ for a family of exceptional groups in three decades. We obtain a complete solution for $F_4(q)$, $E_6(q)$ and ${}^2\!E_6(q)$. Prior to this paper, a complete answer was only known for $F_4(2)$ \cite{nortonwilson1989}, $E_6(2)$ \cite{kleidmanwilson1990a} and ${}^2\!E_6(2)$ \cite[p.191]{atlas}. (No proof is given there, but there is a proof given recently by Wilson \cite{wilson2018un}, which appeared as this work was underway. The proof here depends on neither \cite{atlas} nor \cite{wilson2018un}.) Note also that our proof depends on neither \cite{aschbacherE6Vun} nor \cite{magaardphd} (although we use ideas from \cite{aschbacherE6Vun} in Section \ref{sec:sl28}), and thus all parts of the determination of the maximal subgroups of these groups, including any dependencies, are published.

\begin{thm}\label{thm:maxsubgroups} Let $\bar G$ be an almost simple group with socle $G$ one of $F_4(q)$, $E_6(q)$ and ${}^2\!E_6(q)$. All maximal subgroups of $\bar G$ are known, and are given in Tables \ref{tab:f4curlyS}, \ref{tab:f4othermaximalsqodd} and \ref{tab:f4othermaximalsqeven} for $F_4(q)$, Tables \ref{tab:e6curlyS} and \ref{tab:e6othermaximals} for $E_6(q)$, and Tables \ref{tab:2e6curlyS} and \ref{tab:2e6othermaximals} for ${}^2\!E_6(q)$.
\end{thm}

Tables \ref{tab:f4curlyS} to \ref{tab:2e6curlyS} give the members of $\mathcal S$ for these groups, and is what is done in this paper. The other maximal subgroups, in Tables \ref{tab:f4othermaximalsqodd} to \ref{tab:2e6othermaximals}, were already known, either directly via the work of Magaard and Aschbacher, or as part of a wide-ranging programme that classifies all maximal subgroups outside $\mathcal S$ for all exceptional groups of Lie type,\footnote{This is not technically true in one particular case: for the simple groups $E_7(q)$ for $p\geq 5$, there are simple subgroups $\PSL_3(q)$ and $\PSU_3(q)$, and their normalizers in the simple group are not determined in \cite{liebeckseitz2004}, merely in the adjoint group. This remains unresolved in \cite{craven2021un}, but all other maximal subgroups outside of $\mathcal S$ are known for all exceptional groups of Lie type, in particular including the groups we consider here.} primarily by Liebeck and Seitz; see Section \ref{sec:remainingmaximals} for full details. (The format of these tables follows that of \cite{bhrd}. The first column gives the subgroup of the simple group $G$, the second the prime $p$, the third the prime power $q$, the fourth the number of classes in the simple group, and the fifth the stabilizer of each class in $\Out(G)$. With this information one can compute the normalizer and number of classes in any almost simple group $\bar G$ with socle $G$, and of course if it is maximal in $\bar G$.)

As a remark, these tables make clear that there is a version of Ennola duality for subgroups. Recall that Ennola duality is broadly the principle that information about representations of a twisted group of Lie type at $q$ can be obtained from the untwisted group at `$-q$', at least in terms of numerical and combinatorial information. Tables \ref{tab:e6curlyS} and \ref{tab:2e6curlyS} display the same duality. This is explored more in \cite{craven2022}.

\bigskip

We should also mention that, in the process of determining the maximal subgroups of $F_4(q)$ for $q$ even, a new maximal subgroup of the large Ree groups was found. (see Remark \ref{rem:psl213}). There are exactly three conjugacy classes of maximal subgroups $\PGL_2(13)$ in ${}^2\!F_4(8)$, permuted transitively by the field automorphism. We explicitly mention this as this subgroup was erroneously excluded in \cite{malle1991}. In addition, another error came to light during the production of these results, in the main theorem of \cite{cohenwales1997}: the Lie primitive subgroups $\PSL_2(11)$ in $F_4(\C)$ are missing from their list (see Section \ref{sec:e6;psl211}, and Tables \ref{tab:e6curlyS} and \ref{tab:2e6curlyS}).

All of the subgroups in $\mc S$ have previously been constructed in the algebraic group, with many in \cite{cohenwales1997}, with the exceptions of $\SL_2(8)$ in $F_4$ for $p=7$ (which was erroneously ruled out in \cite{magaardphd}) and $\PSL_2(11)$ in $E_6$ (which was erroneously ruled out in both \cite{aschbacherE6Vun} and \cite{cohenwales1997}). What is new is in many cases the structure of the normalizer of the simple group, and in almost all cases the number of conjugacy classes, the actions of the outer automorphisms, and the specific finite groups into which these subgroups embed.

\medskip

The following is a somewhat abbreviated history of the classification of the maximal subgroups of the simple groups in question. In the same vein as the work of Kleidman \cite{kleidman1988,kleidman1988b}, Aschbacher and his student Magaard set about using the geometry of the trilinear form on the minimal module for $E_6$ to obtain a list of the maximal subgroups of the exceptional groups $F_4(q)$ and $E_6(q)$, more or less entirely by hand. Magaard's thesis \cite{magaardphd} was never published, and although Aschbacher produced four papers on $E_6$, starting in this journal more than 30 years ago with \cite{aschbacher1987}, the fifth \cite{aschbacherE6Vun}, which attempted a full classification, was never published either. Both of these papers leave open the question of whether several almost simple groups were maximal subgroups, and the number of conjugacy classes. At the same time, Cohen and Wales \cite{cohenwales1997}, using hand and computational techniques, started to construct various subgroups of $F_4(\C)$ and $E_6(\C)$, which form the bulk of the maximal subgroups in Tables \ref{tab:f4curlyS} to \ref{tab:2e6curlyS}.

Later, using techniques from algebraic groups, a number of authors, most notably Liebeck and Seitz, set about classifying the maximal subgroups of \emph{all} exceptional groups. These classifications were more broad brush than Aschbacher--Magaard's, but also dealt with $E_7$ and $E_8$, which are not amenable to the geometric methods of Aschbacher--Magaard. The first milestone of this approach (although it is later chronologically) is \cite{liebeckseitz2004}, which gave the positive-dimensional maximal subgroups of the exceptional algebraic groups, and furnished us with `most' of the maximal subgroups of the finite exceptional groups $G$. A result of Borovik \cite{borovik1989} and Liebeck--Seitz \cite{liebeckseitz1990} (together with the lists from \cite{clss1992,liebecksaxlseitz1992}) reduced the problem of classifying maximal subgroups of $G$ to classifying almost simple maximal subgroups $H$.

If $H$ is Lie type in the same characteristic as $G$ (called a `generic' maximal subgroup) then results of a number of authors (see \cite{liebeckseitz1998} and the references therein) prove that either $H=H(q)$ is of semisimple rank at most half of that of $G$ and $q$ is `small' (at most $248$ for $E_6$), or $H$ is the fixed points of a maximal positive-dimensional subgroup of the corresponding algebraic group (and as mentioned earlier, these are known). For $G$ one of $F_4$, $E_6$ and ${}^2\!E_6$, in \cite{craven2015un,craven2019un} the former case was removed, and so all such maximal subgroups for these groups are known.

A complete list of all `non-generic' almost simple groups that embed in exceptional algebraic groups was obtained by Liebeck and Seitz in \cite{liebeckseitz1999}, although most of these cannot be maximal. Using this list, Litterick \cite{litterickmemoir} computed the possible composition factors of one of these on both the minimal and adjoint modules, and hence eliminated several options for a maximal subgroup. More importantly for us here, \cite{litterickmemoir} eliminates many possible actions, leaving very few sets of composition factors for the action of $H$ on minimal and adjoint modules. The author considered the case where $H$ is alternating in \cite{craven2017}. This paper considers all other options for $H$, using the lists in \cite{litterickmemoir} as a jumping-off point.

\medskip

A few words about our techniques, particularly the use of computers. Theoretical approaches, most notably Aschbacher's work \cite{aschbacherE6Vun}, make the most headway when the action of the potential maximal subgroup on the minimal module is `easy', in the sense that it is a relatively simple task to write down the action of generators of the group on basis elements of the space, or it stabilizes some geometrically important subspace. One of the highlights of \cite{aschbacherE6Vun}, which we cannot reproduce using purely computational ideas, is understanding the action of an outer diagonal automorphism of $E_6(q)$ on classes of $\SL_2(8)$ subgroups. The group $\SL_2(8).3\cong {}^2\!G_2(3)$ acts on the minimal module for $E_6$ as the Steinberg representation, so element actions are very combinatorial. This allows a detailed examination of the invariant symmetric trilinear forms. Once one moves to other Lie primitive groups, for example $\PSL_2(11)$, which just acts as a sum of simple modules of dimensions $5$, $10$ and $12$, it looks to be a hard problem to compute the forms. The simple modules are now harder to write down, and in effect we end up using a computer anyway. Using Aschbacher's methods to resolve the dozen or so missing groups from \cite{aschbacherE6Vun}---by hand---appears a daunting task.

We work with the trilinear form here as well when we employ the most commonly used technique in this paper. We choose a copy $G$ of (the simply connected version of) $E_6(k)$ inside $\GL_{27}(k)$, and choose a copy of our putative subgroup $H$ of $G$, also inside $\GL_{27}(k)$. We find some subgroup $L$ of $H$ that we can already classify conjugacy classes of inside $G$, and arrange our copy of $H$ so that $L\leq G$.

If $N$ denotes the normalizer of $L$ inside $\GL_{27}(k)$, then we will determine which $n\in N$ satisfy $H^n\leq G$. Using the trilinear form $f$ corresponding to $G$, this reduces to solving systems of polynomial equations. Producing the polynomials requires explicit descriptions of $H$, $L$, $N$ and $f$; it would be essentially impossible to do this by hand for most groups $H$, although solving the resulting systems of polynomials by hand is certainly possible for the groups we consider here, as equations are not particularly complicated. Section \ref{sec:trilinear} describes the process by which the polynomials are generated.

The computer use is more or less limited to determining the number of conjugacy classes of various subgroups in the algebraic group. For the finite groups and stability under outer automorphisms, the calculations are done by hand. All computations mentioned here were performed using Magma \cite{magma}. Difficult or lengthy computations are completed in the supplementary materials, and are described in this way. Shorter computations, for example counting classes of subgroups of certain Weyl groups, which are easy to perform in Magma, are simply described as `a computer check' or similar words.

\medskip

This paper is structured as follows. The next section gathers our notation and gives some general preliminary results on algebraic groups. Section \ref{sec:techniques} gives an overview of the various techniques that are used in this article. Some have appeared in the author's previous papers on maximal subgroups of exceptional groups \cite{craven2017,craven2015un,craven2019un}, others are new. The following two sections go through the candidates for maximal subgroups of $F_4(q)$, and then $E_6(q)$ and ${}^2\!E_6(q)$, determining whether they do or do not belong to $\mc S$. Section \ref{sec:sl28} deals with whether $\SL_2(8)\rtimes 3$ lies in $E_6(q)$ or ${}^2\!E_6(q)$. After this we provide the tables of the other maximal subgroups and prove that they are correct, and the final section confirms that the subgroups in Tables \ref{tab:f4curlyS}, \ref{tab:e6curlyS} and \ref{tab:2e6curlyS} do not have overgroups in $F_4(q)$ and ${}^\varepsilon\!E_6(q)$, by checking all possible obstructions to maximality.

\begin{table}[ht]
\begin{center}
\begin{tabular}{lllll}
\hline Subgroup & $p$ & $q$& No. classes & Stabilizer
\\ \hline $\PSL_4(3).2_2$ & $2$ & $2$ &  $1$ & $\gen \gamma$
\\ ${}^3\!D_4(2).3$ & $p\neq 2$ & $p$ & $1$ & 1
\\ $\PSL_2(8).3$ & $7$ & $7$ & $1$ & 1
\\ $\PSL_2(8)$ & $p\equiv \pm 1\bmod 7$ & $p$ & $1$ & $1$
\\ $\PSL_2(8)$ & $p\neq 2,3$, $p\equiv\pm2,\pm3\bmod 7$ & $p^3$ & $1$ & $\gen\phi$
\\ $\PGL_2(13)$ & $7$ & $7$ & $1$ & 1
\\ $\PGL_2(13)$ & $p\neq 13$, $p\equiv \pm 1\bmod 7$ & $p$ & $3$ & $1$
\\ $\PGL_2(13)$ & $p\neq 2$, $p\equiv\pm2,\pm3\bmod 7$ & $p^3$ & $3$ & $1$
\\ $\PSL_2(17)$ & $p\neq 2$, $p\equiv\pm 1,\pm2,\pm4,\pm8\bmod 17$ & $p$ & $1$ & $1$
\\ $\PSL_2(17)$ & $p\neq 3$, $p\equiv\pm 3,\pm5,\pm6,\pm7\bmod 17$ & $p^2$ & $1$ & $\gen\phi$
\\ $\PSL_2(25).2$ & $p\neq 2,5$ & $p$ & $1$ & 1
\\ $\PSL_2(27).3$ & $7$ & $7$ & $1$ & $1$
\\ $\PSL_2(27)$ & $p\equiv \pm 1\bmod 7$ & $p$ & $1$ & $1$
\\ $\PSL_2(27)$ & $p\neq 3$, $p\equiv \pm 2,\pm 3\bmod 7$ & $p^3$ & $1$ & $\gen{\phi,\gamma}$
\\ \hline
\end{tabular}
\end{center}
\caption{Members of $\mc S$ for $G=F_4(q)$, $q$ a power of $p$. Here, $\phi$ is a generator of the group of field automorphisms, $\gamma$ is a non-trivial graph automorphism or $1$ if $p$ is odd. The group $\PSL_2(25).2$ is the extension by the field-diagonal automorphism, and is non-split.}
\label{tab:f4curlyS}
\end{table}

\begin{table}[ht]
\begin{center}
\begin{tabular}{lllll}
\hline Subgroup & $p$ & $q$& No. classes & Stabilizer
\\ \hline $M_{12}$ & $5$ & $5$ & $4$ & 1
\\ $J_3$ & $2$ & $4$ & $6$ & $\gen\phi$
\\ ${}^2\!F_4(2)'.2$ & $p\equiv 1\bmod 4$ & $p$ & $2e$ & $1$
\\ $\PSL_2(8)$ or $\PSL_2(8).3$ & $p\neq 2,3,7$, $p\equiv 1,2,4\bmod 7$ & $p$ & $2$ or $2e$ & $\gen\delta$ or $1$
\\ $\PSL_2(11)$ & $p\equiv \pm 1\bmod 5$, $p\equiv 1,3,4,5,9\bmod 11$ & $p$ & $2e$ & $\gen{\gamma}$
\\ $\PSL_2(11)$ & $p\equiv \pm 2\bmod 5$ & $p^2$ & $2e$ & $\gen{\gamma}$
\\ $\PSL_2(13)$ & $p\equiv 3,5,6\bmod 7$, $p\equiv \pm2,\pm5,\pm6\bmod 13$ & $p$ & $e$ & $\gen{\gamma}$
\\ $\PSL_2(19)$ & $5$ & $5$ & $1$ & $\gen\gamma$
\\ $\PSL_2(19)$ & $p\equiv \pm 1\bmod 5$, $p\equiv 1,4,5,6,7,9,11,16,17 \bmod 19$ & $p$ & $2e$ & $\gen\gamma$
\\ $\PSL_2(19)$ & $p\equiv\pm 2\bmod 5$ &$p^2$ & $2e$ & $\gen\gamma$
\\ $\PSL_2(13)$ & Nov., $p\equiv 1,2,4\bmod 7$, $p\equiv \pm1,\pm3,\pm4\bmod 13$ & $p$ & $e$ & $\gen\gamma$
\\ $\PSL_2(19)$ & Nov., $p=2$ & $4$ & $6$ & $\gen \gamma$
\\ \hline
\end{tabular}
\end{center}
\caption{Members of $\mc S$ for $G=E_6(q)$, $q$ a power of $p$. Here, $\delta$ is a generator for the group of diagonal automorphisms, $\phi$ is a generator of the group of field automorphisms, $\gamma$ is a non-trivial graph automorphism, and $e=\gcd(3,q-1)$. The subgroups $\PSL_2(13)$ and $\PSL_2(19)$ labelled `Nov.' are novelty maximals that only occur if an almost simple group induces $\gamma$ on $G$. For $\PSL_2(8)$, see Theorem \ref{thm:sl28}.}
\label{tab:e6curlyS}
\end{table}

\begin{table}[ht]
\begin{center}
\begin{tabular}{lllll}
\hline Subgroup & $p$ & $q$ & No. classes & Stabilizer
\\ \hline $Fi_{22}$ & $2$ & $2$ & $3$ & $\gen\phi$
\\ ${}^2\!F_4(2)'.2$ & $p\equiv 3\bmod 4$ & $p$ & $2e'$ & $1$
\\ $\PSL_2(8)$ or $\PSL_2(8).3$ & $p\neq 2,3,7$, $p\equiv 3,5,6\bmod 7$ & $p$ & $2$ or $2e'$ & $\gen\gamma$ or $1$
\\ $\PSL_2(11)$ & $p\equiv \pm 1\bmod 5$, $p\equiv 2,6,7,8,10\bmod 11$ & $p$ & $2e'$ & $\gen{\phi}$
\\ $\PSL_2(13)$ & $p\neq 2$, $p\equiv 1,2,4\bmod 7$, $p\equiv \pm2,\pm5,\pm6\bmod 13$ & $p$ & $e'$ & $\gen{\phi}$
\\ $\PSL_2(19)$ & $p\equiv \pm 1\bmod 5$, $p\equiv 2,3,8,10,12,13,14,15,18 \bmod 19$ & $p$ & $2e'$ & $\gen\phi$
\\ $\Omega_7(3)$ & Nov., $2$ & $2$ & $3$ & $\gen\phi$
\\$\PSL_2(13)$ & Nov., $p\equiv 2,5,6\bmod 7$, $p\equiv \pm1,\pm3,\pm4\bmod 13$ & $p$ & $e'$ & $\gen\phi$

\\ \hline
\end{tabular}
\end{center}
\caption{Members of $\mc S$ for $G={}^2\!E_6(q)$, $q$ a power of $p$. Here, $\delta$ is a generator for the group of diagonal automorphisms and $\phi$ is a generator of the group of field automorphisms, and $e'=\gcd(3,q+1)$. The subgroups $\Omega_7(3)$ and $\PSL_2(13)$ labelled `Nov.' are novelty maximals that only occur if an almost simple group induces $\phi$ on $G$. For $\PSL_2(8)$, see Theorem \ref{thm:sl28}.}
\label{tab:2e6curlyS}
\end{table}

\section{Notation and preliminaries}
\label{sec:prelim}
In this section we review the results that we need for the proofs that follow. Much of these will be familiar from \cite{craven2015un,craven2019un}, with a few important new techniques for attacking `non-generic' subgroups, i.e., simple subgroups that are not of Lie type in the same characteristic as the algebraic group.

We start with some notation from module theory. Let $k$ be an algebraically closed field of characteristic $p\geq 0$. If $H$ is a finite group and $M$ is a $kH$-module then $M^H$ denotes the set of $H$-fixed points of $M$. If $L\leq H$ then $M{\downarrow_L}$ is the restriction of $M$ to $L$. We often need to describe the socle structure of a $kH$-module. We use `$/$' to separate socle layers, so that
\[ A/B,C/D,E\]
denotes a module with socle $D\oplus E$, second socle $B\oplus C$, and third socle $A$. Let $M^*$ denote the dual of $M$. Write $S^i(M)$ and $\Lambda^i(M)$ for the $i$th symmetric and exterior power of $M$ respectively. If $M$ is simple, write $P(M)$ for the projective cover of $M$.

We often use the dimensions of simple $kH$-modules to label them, rather than give them special names. For example, if $H$ has simple modules of dimensions $1$, $5$, $5$, $5$, $10$ and $10$, with only the $1$- and one $5$-dimensional module self-dual, then we would label these modules $1$, $5_1$, $5_2$, $5_2^*$, $10$ and $10^*$.

If $p>0$ and $u\in \GL_n(k)$ has $p$-power order then $u$ is a unipotent element, and acts (up to conjugacy) as a sum of Jordan blocks of various sizes. We write $n_1^{a_1},\dots,n_r^{a_r}$ to describe the sizes, for example $5^2,4,1$ for an element of $\GL_{15}(k)$, as in \cite{lawther1995}.

Let $\mbG$ denote a simple, simply connected algebraic group of type $F_4$, $E_6$, $E_7$ or $E_8$ over an algebraically closed field of characteristic $p$, and let $\sigma$ be a Frobenius endomorphism of $\mbG$, writing $G_{\mathrm{sc}}$ for $\mbG^\sigma$. Thus $G_{\mathrm{sc}}$ is one of the groups $F_4(q)$, $E_6(q)_{\mathrm{sc}}$, ${}^2\!E_6(q)_{\mathrm{sc}}$, $E_7(q)_{\mathrm{sc}}$ and $E_8(q)$ for some $q$ a power of $p$ (so that, for example $E_6(q)_{\mathrm{sc}}$ has a centre of order $\gcd(q-1,3)$). Write $E_6(q)$, and so on, for the simple quotient of $E_6(q)_{\mathrm{sc}}$, and $G=G_{\mathrm{sc}}/Z(G_{\mathrm{sc}})$.  Write $F_q$ for the `standard' Frobenius morphism, which induces the trivial action on the Weyl group. Thus $F_4(q)=(F_4)^{F_q}$, for example. Let $\mb T$ denote a maximal torus of $\mbG$, and assume it is $\sigma$-stable if $\sigma$ is being considered. Let $\bar G$ be an almost simple group with socle $G$, and $G_{\mathrm{ad}}$ be the adjoint version of $G$, often denoted by simple group theorists by $\mathrm{Inndiag}(G)$. The finite groups $G$ and $\bar G$ are the groups whose maximal subgroups we will investigate in this article.

As in \cite{craven2015un,craven2019un} (and possibly originally in \cite{liebeckmartinshalev2005}), write $\Aut^+(\mbG)$ for the group generated by the inner automorphisms, the graph automorphisms for $\mb G$ of type $E_6$ (and the graph morphisms for $F_4$ and $p=2$), and the $p$-power field morphisms of $\mbG$. Thus $\Aut^+(\mbG)$ induces the full group $\Aut(G)$ if $G_{\mathrm{sc}}=\mbG^\sigma$.

As in \cite{craven2015un,craven2019un}, a subgroup $H$ of $\mbG$ is \emph{Lie imprimitive} if $H$ is contained in a proper, positive-dimensional subgroup $\mb X$ of $\mbG$, and \emph{strongly imprimitive} if the subgroup $\mb X$ may be chosen such that it is $N_{\Aut^+(\mbG)}(H)$-stable. If $H$ is strongly imprimitive then $N_{\bar G}(\bar H)$ is contained in $N_{\bar G}(\bar X)$ for $X=\mb X^\sigma$, where $\bar H$ and $\bar X$ are the images of $H$ and $X$ in $G$ respectively.\footnote{If $Z(\mbG)\neq 1$ then $G$ is not a subgroup of $\mbG$, so $H$ is not a subgroup of $G$, and instead we take the quotient $H\cdot Z(\mbG)/Z(\mbG)$.}

We use the `standard' labelling for highest-weight modules, as in the previous entries \cite{craven2015un,craven2019un} in this series, coming from \cite{bourbakilie2}. Write $W(\lambda)$ and $L(\lambda)$ for the Weyl and simple modules with highest weight $\lambda$. We define two modules $M(\mbG)$ and $L(\mbG)$ for $\mbG$ (minimal and Lie algebra), as the Weyl module with the labels in the table below. Write $M(\mbG)^\circ$ and $L(\mbG)^\circ$ for the quotient of $M(\mbG)$ and $L(\mbG)$ by their fixed points, which always has codimension at most $1$. The dimensions of $M(\mbG)^\circ$ and $L(\mbG)^\circ$ for the groups involved are in Table \ref{t:labels}. The modules $M(\mbG)^\circ$ and $L(\mbG)^\circ$ are simple unless $\mbG=F_4$ and $p=2$, when $L(\mbG)$ is a non-split extension of $M(\mbG)$ and $M(\mbG)^\tau$, where $\tau$ is a non-trivial graph morphism of $\mbG$.\footnote{The map $\tau$ may be taken to send $x_r(t)$ to $x_{\rho(r)}(t^{\alpha(r)})$, where $\rho$ is a symmetry of the root system and $\alpha(r)=1$ of $r$ is long and $\alpha(r)=p$ if $r$ is short.} In this case we consider $L(\mbG)$ itself.

\begin{table}[ht]
\begin{center}\begin{tabular}{lllll}
\hline $\mbG$ & Label for $M(\mbG)$ & Label for $L(\mbG)$ & $\dim(M(\mbG)^\circ)$ & $\dim(L(\mbG)^\circ)$
\\ \hline $F_4$ & $W(\lambda_4)$ & $W(\lambda_1)$ & $26-\delta_{p,3}$ & $52$
\\ $E_6$ & $L(\lambda_1)$ or $L(\lambda_6)$ & $W(\lambda_2)$ & $27$ & $78-\delta_{p,3}$
\\ $E_7$ & $L(\lambda_7)$ & $W(\lambda_1)$ &$56$ & $133-\delta_{p,2}$
\\ $E_8$ & $-$ & $L(\lambda_1)$ &$-$ & $248$
\\ \hline
\end{tabular}\end{center}
\caption{Labels and dimensions for minimal and adjoint modules for exceptional groups.}
\label{t:labels}
\end{table}

\medskip

A maximal subgroup of $\bar G$ not containing $G$ itself is either $N_{\bar G}(X)$ for a maximal subgroup $X$ of $G$, called an \emph{ordinary} maximal subgroup, or $N_{\bar G}(X)$ for $X$ a non-maximal subgroup of $G$, called a \emph{novelty} maximal subgroup. The maximal subgroups of $\bar G$ are known except for almost simple subgroups by \cite{borovik1989,liebeckseitz1990}. The quasisimple subgroups for exceptional $\mbG$ are all known, in the sense that there is a list of groups, and a quasisimple group embeds in $\mbG$ if and only if it appears on that list \cite{liebeckseitz1999}. What remains for the classification of maximal subgroups of finite simple exceptional groups is to determine, given $H$ on that list, whether $H$ embeds in $G_{\mathrm{sc}}$ (rather than just $\mbG$), if $N_{\bar G}(H\cdot Z(\mbG)/Z(\mbG))$ can be maximal, and if so to count the number of conjugacy classes of such subgroups. For notational convenience, as did earlier we  write $\bar H$ for the image of $H$, a subgroup of $\mbG$ in $\bar G$, so the subgroup $H\cdot Z(\mbG)/Z(\mbG)$. So $H$ is a subgroup of $\mbG$, and $\bar H$ is a subgroup of $G$ and $\bar G$.

Specifically, let $\ms X$ denote the set of maximal positive-dimensional subgroups of $\mbG$, and let $\ms X^\sigma$ denote the set of fixed points $\mb X^\sigma$ for $\mb X$ a $\sigma$-stable member of $\ms X$. If $Z(\mbG^\sigma)\neq 1$, then extend the notation $\ms X^\sigma$ to include the images of elements of $\ms X^\sigma$ under the quotient map $\mbG^\sigma\to G$. We extend further by defining $\ms X^\sigma$ for $\bar G$ to be $N_{\bar G}(X)$ for $X$ lying in the set $\ms X^\sigma$ associated to $G$. Thus $\ms X^\sigma$ in all cases is the set of maximal subgroups arising from positive-dimensional subgroups of $\mbG$. By \cite[Theorem 2]{liebeckseitz1990} (see also \cite{borovik1989}), a maximal subgroup of $\bar G$ (where $G_{\mathrm{sc}}=\mbG^\sigma$) is one of:
\begin{enumerate}
\item $X$, where $X\in \ms X^\sigma$;
\item a subgroup of the the same type as $\mbG$, e.g., $E_6(p)$ and ${}^2\!E_6(p)$ in $E_6(p^2)$;
\item the Borovik subgroup $(\Alt(5)\times \Alt(6)).2^2$, which only occurs for $\mbG=E_8$;
\item an exotic $r$-local subgroup, given in \cite{clss1992};
\item an almost simple group $N_{\bar G}(\bar H)$ for $\bar H$ a simple subgroup of $G$ not contained in one of the above lists, the collection of which is denoted $\mc S$.
\end{enumerate}
Furthermore, it follows easily from the definition that if $H$ is strongly imprimitive, then $H$ is contained in one of the subgroups from (i). All of the subgroups except those in (v) are known, so here we determine those subgroups in (v) that do not also lie in (i) and (ii).

The list of potential quasisimple subgroups $H$ naturally divides into two classes: Lie type in characteristic $p$ (including the cases where the simple quotient is Lie type but $H$ need not be, e.g., $3\cdot \Omega_7(3)$ for $p=3$), called \emph{generic} subgroups, and all other quasisimple groups, called \emph{non-generic} subgroups. Recently, significant reductions in the list of potential maximal subgroups $N_{\bar G}(\bar H)$ of $\bar G$ have been made by Alastair Litterick \cite{litterickmemoir} and the author \cite{craven2017} for non-generic $H$.

First, from the author's previous papers \cite{craven2015un,craven2019un} (relying heavily on earlier work of others, for example \cite{lst1996,liebeckseitz1998}), we see that if $H$ is in (v) then $H$ is not a generic subgroup (this is for $F_4$ and $E_6$, the result is not yet completely known for $E_7$ and $E_8$). Thus, for $\mbG$ of types $F_4$ and $E_6$, $H$ is one of an alternating group, a sporadic group, and a group of Lie type in characteristic $r\neq p$. In \cite{craven2017} severe restrictions on alternating groups were proved, and in \cite{litterickmemoir} there are more results about non-generic subgroups, and in particular a complete list of all possible sets of composition factors for $M(\mbG){\downarrow_H}$ and $L(\mbG){\downarrow_H}$ was enumerated. A \emph{conspicuous} set of composition factors for $H$ is a pair $(M,L)$ of multisets of irreducible $kH$-modules (or for $E_8$, a single multiset $L$) where the eigenvalues of all $p'$-elements of $H$ on $M$ and $L$ appear as eigenvalues of semisimple classes of $\mbG$ on $M(\mbG)$ and $L(\mbG)$ respectively. (Notice that the pair of multisets offers more information; a conjugacy class of $\mbG$ must have the same eigenvalues as an element of $H$ on both $M(\mbG)$ and $L(\mbG)$ simultaneously.)

\medskip

We need a few results from the theory of algebraic groups. The first allows us to find certain subgroups, usually the Borel subgroup of a group $\PSL_2(r)$ for $p\nmid r$.

\begin{lem}\label{lem:normtorus} Let $L$ be a finite subgroup of $\mbG$.
\begin{enumerate}
\item If $L$ is elementary abelian of order $r^2$ for some $r\neq p$ then $L$ lies inside a maximal torus $\mb T$.
\item If $L$ is a supersoluble $p'$-group then $L$ lies inside $N_\mbG(\mb T)$ for $\mb T$ a maximal torus. Furthermore, if $L\leq \mb G^\sigma$ then $\mb T$ may be chosen to be $\sigma$-stable.
\end{enumerate}
\end{lem}

The second of these is a famous result of Borel and Serre \cite{borelserre1953} (see \cite[Theorem II.5.16, p.210]{springersteinberg1970} for a proof in English, with the added conditions on $\sigma$-stability). A reference for (i) is \cite[II.5.1, p. 206]{springersteinberg1970}.

The next lemma will be used occasionally, with $X$ a $\sigma$-stable $\mbG$-conjugacy class of subgroups.

\begin{lem}[See, for example, {{\cite[Theorem 21.11]{malletesterman}}}]\label{lem:fixedpoints} Let $\mb Y$ be a connected linear algebraic group with a Steinberg endomorphism $\rho$, and write $Y=\mb Y^\rho$. Let $X$ be a non-empty set such that $\mb Y$ and $\rho$ act on $X$, and such that $(x^g)^\rho=(x^\rho)^{(g^\rho)}$ for all $x\in X$ and $g\in \mb Y$. Then $X^\rho\neq \emptyset$, i.e., there are $\rho$-fixed points in $X$.

Furthermore, if for some $x\in X$ the stabilizer $\mb Y_x$ of $x$ is closed, and $\mb Y$ is transitive on $X$, then the $Y$-orbits on the $\rho$-fixed points of $X$ are in one-to-one correspondence with $\rho$-conjugacy classes on $\mb Y_x/\mb Y_x^\circ$. In particular, if $\mb Y_x$ is connected then all $\rho$-fixed points of $X$ are $Y$-conjugate.
\end{lem}

The lemma can be used with centralizers, not just of elements---as in \cite[Theorem 26.7]{malletesterman}, for example---but also of subgroups. If $X$ just consists of all $\mbG$-conjugates of a subgroup $L$ then the stabilizer is $N_\mbG(L)$. So instead we let $X$ consist of all $\mbG$-conjugates of a generating set for $L$, and then $\mbG_x$ is $C_\mbG(L)$. If this is connected (usually in our cases $C_\mbG(L)$ is either $Z(\mbG)$ or a torus, depending on $L$) then all copies of $L$ inside $G_{\mathrm{sc}}$ are $G_{\mathrm{sc}}$-conjugate, and if the centralizer is $Z(\mb G)$ then the centralizer is connected in $\mbG/Z(\mbG)$.

\begin{cor}\label{cor:centsubgroups}
Let $\mb Y$ be a connected linear algebraic group with a Steinberg endomorphism $\rho$, and write $Y=\mb Y^\rho$. Let $L$ be a subgroup of $Y$. All $\mb Y$-conjugates of $L$ contained in $Y$ are $Y$-conjugate if $C_{\mb Y}(L)$ is connected, and if $C_{\mb Y}(L)\leq Y$ then the $\mb Y$-conjugates of $L$ contained in $Y$ fall into a divisor of $n$ distinct $Y$-classes, where $n$ is the number of conjugacy classes of $C_{\mb Y}(L)$.
\end{cor}
\begin{proof} Let $l_1,\dots,l_d$ be a generating set for $L$, and let $X$ denote the set
\[ \{(l_1^g,l_2^g,\dots,l_d^g)\mid g\in \mb Y\}.\]
Clearly $\mb Y$ acts transitively on this, and also $X$ is $\rho$-stable, since $l_i^\rho=l_i$, so
\[ (l_1^g,l_2^g,\dots,l_d^g)^\rho=(l_1^{g^\rho},l_2^{g^\rho},\dots,l_d^{g^\rho})\in X.\]
Thus by Lemma \ref{lem:fixedpoints}, the $Y$-classes on $X^\rho$ are in bijection with the $\rho$-classes of $A=C_{\mb Y}(L)/C_{\mb Y}^\circ(L)$. If $C_{\mb Y}(L)$ is connected then this means that there is exactly one class. But if $L_1$ is $\mb Y$-conjugate to $L$ and contained in $\mb Y^\rho$, then $L_1$ is generated by some tuple in $X$, and this tuple is also in $X^\rho$. Hence $L$ and $L_1$ are $Y$-conjugate.

Now suppose that $C_{\mb Y}(L)\leq Y$, so that $A=C_{\mb Y}(L)$. Notice that $g^\rho=g$ for $g\in A$, so the $\rho$-classes on $A$ are exactly the usual conjugacy classes on $A$. We finally note that two elements $x$ and $y$ of $X$ generate the same subgroup $L_1$ if and only if there exists an element $g\in N_{\mb Y}(L_1)$ mapping $x$ to $y$; say there are $m$ of these. This condition is conjugation-equivariant, so the number of $Y$-classes of subgroups is $n/m$. The result is proved. 
\end{proof}

If $H$ is a $\sigma$-stable quasisimple subgroup of an exceptional group $\mbG$ and $N_{\bar G}(\bar H)$ is an almost simple maximal subgroup of $\bar G$ then $C_{\bar G}(\bar H)=1$. We would like to deduce that $C_\mbG(H)=Z(\mbG)$, but \emph{a priori} this is not clear, and it certainly is not true in general that if $\bar H$ is an arbitrary subgroup of $G$ with trivial centralizer then $C_\mbG(H)=Z(\mbG)$, for example when $q=2$ and the centralizer is a torus. If $X$ is a group and $Z$ is a central subgroup of $X$, let $Y$ be a subgroup of $X$ with image $\bar Y$ in $X/Z$. Of course, $N_X(Y)/Z=N_{X/Z}(\bar Y)$, and $C_X(Y)/Z\leq C_{X/Z}(\bar Y)$. If $|Z|$ is finite and $\bar Y$ is generated by elements of order prime to $|Z|$, then $C_X(Y)/Z=C_{X/Z}(\bar Y)$. If $\bar Y$ is simple, the derived subgroup $Y'$ of $Y$ is a central extension of $\bar Y$, and $Y$ is a central product of $Y'$ and $Z$. If $\bar Y$ is simple and $|\bar Y|$ is divisible by a prime not dividing $|Z|$, we also have $C_X(Y')=C_X(Y)$.

Applied to our case where $H$ is quasisimple with $Z(H)\leq Z(\mbG)$, $C_\mbG(H)/Z(\mbG)=C_{\mbG/Z(\mbG)}(\bar H)$, so we must consider the situation where $C_\mbG(H)^\sigma=Z(\mbG)$ while $C_\mbG(H)>Z(\mbG)$.

So let $H$ be a quasisimple subgroup of $\mbG$ with $C_\mbG(H)^\sigma=Z(\mbG)$ and consider $C_\mbG(H)$, which is $\sigma$-stable. If this is positive-dimensional then $H$ is contained in a proper, positive-dimensional subgroup of $\mbG$, namely $H\cdot C_\mbG(H)$, and so therefore $H$ is contained in a member of $\ms X^\sigma$ (and the same for $\bar H$). On the other hand, if $C_\mbG(H)$ is finite then $H\leq C_\mbG(C_\mbG(H))$. If $\mbG$ does not have type $E_8$ then this double centralizer is always positive-dimensional by \cite[Lemma 1.3(a)]{liebeckseitz1990}, so again $H$ is contained in a member of $\ms X^\sigma$.

Since we are considering only $F_4$ and $E_6$ here, we do not need to analyse the situation further to state that if $\bar H$ is a simple subgroup of $G$ with $C_G(\bar H)=1$, then either $C_\mbG(H)=Z(\mbG)$ or $N_{\bar G}(\bar H)$ is contained in a member of $\ms X^\sigma$. In particular if $N_{\bar G}(\bar H)$ lies in $\mc S$ then $C_\mbG(H)=Z(\mbG)$, so we may use Corollary \ref{cor:centsubgroups}. Taking $\mb Y=\mbG/Z(\mbG)$, this shows that all $\mbG/Z(\mbG)$-conjugates of a simple subgroup $\bar H$ of $G$ in $\mc S$ are $G_{\mathrm{ad}}$-conjugate. The same holds for the subgroups $H$ of $\mbG$, where the $G_{\mathrm{sc}}$-classes of conjugates in a single $\mbG$-class are fused by $\mathrm{Outdiag}(G)$.

One of our tasks in computing the maximal subgroups of $\bar G$ is to determine the number of $G$-classes and the stabilizer in $\Out(G_{\mathrm{sc}})$ of a quasisimple group $H$ in $\mc S$. The outer diagonal part of this stabilizer very easy in almost all situations. This is combined with the fact that any two subgroups in $\mc S$ that are $\mbG$-conjugate are conjugate in $G_{\mathrm{ad}}$.

\begin{cor}\label{cor:diagaut} Let $\mbG$ have type $E_6$, and let $H$ be a $\sigma$-stable quasisimple subgroup of $\mbG$. If $|N_\mbG(H)/HC_\mbG(H)|$ is not divisible by $3$ then the stabilizer of $H$ in $\Out(G_{\mathrm{sc}})$ contains no non-trivial outer diagonal automorphism.
\end{cor}
\begin{proof} This is clear. If any elements of $\Outdiag(G)$, which has order $1$ or $3$, normalize $H$ they must correspond to elements of $N_\mbG(H)$. If $3$ does not divide $|N_\mbG(H)/HC_\mbG(H)|$ then any such element must be trivial.
\end{proof}

Examining Tables \ref{tab:e6curlyS} and \ref{tab:2e6curlyS}, we see that the only option for $H$ where $|\Out(H)|$ has order divisible by $3$ is when $H\cong \PSL_2(8)$, and then $N_\mbG(H)$ does include the outer automorphism of order $3$. This case is much more delicate than the others and we devote Section \ref{sec:sl28} to it. (This situation is considerably more complex in \cite{craven2021un} for $\mbG$ of type $E_7$ and $p$ odd, because there $|\Outdiag(G)|$ has order $2$ and there are several subgroups $H$ with $N_\mbG(H)/HC_\mbG(H)$ of order $2$.)

\medskip

Next, we describe some line stabilizers for $M(F_4)^\circ$ and $M(E_6)$.

\begin{lem}\label{lem:linestabs} Let $\mbG$ be one of $F_4$ and $E_6$. Suppose that a subgroup $H$ of $\mbG$ stabilizes a line on $M(\mbG)^\circ$.
\begin{enumerate}
\item If $\mbG$ is of type $F_4$ then $H$ lies in a maximal parabolic subgroup, in a subgroup of type $B_4$, or in a subgroup $D_4.3$ ($D_4.\Sym(3)$ for $p=3$). If $p\neq 3$ then $B_4$ and $D_4.3$ act on $M(F_4)$ as
\[ 1\oplus 9\oplus 16\quad \text{and} \quad 1_a\oplus 1_a^*\oplus 24\]
respectively, where $1_a$ is a non-trivial $1$-dimensional representation, and $24$ is the sum of the three $8$-dimensional summands for the action of $D_4$. If $p=3$ then $B_4$ and $D_4.\Sym(3)$ act on $M(F_4)^\circ$ as
\[ 9\oplus 16\quad\text{and}\quad 1_b\oplus 24\]
respectively, where $1_b$ is the non-trivial $1$-dimensional representation and $24$ is as before.

In particular, if $H$ centralizes a line on $M(F_4)^\circ$ and $p\neq 3$, or has no subgroup of index $2$ or $3$, then $H$ is contained in a maximal parabolic subgroup or $B_4$.

\item If $\mbG$ is of type $E_6$ then $H$ lies in a subgroup of type $F_4$, a $D_5T_1$-parabolic subgroup, or a subgroup $\mb U\cdot B_4T_1$, where $\mb U$ is unipotent of dimension $16$.\footnote{Only one of the two classes of $D_5T_1$-parabolic subgroups stabilizes a line on $M(E_6)$ (with the other stabilizing a hyperplane). The subgroup $\mb U\cdot B_4T_1$ that stabilizes a line does not lie in this class, but in the other one.}
If $p\neq 3$ then the subgroup $F_4$ acts on $M(E_6)$ as
\[ 1\oplus M(F_4),\]
and if $p=3$ then it acts as $1/M(F_4)^\circ/1$. The $D_5$-parabolic subgroup and the $B_4$-type subgroup act as
\[ 10/16/1\qquad\text{and}\qquad 1/16/1,9\]
respectively (with the latter structure only valid for $p$ odd).
\end{enumerate} 
\end{lem}
\begin{proof} The second part follows immediately from \cite[Lemma 5.4]{liebecksaxl1987}. For the first, by \cite[(B.1)]{cohencooperstein1988} the statement about the line stabilizers holds. The action of $B_4$ is well known, for example \cite[Theorem 3.1]{thomas2016}, and from this the action of $D_4$ is clear. To obtain the action of $D_4.3$ from this for $p\neq 3$, we note that $\mb X=D_4.\Sym(3)$ does not stabilize a line on $M(F_4)$, so must act as the $2$-dimensional simple module on the fixed-point space $M(F_4)^{D_4}$. (The simple module of dimension $24$ is obviously a summand, since the triality graph automorphism acts transitively on the three $8$-dimensional factors for $D_4$.)

When $p=3$, $\mb X$ stabilizes the line, and must invert it because $\mb X$ still only stabilizes a single line on $M(F_4)$, whence it has both $1$-dimensional composition factors on $M(F_4)$ (as there is a non-split extension between the two simple modules for $\Sym(3)$, but no non-split self-extension). Since $F_4$ itself has a trivial composition factor on $M(F_4)$, $\mb X$ must act non-trivially on the $1$-dimensional factor in $M(F_4)^\circ$.
%
%
%
%
%
\end{proof}

Finally, we have a small lemma on some polynomials and their splitting fields. These polynomials will appear as the minimal polynomials for the irrationalities in the characters that we examine.

\begin{lem}\label{lem:irrationalities} Let:
\begin{enumerate}
\item $f_1(x)=x^2-x-4$;
\item $f_2(x)=x^2+x+3$;
\item $f_3(x)=x^2+x-1$;
\item $f_4(x)=x^2+x+5$;
\item $f_5(x)=x^3-x^2-2x+1$.
\end{enumerate}
The splitting fields for the first four polynomials over $\F_p$ are $\F_p[\sqrt d]$ for $d=17,-11,5,-19$. The element $\sqrt d$ belongs to $\F_p$ if and only if:
\begin{enumerate}
\item $p\equiv 0,\pm1, \pm2,\pm4,\pm8\bmod 17$,
\item $p\equiv 0,1,3,4,5,9\bmod 11$;
\item $p\equiv 0,\pm 1\bmod 5$;
\item $p\equiv 0,1,4,5,6,7,9,11,16,17\bmod 19$
\end{enumerate}
respectively.

In the final case, $f_5(x)$ splits over $\F_q$ if and only if $q\equiv 0,\pm 1\bmod 7$.\end{lem}
\begin{proof} We just give solutions to the equations in the respective field, assuming $p\neq 2$:
\begin{enumerate}
\item $(\sqrt{17}+1)/2$;
\item $(\sqrt{-11}-1)/2$;
\item $(\sqrt 5-1)/2$;
\item $(\sqrt{-19}-1)/2$;
\item $-(\zeta_7+\zeta_7^{-1})$ (where $\zeta_7$ is a primitive $7$th root of unity).
\end{enumerate}
(If $p=2$ then the first case is $x(x-1)$, and the next three are all $x^2+x+1$, which is irreducible over $\F_2$. Thus the result holds for $p=2$ as well.) In the first four cases, if $d=0$ in $\F_p$ then the polynomial splits, so we assume this is not the case. The congruences for the primes follows from quadratic reciprocity in the first four cases.

In the final case, note that the solution certainly lies in $\F_q$ if $q\equiv 1\bmod 7$. For other congruences, embed $\F_q$ in $\F_{q^a}$ where $q^a\equiv 1\bmod 7$ and note that $\zeta_7+\zeta_7^{-1}$ remains invariant under the map $x\mapsto x^q$ if and only if $q\equiv\pm1\bmod 7$.
\end{proof}

\section{Techniques used in the proof}
\label{sec:techniques}
This section surveys the various ideas that we will use to understand subgroups of exceptional algebraic groups and their finite versions. The first subsection discusses strong imprimitivity, a sufficient condition for $N_{\bar G}(\bar H)$ to never be a maximal subgroup of $\bar G$. We then talk a little about how to go from $\mbG$-conjugacy to $G$-conjugacy. After that, we discuss the trilinear form, as used in \cite{craven2019un}. The subsection after discusses novelty maximal subgroups, and finally we consider the Lie algebra structure of $L(\mbG)$ and subalgebras, as used in \cite{craven2015un} to eliminate certain subgroups $\PSL_2(q)$ in defining characteristic.

\subsection{Strong imprimitivity}

In \cite{craven2015un,craven2019un}, we developed several techniques to prove that a subgroup $H$ of $\mbG$ is strongly imprimitive. The easiest of these is the following (see \cite[Propositions 4.5 and 4.6]{craven2015un}).

\begin{prop}\label{prop:fix1space} Suppose that $H$ is a finite subgroup of $\mbG$ possessing no subgroup of index $2$. If $H$ centralizes a line or hyperplane on $M(\mbG)^\circ$ or $L(\mbG)^\circ$ then $H$ is strongly imprimitive.
\end{prop}

The best way to prove that $H$ centralizes a line on a module is to use the concept of pressure, which is by now well established. If $M$ is a $kH$-module then the \emph{pressure} of $M$ is
\[ \sum_{V\in \cf(M)} \dim(H^1(H,V))-\delta_{V,k},\]
where $\cf(M)$ is the multiset of composition factors of $M$. This means that we add up the $1$-cohomologies of all composition factors and subtract the number of trivial factors.

\begin{prop}\label{prop:pressure} Let $H$ be a finite group such that $O^p(H)=H$, and let $M$ be a $kH$-module. If $M$ has negative pressure then $H$ stabilizes a line on $M$. If $\dim(H^1(H,V))=\dim(H^1(H,V^*))$ for all simple $kH$-modules $V$, $M$ has at least one trivial composition factor, and $M$ has pressure $0$ then $H$ stabilizes either a line or a hyperplane on $M$.
\end{prop}

This appears in for example \cite[Lemma 1.8]{craven2017}, \cite[Lemma 2.2]{craven2015un} or \cite[Proposition 3.6]{litterickmemoir}, and originally \cite[Lemma 1.2]{lst1996}. In \cite[Section 4]{craven2015un} a number of other conditions were produced that guarantee strong imprimitivity. For example, this is \cite[Proposition 4.2]{craven2015un}.

\begin{prop} Let $H$ be a finite subgroup of $\mbG$. Either $H$ is strongly imprimitive or $H$ is $\mbG$-irreducible, i.e., $H$ is not contained in a maximal parabolic subgroup of $\mbG$.
\end{prop}

We also see in that section the following. A module $V$ is \emph{graph-stable} if it is semisimple, has at most two composition factors, and for $\tau$ a (possibly trivial) graph morphism in $\Aut^+(\mbG)$, the composition factors of $V^\tau$ are the same as those of $V$, up to Frobenius twist. Instructive examples are any semisimple module for $E_7$, $M(E_6)\oplus M(E_6)^*$, and $L(\lambda_1)\oplus L(\lambda_4)$ for $F_4$ and $p=2$ (and odd primes as well, of course), but not $L(\lambda_2)\oplus L(\lambda_4)$ for $p=2$.

(The example $L(\lambda_1)\oplus L(\lambda_4)$ is why we need to add the condition `up to Frobenius twist' in the definition. Since `the' graph morphism for $F_4$ squares to a Frobenius endomorphism, in general one cannot remove this condition.) Using the identification of a module with its Frobenius twist via a suitable semilinear map, we can let $\tau$ permute the subspaces of $V$.

For the purposes here, we restrict our graph-stable modules to have at most two composition factors (if we were dealing with $D_4$ we would need three). This is just to simplify the situation, and the only issue with extending the definition is keeping track of the maps and the weights.

\begin{prop}\label{prop:intorbit} Let $V$ be a graph-stable module for $\mbG$, and let $H$ be a subgroup of $\mbG$. If there exists an $N_{\Aut^+(\mbG)}(H)$-orbit $\mc W$ of subspaces of $V$ that is stabilized by both $H$ and a positive-dimensional subgroup $\mb X$ of $\mbG$, but not by $\mbG$ itself, then $H$ is strongly imprimitive.
\end{prop}

A result of this type started life in \cite[Proposition 1.12]{liebeckseitz1998}, then was generalized in \cite[Section 4.2]{litterickmemoir}. The results of that section were distilled into \cite[Proposition 4.3]{craven2015un}, and this is a special case of that proposition, but this is entirely the work of \cite{liebeckseitz1998,litterickmemoir}. A corollary of this, used frequently in the author's work, is as follows. A subgroup $H$ is a \emph{blueprint} for a graph-stable $k\mbG$-module $V$ if there exists a proper, positive-dimensional subgroup $\mb X$ of $\mbG$, such that $H$ and $\mb X$ stabilize the same subspaces of $V$.

\begin{cor}\label{cor:blueprint} If $H$ is a blueprint for $V$, then either $H$ and $\mbG$ stabilize the same subspaces of $V$ or $H$ is strongly imprimitive.
\end{cor}

Note the fact that if $H$ is a blueprint for $V_1$ and $V_2$, then $H$ need not be a blueprint for $V_1\oplus V_2$. This is important when considering groups with a non-trivial graph morphism.

\medskip

If $p=3$ and $\mbG=E_6$, then $\mb X=F_4$ acts on $M(E_6)$ uniserially, with top and socle the trivial module and heart $M(F_4)^\circ$, as opposed to $L(0)\oplus M(F_4)$ for other primes (see Lemma \ref{lem:linestabs}). The next lemma shows that if $H\leq \mb X$ splits this extension then $H$ is strongly imprimitive.

\begin{lem}\label{lem:nocohomologyonminimal} Let $\mbG=E_6$ and $p=3$, and let $F_4=\mb X\leq \mbG$. If $H$ is a finite subgroup of $\mb X$ with no subgroups of index $2$, and the action of $H$ on $M(E_6)$ is the sum of two trivial modules and $M(F_4)^\circ{\downarrow_H}$, then $H$ is strongly imprimitive in $\mb X$.
\end{lem}
\begin{proof} If $M(E_6)^H$ has dimension at least $3$ then $(M(F_4)^\circ)^H\neq 0$, so $H$ is strongly imprimitive in $\mb X$ by Proposition \ref{prop:fix1space}. Thus we may assume that $M(E_6)^H$ has dimension exactly $2$. The centralizer $\mb Y$ of this subspace is positive-dimensional by dimension counting ($\dim(\mbG)=78$ and $\dim(M(E_6))=27$) so $H$ is contained in a positive-dimensional subgroup of $\mbG$. Furthermore, any inner, diagonal and field morphism in $N_{\Aut^+(\mbG)}(H)$ stabilizes $M(E_6)^H$, so stabilizes $\mb Y$ (but we can say nothing about graph morphisms in $\Aut^+(\mb G)$).

Finally, $\mb Y$ must be contained in $\mb X$, whence $H$ is strongly imprimitive in $F_4$. (Note that every element of $\Aut^+(\mb X)$ extends to an element of $\Aut^+(\mbG)$ not involving a graph automorphism.)
\end{proof}

\subsection{Actions on \texorpdfstring{$M(\mbG)$ and $L(\mbG)$}{M(G) and L(G)}}
\label{sec:actions}

Let $H$ be a finite quasisimple subgroup of $\mbG$ with $Z(H)\leq Z(\mbG)$, and suppose that $H$ is not of Lie type in characteristic $p$. The possible sets of composition factors for $M(\mbG){\downarrow_H}$ and $L(\mbG){\downarrow_H}$ are collated in tables in \cite[Chapter 6]{litterickmemoir}. As pressure is also used in \cite{litterickmemoir}, the first test is whether the pressure of either $M(\mbG)^\circ$ or $L(\mbG)^\circ$ is negative (or possibly pressure $0$ if the other conditions are met). If this is the case, then $H$ is strongly imprimitive by Propositions \ref{prop:fix1space} and \ref{prop:pressure}, so we may immediately exclude this. This was done in \cite{litterickmemoir}, and in the tables in that paper, exactly those rows labelled `\textbf{P}' consist of possible embeddings $H$ where both $M(\mbG)^\circ{\downarrow_H}$ and $L(\mbG)^\circ{\downarrow_H}$ have non-negative pressure.

Thus we may always assume that we have one of the sets of composition factors from those tables when analysing the potential Lie primitive, and then not strongly imprimitive subgroups. Despite having positive pressure, many of these sets of composition factors end up yielding subgroups that always stabilize lines on one of the two modules, and some sets of factors do not correspond to embeddings of subgroups at all. Even if the subgroup does exist, and does not stabilize a line on either module, then it might still be strongly imprimitive.

However, the tables in \cite{litterickmemoir} offer a very good first approximation to the list of maximal subgroups we will find in this article.

\subsection{Conjugacy in the algebraic and finite group}
\label{sec:conjfinite}

Many of our potential maximal subgroups $H$ are of the form $\PSL_2(r)$ for some prime power $r$, where $p\nmid r$. A Borel subgroup $L$ of such a group has the form $L_0\rtimes (r-1)/2$, where $L_0$ is an elementary abelian group of order $r$ and $(r-1)/2$ is a cyclic group (unless $r$ is even, in which case the acting group is $(r-1)$). If $L_0$ is non-toral then $L_0$ is one of very few possibilities, and we can work from there. If $L_0$ does lie inside a torus, then we normally can place the whole of $L$ inside $N_\mbG(\mb T)$ for some maximal torus $\mb T$ of $\mbG$. Our first few results attempt to understand the conjugacy classes of such subgroups $L$ in the latter case.

We start by determining the number of classes of complements to $\mb T$ (and its finite version when maximally split) in $\gen{\mb T,w}$, where $w$ is (a preimage of) an element of the Weyl group.

\begin{lem}\label{lem:algconjweylelements} Let $\bar L$ be a cyclic subgroup of the Weyl group $W(\mbG)$, and suppose that $\bar L$ has no trivial constituent on the reflection representation of $W(\mbG)$, or equivalently $C_{\mb T}(\bar L)$ is finite. If $L_1$ and $L_2$ are complements to $\mb T$ in the preimage of $\bar L$ in $N_\mbG(\mb T)$, then $L_1$ and $L_2$ are $\mb T$-conjugate.
\end{lem}
\begin{proof} Note that $\mb T$ is connected. If $L_1=\gen w$ then $C_{\mb T}(w)$ is finite. Conjugation by $w$ is an endomorphism on $\mb T$ with finitely many fixed points, so by the Lang--Steinberg theorem the map $t\mapsto t^{-1}t^w=[t,w]$ is surjective. Thus if $t\in \mb T$ then there exists $u\in \mb T$ such that $t=u^{-1}u^w$. Hence $tw^{-1}=(w^{-1})^u$, so $wt$ and $w$ are always conjugate via an element of $\mb T$.

This shows that any two elements of the coset $w\mb T$ are $\mb T$-conjugate, and therefore any two complements to $\mb T$ in $\gen w\mb T$ are $\mb T$-conjugate, as claimed.
\end{proof}

If we descend to the finite group $G_{\mathrm{sc}}$ then the precise structure of the centralizer $C_{\mb T}(\bar L)$ in the previous lemma becomes important. Of course, $Z(\mbG)$ always lies in the centralizer (if it lies in $G_{\mathrm{sc}}$). Notice that $Z(\mbG)\leq G_{\mathrm{sc}}$ if and only if $G$ possesses non-trivial diagonal automorphisms, and then $Z(G)$ and $G_{\mathrm{ad}}/G$ have the same order.

\begin{lem}\label{lem:finiteconjweylelements} Let $T$ be a finite abelian group and let $w$ be an element of $\Aut(T)$. Let $H$ be the group $T\rtimes \gen w$. There are at most $n$ conjugacy classes of complements to $T$ in $H$, where $n$ is the largest divisor of $|C_T(w)|$ that divides a power of $o(w)$ (in other words, it is the $\pi$-part of $|C_T(w)|$, where $\pi$ is the set of primes dividing $o(w)$).
\end{lem}
\begin{proof} Write $T=T_0\times T_1$, where $\gcd(o(w),|T_1|)=1$ and $|T_0|$ divides a power of $o(w)$. By the Schur--Zassenhaus theorem, all complements to $T_1$ in $\gen{T_1,w}$ are $T_1$-conjugate, so it suffices to assume that $T=T_0$, so $n=|C_T(w)|$.

Note that the map $t\mapsto [t,w]$ is a homomorphism as $T$ is abelian---this is because $[tu,w]=[t,w]^u[u,w]$ in general---and the kernel of the map is $C_T(w)$, and the image is $[T,w]$, a subgroup of $T$. Note that $wt$ and $wu$ are conjugate if and only if $t^{-1}u$ lies in $[T,w]$, so the coset $wT$ splits into exactly $|C_T(w)|$ classes.

As not all of them need consist of elements of order $o(w)$, this means that there are at most $n$ classes of complements, as needed.
\end{proof}

Note that in general, if $w$ has order $n$, then $(tw)^n$ is the product of the elements $t^{(w^i)}$ for $0\leq i\leq n-1$, so this can be constructed inside $\mb T$. Also note that, if one knows that $C_{\mb T}(w)$ is finite, one may check whether $C_{\mb T}(w)=Z(\mbG)$ by simply checking whether $C_Q(w)\leq Z(\mbG)$ for the elementary abelian subgroups $Q$ of $\mb T$ for the primes dividing $o(w)$, and the subgroup of elements of order dividing $9$ if $|Z(\mbG)|=3$ and $3\mid o(w)$.

\medskip

Suppose we are given the action of $L$ on $M(\mbG)$ and $L(\mbG)$, and we wish to know if this determines $L$ up to $\mbG$-conjugacy. The Brauer character values of $L$ on $M(\mbG)$, and perhaps $L(\mbG)$ also, usually determine the $\mbG$-class of both $L_0$ and a complement $\gen w$ to $L_0$ in $L$. The purpose of the above results is to show what extra information is needed to see that if $L_0$ and $\gen w$ are both determined up to $\mbG$-conjugacy, then $L\cong L_0\rtimes \gen w$ is as well. 

Suppose that $L$ and $\hat L$ are two (isomorphic) subgroups of $\mbG$ with isomorphic module actions on $M(\mbG)$ and $L(\mbG)$, with derived subgroups $L'$ and $\hat L'$ isomorphic to $L_0$ and $\mbG$-conjugate, with $L=L'\rtimes \gen w$. Then we may assume that $L'=L\cap \hat L$. If $L'$ is regular, then $C_\mbG^\circ(L')=\mb T$ for some maximal torus $\mb T$. Any element of $\mbG$ that normalizes $L'$ normalizes $C_\mbG^\circ(L')$, whence $L$ and $\hat L$ both lie in $N_\mbG(\mb T)$. Then Lemma \ref{lem:algconjweylelements} shows that $L$ and $\hat L$ are $\mb T$-conjugate.

Thus if $L'$ is regular and $C_{\mb T}(w)$ is finite, then $L$ is unique up to $\mbG$-conjugacy. We explicitly state this as a lemma now.

\begin{lem}\label{lem:torusisgood} Let $r$ be a prime power with $p\nmid r$, and let $L$ and $\hat L$ be subgroups of $\mbG$ that are isomorphic to a Borel subgroup of the group $\PSL_2(r)$. Suppose that $L$ and $\hat L$ have isomorphic actions on $M(\mbG)$ and $L(\mbG)$. Suppose that the derived subgroups $L'$ and $\hat L'$ are $\mbG$-conjugate, and complements $\gen w$ and $\gen{\hat w}$ to $L'$ and $\hat L'$ in $L$ and $\hat L$ respectively are also $\mbG$-conjugate. If $L'$ is regular, and $C_{\mb T}(w)$ is finite for some maximal torus $\mb T$ with $L\leq N_\mbG(\mb T)$, then $L$ and $\hat L$ are $\mbG$-conjugate.
\end{lem}

If one of these conditions does not hold, we will often have to work harder to prove uniqueness of $L$ up to $\mbG$-conjugacy. For an example where the hypotheses, and conclusion, of this lemma do not apply, see Proposition \ref{prop:frobenius55} later, which classifies subgroups $11\rtimes 5$ in $E_6$. In that situation, although $L'$ and $\gen w$ \emph{are} determined up to $\mbG$-conjugacy, $L$ \emph{is not}.

\bigskip

This also seems the right place to place two fundamental results of Larsen and Serre, which allow us to move between algebraic groups for different primes. The first result is Serre's. It allows one to move subgroups of an algebraic group in characteristic $0$ to subgroups in characteristic $p$, maintaining their Brauer characters on all highest weight modules. Since it is a summary of the results of Section 5 of a paper, we give a brief proof using results from that section.

\begin{thm}[Serre {{\cite[Section 5]{serre1996}}}]\label{thm:serre} Let $p$ be a prime and let $k$ be an algebraically closed field of characteristic $p$. Let $\mb X_0$ denote a simple algebraic group of adjoint type over an algebraically closed field of characteristic $0$, and let $H_0$ be a finite subgroup of $\mb X_0$ such that $O_p(H_0)=1$. If $\mb X_p$ denotes a simple algebraic group of the same type as $\mb X_0$ over an algebraically closed field of characteristic $p$, then there exists a finite subgroup $H_p$ of $\mb X_p$ such that $H_0\cong H_p$, and for all highest weights $\lambda$, $W(\lambda){\downarrow_{H_p}}$ for $\mb X_p$ is a reduction modulo $p$ of $W(\lambda){\downarrow_{H_0}}$ for $\mb X_0$.
\end{thm}
\begin{proof} We summarize how to deduce this result from the results of \cite[Section 5]{serre1996} here, because there is no direct theorem in that paper that looks like this.

Let $K$ be an algebraic number field, with ring of integers $\mc O$, with a residue field $\F$ of characteristic $p$. If $H_0$ is a subgroup of $\mb X(K)$ that can be conjugated into $\mb X(\mc O)$, then the map $\mb X(\mc O)\to \mb X(\F)$ restricts to $H_0$ with kernel a $p$-group by \cite[Lemme 4]{serre1996}. As $O_p(H_0)=1$, this means that we obtain a subgroup $H_p$ of $\mb X(\F)$ with the same Brauer character on all highest weight modules.

Not all subgroups of $\mb X(K)$ may be conjugated into $\mb X(\mc O)$ (called a \emph{good reduction} in \cite{serre1996}), but by \cite[Proposition 8]{serre1996}, given $H_0$ there is a finite, totally ramified field extension $K'/K$ such that $H_0$ has a good reduction for $\mb X(K')$. We thus map $H_0$ into $\mb X(\mc O_{K'})$ and then into $\mb X(\F')$ for some field $\F'$ of characteristic $p$.

Any finite subgroup of $\mb X_0$ is conjugate to one over $\mb X(K)$ for some algebraic number field $K$, and this therefore completes the proof.
\end{proof}

This says nothing about whether the subgroups $H_0$ and $H_p$ in the theorem are Lie primitive. Indeed, one may be Lie primitive while the other is not: in this paper, $\PSL_2(8)$ is Lie primitive in $F_4(\C)$, but modulo $7$ it is contained in $G_2$. Another phenomenon is that the $\PSL_2(8)$ subgroup of $D_4$ inside $F_4$ has the same composition factors on $L(F_4)$ as the Lie primitive one (see \cite[Table 6.15]{litterickmemoir}). Similarly, there is a Lie primitive subgroup $M_{12}$ in $E_7$ in characteristic $5$ acting with composition factors of dimensions $32,12^2$ on $M(E_7)$, which are the same as the copy of $M_{12}$ inside $D_6$ in characteristic $0$. However, the author is not aware of any example where a Lie imprimitive subgroup $H$ in characteristic $0$ has the same characters (modulo $p$) as a Lie primitive subgroup in characteristic $p$ on both $M(\mbG)$ and $L(\mbG)$, never mind all highest-weight modules.

Larsen's result however establishes this bijection between Lie primitive subgroups in characteristics $p$ and $0$ whenever $p\nmid |H|$.

\begin{thm}[Larsen {{\cite[Theorem A.12]{griessryba1998}}}]\label{thm:larsen} Let $\mb X_0$, $\mb X_p$ be as in Theorem \ref{thm:serre}, let $H$ be a finite group and suppose that $p\nmid |H|$. There is a one-to-one correspondence between conjugacy classes of subgroups of $\mb X_0$ isomorphic to $H$ and conjugacy classes of subgroups of $\mb X_p$ isomorphic to $H$. Such a correspondence can be chosen to respect Lie primitivity.
\end{thm}

Although the condition on Lie primitivity is not mentioned in Larsen's result, it holds simply because Theorem \ref{thm:larsen} is proved for all split group schemes, so induction on the dimension of the algebraic group proves the claim.

\subsection{Trilinear form}
\label{sec:trilinear}

As is well known, the minimal module $M(E_6)$ of $E_6$ supports a unique (up to scalar) symmetric trilinear form. Hence one might be able to prove that a given finite subgroup $H$ of $\GL_{27}(k)$ lies in a given copy of $\mbG=E_6(k)$ by showing that $H$ leaves this trilinear form invariant, or equivalently that the form yielding $\mbG$ lies inside the space of $H$-invariant symmetric trilinear forms.

In particular, we immediately see the following, which is occasionally useful. (It is not very often useful because the condition rarely holds.) It shows that a certain condition on trilinear forms implies uniqueness up to conjugacy in $\mbG$.

\begin{lem}\label{lem:uniquesymform} Let $\mbG$ have type either $F_4$ or simply connected $E_6$. Suppose that $H$ and $\hat H$ are isomorphic subgroups of $\mbG$, and have isomorphic actions on the minimal module $M(\mbG)$. Let $M$ denote the restriction of $M(\mbG)$ to $H$. If $S^3(M)^H$ is $1$-dimensional then $H$ and $\hat H$ are $\mbG$-conjugate.
\end{lem}

More often we have the following setup, which is also described in more detail in \cite[Section 13]{craven2019un}. Let $V$ be a $k$-vector space, and equip $V$ with a symmetric trilinear form whose symmetry group is $\mbG$, either $F_4$ or (simply connected) $E_6$, so that $V\cong M(\mbG)$. (By using this setup we imply that we are not fixing once and for all a single trilinear form, and we will use different forms when constructing different subgroups of $\mbG$.) Let $L$ be a subgroup of $\mbG$, and suppose that we are given $L$ explicitly as linear transformations of $V$, which of course stabilize the trilinear form. Let $\mb X$ denote the normalizer of $L$ in $\GL(V)$. Let $H$ denote \emph{any} overgroup of $L$ in $\GL(V)$, and let $\mc H$ denote the set of conjugates of $H$ that contain $L$. We wish to understand the set $\mc H$, and in particular which elements of $\mc H$ lie in $\mb G$.

The group $\mb X$ acts by conjugation on $\mc H$, and if we assume that all subgroups of $H$ isomorphic to $L$ are $H$-conjugate (for example, if $H$ is of Lie type and $L$ is a Borel subgroup) then $\mb X$ acts transitively on $\mc H$. Of course, the stabilizer of $H\in \mc H$ is $N_{\mb X}(H)$, so we need only consider cosets of $\mb X$ by $N_{\mb X}(H)$. These cosets are often in our cases covered by elements of the centralizer $C_{\GL(V)}(L)$. This group can be described as a multi-parameter matrix with coefficients in a field $k$, so we replace a general element of $\mb X$ by a matrix $g$ with entries in a polynomial ring over $k$.

We then choose $h\in H$: we see that $h^g$ lies in $\mbG$ if and only if, for all $u,v,w\in V$, we have $f(u^{h^g},v^{h^g},w^{h^g})-f(u,v,w)=0$. Since inverting a matrix with polynomial entries is difficult, so we replace $u$ by $u^g$, and so on, and we instead compute
\[ f(u^{hg},v^{hg},w^{hg})-f(u^g,v^g,w^g).\]
This is zero for all $u,v,w$ if and only if $h^g$ lies in $\mbG$. Thus we obtain a series of polynomial equations which must be solved. Solving them yields all possible $g$ that conjugate our fixed $h$ into $\mbG$. One then needs to take orbits under left multiplication by elements of $C_{\GL(V)}(H)$, yielding all possible conjugation \emph{maps} that send $h$ into $\mbG$. Finally, we may act on the set of all $h^g$ by $N_\mbG(L)$. The number of orbits under this simultaneous left and right multiplication is the number of $\mbG$-conjugacy classes of subgroups $H$ (that contain a $\mbG$-conjugate of $L$, and that act with a given representation on $M(\mbG)$).

This method was used to great success in \cite[Section 13]{craven2019un}, solving the outstanding cases of the subgroups $\PSL_3(3)$ and $\PSU_3(3)$ for $p=3$ acting irreducibly on $M(E_6)$. We will employ it here to resolve many similar cases. Each time we do, we present in the supplementary materials explicit matrices and explicit computations with the trilinear form to show our claims, and here do the theoretical work. One particularly important step is to show that $L$ is unique up to $\mbG$-conjugacy (if it is), so that we may give an explicit description of it.

\subsection{Novelty maximals}
\label{sec:novelty}
Let $H$ be a finite subgroup of $\mbG$, and suppose that $H$ is $\sigma$-stable. This section explains how $N_{\bar G}(\bar H)$ could be a novelty maximal subgroup of an almost simple group $\bar G$, so we assume that $N_G(\bar H)$ is \emph{not} a maximal subgroup of $G$. This applies to any almost simple group $\bar G$, not just a group of Lie type. It is also well known, and appears in \cite{wilson1985}, although our version comes mainly from \cite[Section 1.3.1]{bhrd}.

Let $X$ be a subgroup of $G$ and let $K$ be another subgroup such that $N_G(X)<N_G(K)<G$. We assume that $X<K$ as well, but one can do without this condition. Let $\bar G$ be an almost simple group with socle $G$. We want to understand when $N_{\bar G}(X)$ is not contained in $N_{\bar G}(K)$.

First, we may assume that $\bar G=G\cdot N_{\bar G}(X)$, or in other words $N_{\bar G}(X)/N_G(X)\cong \bar G/G$, as otherwise $N_{\bar G}(X)$ is contained in a proper subgroup of $\bar G$ containing $G$. Write $A$ for $\bar G/G$, and view it as a subgroup of $\Out(G)$. If $\bar G\neq G\cdot N_{\bar G}(K)$ then clearly $N_{\bar G}(X)\not\leq N_{\bar G}(K)$. We say that $H$ is a \emph{type I novelty} with respect to $K$ in this case. This is easy to test: we compute the order of $N_{\bar G}(K)/N_G(K)$ and check if it is $|A|$. If it is not, $X$ is a type I novelty with respect to $K$.

Thus we suppose that $X$ is not a type I novelty with respect to $K$. Let $g\in N_{\bar G}(X)$. If $g\not\in N_{\bar G}(K)$ then $K^g\neq K$. However, since $K$ is $A$-stable, there is some $x\in G$ such that $K^{gx}=K$, hence $gx\in N_{\bar G}(K)$. Note that $X$ and $X^{gx}=X^x$ are both subgroups of $K$, so if there exists $y\in N_G(K)$ such that $X^{xy}=X$, then 
\[xy\in N_G(X)\leq N_G(K)\leq N_{\bar G}(K).\]
Since $gxy\in N_{\bar G}(K)$ this implies that $g\in N_{\bar G}(K)$, a contradiction. Thus if $g\not\in N_{\bar G}(K)$ then $X$ and $X^{gx}$ are not $N_G(K)$-conjugate. However, $gx$ and $g$ induce the same outer automorphism on $K$, so there exists an element of $A$ that, in its induced action on $K$, does not stabilize the $N_G(K)$-conjugacy class containing $X$.

Conversely, if there exists $g\in N_{\bar G}(K)$ such that $X$ and $X^g$ are not $N_G(K)$-conjugate, then no element of the coset $N_G(K)g$ can lie in $N_{\bar G}(X)$. Thus if $h\in N_{\bar G}(X)$ satisfies $Gh=Gg$ (which must exist as $\bar G=GN_{\bar G}(X)$), then $h$ cannot lie in $N_{\bar G}(K)$. This is a \emph{type II novelty}.

Thus we have the following.

\begin{lem}\label{lem:findnovelty} Let $X$ be a subgroup of the finite simple group $G$, and let $\bar G$ be an almost simple group with socle $G$. Suppose that $\bar G=G\cdot N_{\bar G}(X)$, and that $N_G(X)$ is not a maximal subgroup of $G$. Then $N_{\bar G}(X)$ is a maximal subgroup of $\bar G$ if and only if, for any subgroup $K$ containing $X$ and such that $N_G(X)<N_G(K)<G$, one of the following holds:
\begin{enumerate}
\item $\bar G\neq G\cdot N_{\bar G}(K)$;
\item $\bar G=G\cdot N_{\bar G}(K)$, and $N_{\bar G}(K)$ fuses more than one $N_G(K)$-conjugacy class of subgroups, one of which contains $X$.
\end{enumerate}
\end{lem}

This lemma will obviously be useful in determining if there are novelty maximal subgroups in certain situations.

\subsection{Lie algebras}
\label{sec:liealgebra}
Let $\mb X$ be a proper reductive subgroup of $\mbG$, and let $H$ be a known subgroup of $\mb X$, so that $H$ embeds in $\mbG$ with a specific action on $L(\mbG)$ arising from the embedding $H\leq \mb X\leq \mbG$. One of the summands of this action will be $L(\mb X){\downarrow_H}$. In some cases, it will be the case that \emph{any} subgroup $\hat H\cong H$ of $\mbG$ that embeds with the same composition factors on $L(\mbG)$ possesses a summand isomorphic to $L(\mb X){\downarrow_H}$, and in this case one could attempt to prove that $\hat H$ lies in a conjugate of $\mb X$. We need a way of proving this, starting with $H$; i.e., first that a specific summand of $L(\mbG){\downarrow_H}$ is a Lie subalgebra of $L(\mbG)$, and second that it is a reductive Lie algebra and to identify it. This result amalgamates \cite[Lemmas 1, 5, 6 and 10]{ryba2002}, with the last part specific to our applications of it, as we always use it to identify a $\mathfrak{g}_2$ subalgebra.

\begin{lem}\label{lem:liesubalgebra} Let $H$ be a subgroup of $\mbG$ and let $W$ be a $kH$-submodule of $L(\mbG)$.
\begin{enumerate}
\item \label{part1} Suppose that 
\[\Hom_{kH}(\Lambda^2(W),W)=\Hom_{kH}(\Lambda^2(W),L(\mbG)),\]
i.e., that the image of any $H$-map $\Lambda^2(W)\to L(\mbG)$ is contained in $W$. Then $W$ is a Lie subalgebra of $L(\mbG)$.
\item If, in addition, $W$ is simple and not induced from some subgroup of $H$, then $W$ is a simple Lie algebra.
\item If, in addition, there is no $kH$-module quotient of $L(\mbG)/W$ isomorphic to $W^*$, then $W$ has a non-singular trace form.
\item If, in addition, $p\geq 5$ and $\dim(W)=14$, then $W$ is a Lie algebra of type $G_2$.
\end{enumerate}
\end{lem}

After this lemma, we suppose we have a finite subgroup $H$ that stabilizes a $14$-dimensional simple summand of $L(\mbG)$, and this has been shown to be a Lie algebra $\mathfrak{g}_2$. We want to push $H$ into an algebraic subgroup $G_2$ of $\mbG$.

Perhaps the easiest way to do this is to use the recent classification of maximal subalgebras of $L(\mbG)$ in good characteristic, due to Premet and Stewart \cite{premetstewart2019}. Although their results are stronger than this, we see that for $\mbG$ either $F_4$ or $E_6$, any $\mathfrak{g}_2$ subalgebra of $L(F_4)$ or $L(E_6)$ (for $p\geq 5$) is contained in a maximal subalgebra, and this is $\mathrm{Lie}(\mb X)$ for some maximal connected subgroup $\mb X$ of $\mb G$.
%

\begin{prop}\label{prop:maxg2} Let $\mathfrak h$ be a $\mathfrak g_2$ subalgebra of the Lie algebra $L(\mbG)$, where $\mbG$ is either $F_4$ or $E_6$, and suppose that $p\geq 5$. If $H\leq \mbG$ stabilizes $\mathfrak h$ and acts on $L(\mbG)$ with factors of dimension either $26,14,6,6$ if $\mbG=F_4$, or $32,32,14$ if $\mbG=E_6$, then $\mathfrak h$ is a maximal $\mathfrak g_2$ and $H$ is strongly imprimitive.
\end{prop}
\begin{proof} Since $H$ has no trivial composition factors on $L(\mbG)$, neither does $\mathfrak h$, and therefore $\mathfrak h$ is not contained in the subalgebra of a parabolic.

In both cases of $F_4$ and $E_6$, none of the maximal semisimple subalgebras apart from $\mathfrak g_2$ has compatible composition factors (see the tables in for example \cite{thomas2016}). Thus $\mathfrak h$ can only be contained in a maximal $\mathfrak g_2$ (and hence is equal to it) and since $H$ stabilizes it, $H$ lies in a $G_2$ maximal subgroup. Both $H$ and $G_2$ stabilize a unique $14$-space of $L(\mbG)$, and thus $H$ is strongly imprimitive by Proposition \ref{prop:intorbit}.
\end{proof}

The next proposition is used to easily prove existence and/or uniqueness for certain subgroups of exceptional algebraic groups that act irreducibly on the Lie algebra.

\begin{prop}\label{lem:rybaproperty} Let $H$ be a finite group and let $V$ be a simple $kH$-module. Suppose that there is no subgroup $L$ of $H$ such that $V$ is induced from a $1$-dimensional module for $L$.
\begin{enumerate}
\item If $\Hom_{kH}(\Lambda^2(V),V)$ is non-zero and $\Hom_{kH}(\Lambda^3(V),V)=0$ then $V$ carries the structure of a simple, $H$-invariant Lie algebra.
\item If $\Hom_{kH}(\Lambda^2(V),V)=k$ then there exists (up to scalar) at most one non-zero Lie product on $V$. If, in addition, there exists a simple algebraic group $\mb X$ such that $H$ embeds in $\mb X$ and $L(\mb X){\downarrow_H}\cong V$, all such subgroups $H$ are $\Aut(\mb X)$-conjugate.
\item If $\Hom_{kH}(\Lambda^2(V),V)=k$ and $\Hom_{kH}(\Lambda^3(V),V)=0$, and $V$ is the Lie algebra of the simple algebraic group $\mb X$, then $H$ embeds in $\mb X$ acting on $L(\mb X)$ as $V$, and any two such subgroups are $\Aut(\mb X)$-conjugate.
\end{enumerate}
\end{prop}

The first part is \cite[Lemma 2]{ryba2002}. The second part is clear since there is a unique alternating bilinear product on $V$, so certainly at most one non-abelian, $H$-invariant Lie structure on $V$. (Note that $\Aut(\mb X)$ is needed here, rather than just $\mb X$, as all of $\Aut(\mb X)$ acts on $L(\mb X)$; this is important for $\mb X$ of type $E_6$, for instance.) The third part is simply the first and second taken together.

If $\Hom_{kH}(\Lambda^2(V),V)=k$ then we say that $V$ \emph{has the Ryba property}, and if in addition $\Hom_{kH}(\Lambda^3(V),V)=0$ we say that $V$ \emph{has the strong Ryba property}. Since all subgroups acting irreducibly on the Lie algebra of an exceptional algebraic group are known \cite{liebeckseitz2004a}---the author has classified them up to conjugacy as well in papers in preparation---the strong Ryba property is less important, but we mention it because it offers quick, and indeed faster than often in the literature, proofs of existence and uniqueness.

\subsection{Actions of outer automorphisms}
\label{sec:outerauto}

Let $H$ be a quasisimple subgroup of $\mbG$, with $\bar H$ its image in $\mbG/Z(\mbG)$. Much of the work in Sections \ref{sec:f4} and \ref{sec:e6} is on understanding the actions of outer automorphisms of $\mbG$ on the $\mbG$-class of $H$ and on $H$ and $\bar H$ themselves. First, $\sigma$ centralizes $H$ if and only if $H\leq G_{\mathrm{sc}}$, so we need to understand the action of $\sigma$ to determine whether $\bar H\leq G$. We need to know if a given Steinberg endomorphism $\sigma$ stabilizes the $\mbG$-class of $H$, and then if it normalizes a $\mbG$-conjugate of $H$, and then if it centralizes a $\mbG$-conjugate of $H$. We also need to understand stability of $\bar H$ under all outer automorphisms of $G$, i.e., determine $N_{\Aut(G)}(\bar H)$, because that is how the maximal subgroups of $\bar G$ are determined.

We start with Steinberg endomorphisms of $\mbG$, so let $\sigma\in\Aut^+(\mbG)$ denote such a map. There are four possible actions for $\sigma$ on $H$. First, $\sigma$ could map $H$ into a different $\mbG$-conjugacy class of subgroups, i.e., fuse conjugacy classes of subgroups. For the other three options, $\sigma$ stabilizes the class containing $H$; by Lemma \ref{lem:fixedpoints}, this means that $\sigma$ normalizes some $\mbG$-conjugate of $H$, so we assume that $\sigma$ normalizes $H$.

The second possible action is that $\sigma$ acts on $H$ as an element of $\Aut(H)$ that does not lie in $\Aut_{\mbG}(H)$. The third is that $\sigma$ acts non-trivially on $H$, but as an element of $\Aut_{\mbG}(H)$, and the fourth is that $\sigma$ centralizes $H$, i.e., $H$ lies in the fixed points $\mbG^\sigma$.

The third and fourth possibilities are the same, in that if $\sigma$ acts as an element of $\Aut_\mbG(H)$ then $\sigma$ centralizes some other $\mbG$-conjugate of $H$. This follows since, in the group $\mbG\gen\sigma$, all elements $g\sigma$ in the coset of $\sigma$ are $\mbG$-conjugate, so there is `only one' option for $\sigma$ (see \cite[Corollary 21.8]{malletesterman}). Thus if $\sigma g$ centralizes $H$ then $\sigma$ centralizes $H^x$, where $x\in \mbG$ conjugates $\sigma g$ to $\sigma$ in $\mbG\gen\sigma$. Thus, if $\sigma$ now denotes a `minimal' Steinberg endomorphism, i.e., one that is not a power of any other Steinberg endomorphism, to compute the minimal $q$ such that $H\leq \mbG$ embeds in $G_{\mathrm{sc}}$, it suffices to determine the smallest power of $\sigma$ that stabilizes the class containing $H$ and acts as an element of $\Aut_\mbG(H)$ on (a $\mbG$-conjugate of) $H$.

We will find that there is always a single $\Aut(\mbG)$-conjugacy class of simple subgroups $H$ with any given character on $M(\mbG)$ and $L(\mbG)$, i.e., $H$ will be determined up to $\mbG$-conjugacy by its Brauer characters on $M(\mbG)$ and $L(\mbG)$ (for $E_6$ there might be two $\mbG$-classes swapped by the graph automorphism). Thus by studying the irrationalities in these characters one can determine whether $\sigma$ stabilizes the class containing $H$. In addition, in most of our cases (but not all, for example $\PSL_2(25)$ in $F_4$), outer automorphisms of $H$ permute the character values of $M(\mbG){\downarrow_H}$ and $L(\mbG){\downarrow_H}$, so it can be detected whether $\sigma$ acts as an outer or inner automorphism on $H$ just by looking at the character values on these modules.

If $\mbG$ is $F_4$ then every outer automorphism in $\Aut^+(\mbG)$ and $\Out(G)$ is Steinberg, so we may assume from now on that $\mbG$ has type $E_6$.

\medskip

If $\gamma$ is a graph automorphism of $E_6$ then again, $\gamma$ either stabilizes the class containing $H$ (and we say that $H$ is \emph{graph-stable}) or fuses two classes. This time, the action of $H^\gamma$ on $M(E_6)$ is the dual of that of $H$ on $M(E_6)$, and this is usually enough to decide whether $H$ is graph-stable (up to conjugacy) or not.

If the class containing $H$ is graph-stable then we also wish to know the outer automorphism of $\bar H$ induced by $\gamma$, i.e., the group $N_{\bar G}(\bar H)$ where $\bar G=\gen{G,\gamma}$. (Note that since there are two classes of graph automorphism in $\Aut(G)$, it may be the case that a particular $\gamma$ can stabilize the class of $\bar H$ without normalizing any element of the class, so we cannot simply ask for the action of $\gamma$ on a conjugate of $\bar H$, but ask for $N_{\bar G}(\bar H)$.) In our cases it is easy to determine which automorphism of $\bar H$ is induced, either because $\Out(\bar H)$ is small or because we can explicitly compute the result for some small field.

\medskip

Having determined when $H\leq \mbG$ embeds in $G_{\mathrm{sc}}$, we now want to know how the $\mbG$-classes of subgroups $H$ break up into $G$- and $G_\mathrm{ad}$-classes of subgroups $\bar H$. First, if $\bar H$ is simple, and $\phi$ is an automorphism of an almost simple group $X$ with socle $\bar H$, and $\phi$ centralizes $\bar H$, then $\phi=1$. To see this, write $A=\Aut(X)$, and note that $\bar H$ is characteristic in $X$ and $Z(X)=1$, so $\bar H\normal A$. Thus $C_A(\bar H)\normal A$, and $C_A(\bar H)\cap X=1$. Hence $[C_A(\bar H),X]=1$, and hence $C_A(\bar H)=1$ as well, since no elements of $A$ centralize $X$. Thus, if $\phi$ centralizes $\bar H$ it centralizes $N_X(\bar H)$.  Thus, if $\bar H\leq G$ then $N_\mbG(H)/Z(\mbG)=N_{\mbG/Z(\mbG)}(\bar H)\leq G_{\mathrm{ad}}$.

We also note that, as with Corollary \ref{cor:centsubgroups}, since the centralizer of $\bar H$ is trivial in $\mbG/Z(\mbG)$, hence connected, the $\mbG/Z(\mbG)$-conjugates of $\bar H$ lying in $G_\mathrm{ad}$ form a single $G_\mathrm{ad}$-conjugacy class.

\medskip

The final automorphisms of $G$ are diagonal automorphisms for $\mbG$ of type $E_6$. Corollary \ref{cor:diagaut} is used to prove that diagonal automorphisms fuse three $G$-classes of subgroups $\bar H$ into a single $G_{\mathrm{ad}}$-class, except possibly for $\bar H\cong \SL_2(8)$. Section \ref{sec:sl28} is devoted to understanding the action of diagonal automorphisms on this group.

\begin{rem} For $G=E_7(q)$ there is also the issue of diagonal automorphisms. For that group, the diagonal automorphism has order $2$, of course, which is much more likely to divide $|N_\mbG(H):HZ(\mbG)|$. In that case (depending on $H$) it can be a significant problem to determine the action of diagonal automorphisms. At the time of writing, one of these -- the normalizer of the subgroup $A_2$ in simply connected $E_7$ -- is still unresolved.
\end{rem}

\section{Subgroups of \texorpdfstring{$\mbG=F_4$}{G=F4}}
\label{sec:f4}

In this section, we let $\mbG$ denote an algebraic group of type $F_4$ in characteristic $p$, and if $p>0$ let $\sigma$ denote a Frobenius endomorphism, with fixed points $\mbG^\sigma=G=F_4(q)$ for some $q$ a power of $p$. Since $Z(\mbG)=1$, $H\cong \bar H$, so $H$ is a simple group, not just a quasisimple group. Thus in this section we do not need to use $\bar H$ and simply consider $H$.

The simple groups $H$ that we consider are given in \cite[Tables 10.1--10.4]{liebeckseitz1999}, and are repeated in Table \ref{t:f4tocheck}. In each case we must determine if there are Lie primitive examples, compute their normalizers, determine the $\mbG$- and $G$-conjugacy classes, and whether there are any maximal such subgroups.
\begin{table}
\begin{center}
\begin{tabular}{ll}
\hline Prime & Group $H$
\\\hline $p\nmid |H|$ & $\PSL_2(q)$, $q=4,7,8,9,13,17,25,27$, $\PSL_3(3)$, $\PSU_3(3)$, ${}^3\!D_4(2)$
\\ $p=2$ & $\Alt(n)$, $n=7,9,10$, $\PSL_2(q)$, $q=13,17,25,27$, $\PSL_3(3)$, $\PSL_4(3)$, $J_2$
\\ $p=3$ & $\PSL_2(q)$, $q=4,7,13,17,25$, ${}^3\!D_4(2)$
\\ $p=5$ & $\Alt(7)$, $\PSL_2(9)$
\\ $p=7$ & $\Alt(7)$, $\PSL_2(q)$, $q=8,13,27$, $\PSU_3(3)$,  ${}^3\!D_4(2)$
\\ $p=11$ & $M_{11}$, $J_1$
\\ $p=13$ & $\PSL_2(q)$, $q=25,27$, $\PSL_3(3)$, ${}^3\!D_4(2)$
\\ \hline
\end{tabular}
\end{center}
\caption{Simple groups inside $F_4$ that are not of Lie type in defining characteristic.}
\label{t:f4tocheck}
\end{table}

We go through each group in turn, collecting them into families: alternating, sporadic, Lie type other than $\PSL_2(r)$, and finally $\PSL_2(r)$. In each case, as discussed in Section \ref{sec:actions}, we will consult the tables in \cite{litterickmemoir} and determine, for each row labelled with `\textbf{P}', whether such an embedding yields a maximal subgroup (potentially a novelty). (Of course, we do not do this if the subgroup has already been considered elsewhere.)

For computing the number of $\mbG$-conjugacy classes, by Theorem \ref{thm:larsen}, if $H$ is a quasisimple group, then we need only consider primes dividing $|H|$, and then all other primes and $p=0$ are a single case. We will not mention this again in this section, as it will be used in almost every case for $H$.

We will also use the ideas from Section \ref{sec:outerauto} extensively when computing the actions of outer automorphisms of $\mbG$ and $G$. The first case where this is done is in Section \ref{sec:psl28inf4}, and we provide more details then than in subsequent computations.

The aim of this section is to prove that if $N_{\bar G}(H)$ is maximal in $\bar G$ then $H$ appears in Table \ref{tab:f4curlyS}. We delay until Section \ref{sec:proofmaximal} the proof that the subgroups here are indeed maximal, by showing that there are no overgroups inside $G$.

\subsection{Alternating groups}

Most alternating groups were proved to be Lie imprimitive in \cite{craven2017}. This left the cases $H\cong \Alt(6)$ and $p=0,3$, with specific actions on $M(F_4)$ in both cases. The case $p=3$ was settled in \cite{craven2015un}, and the case $p=0$ was settled in \cite{pacheraphd}. (Note that Cohen--Wales \cite{cohenwales1997} do not prove that $H$ is always Lie imprimitive.) In all cases we see that $H$ is strongly imprimitive.

\begin{prop}\label{prop:f4;alt} Let $H\cong \Alt(n)$ be a subgroup of $\mbG$ for some $n\geq 5$. Then $H$ is strongly imprimitive.
\end{prop}
\begin{proof} If $H$ stabilizes a line on $M(F_4)^\circ$ or $L(F_4)$ then $H$ is strongly imprimitive by Proposition \ref{prop:fix1space}, and if $H$ stabilizes a $2$-space on $M(F_4)^\circ$ (i.e., $n=5$, $p=2$, and $H\cong \SL_2(4)$) then $H$ is strongly imprimitive by \cite[Proposition 4.7]{craven2015un}. We therefore assume that neither of these statements holds, and use results from \cite{craven2017} to find all possibilities for $n$ and $p$.

If $n=5$ then by \cite[Propositions 5.1--5.4]{craven2017}, $p\neq 2,3,5$, where $H$ lies in $A_2A_2$, and acts on $M(F_4)$ with factors $5^3,4^2,3_1$ (see \cite[Table 6.7]{litterickmemoir}). The $A_2A_2$ acts semisimply on $M(F_4)$ with modules $(10,10)$, $(01,01)$ and $(00,11)$. Then $H$ lies inside a diagonal $A_2$, and indeed inside a diagonal $A_1$ irreducible subgroup $\mb X$ of $A_2$ (since $2\cdot\Alt(5)$ has a $2$-dimensional simple module). Finally, $L(\mb X)$ is a submodule of $L(A_2)$, which is a summand of the diagonal $A_2$, as we see above. Thus $\mb X$ stabilizes $3_1$ as well, and so $H$ is strongly imprimitive by Proposition \ref{prop:intorbit}.

If $n=6$ then from \cite[Section 6]{craven2017} we see that for characteristic $0$ there is a unique set of composition factors on $M(F_4)$ and $L(F_4)$ that have no trivial factors. This is proved to be Lie imprimitive and unique up to $\mbG$-conjugacy in \cite[Theorem 5.2.1]{pacheraphd}. In fact, strong imprimitivity is also shown there.

Examining \cite[Propositions 6.1--6.3]{craven2017}, we see that the other option is $p=3$. Here $H\cong \PSL_2(9)$ and $H$ is strongly imprimitive by \cite[Proposition 10.3]{craven2015un}.
\end{proof}

\subsection{Sporadic groups}

The three cases were considered in \cite{litterickmemoir}, and $H$ was found to be strongly imprimitive. (For $M_{11}$ and $J_1$ and $p=11$, $H$ always stabilizes a line on $M(F_4)$. For $J_2$ and $p=2$, $H$ stabilizes a line on $L(F_4)$. In all three cases this is simply using Proposition \ref{prop:pressure}. Hence $H$ is strongly imprimitive by Proposition \ref{prop:fix1space}.)

\begin{prop} If $H$ is a subgroup of $\mbG$ isomorphic to one of $M_{11}$, $J_1$ and $J_2$, then $H$ is strongly imprimitive.
\end{prop}

\subsection{Cross-characteristic subgroups not \texorpdfstring{$\PSL_2(r)$}{PSL(2,r)}}

In this subsection we must deal with $\PSL_3(3)$, $\PSL_4(3)$, $\PSU_3(3)$ and ${}^3\!D_4(2)$. Recall that we delay proving maximality until Section \ref{sec:proofmaximal}.

\subsubsection{\texorpdfstring{$H\cong \PSL_3(3)$}{H=PSL(3,3)}}

We will prove the following.

\begin{prop} Let $p\neq 3$ and let $H\cong \PSL_3(3)$ act irreducibly on $M(F_4)$. Then $H$ is Lie primitive, unique up to $\mbG$ and $G$-conjugacy (if $H\leq G$), contained in a unique exotic $3$-local subgroup $3^3\rtimes \PSL_3(3)$, and $N_{\bar G}(H)$ is never a maximal subgroup of $\bar G$.
\end{prop}
\begin{proof} In \cite[Theorem A]{cohenwales1995}, Cohen and Wales show that $H\cong \PSL_3(3)$ is always contained in a unique exotic $3$-local subgroup $J=3^3.\PSL_3(3)$ for any $G$ and $p\neq 2,3$ (and this is the normalizer of the $3^3$, see \cite[Table 1]{clss1992} for example). Since $H$ normalizes a unique $3^3$, so does $N_\mbG(H)$, and thus $H=N_\mbG(H)$. Finally, since a $\mbG$-conjugate of $J$ is centralized by any Frobenius endomorphism $\sigma$ (see \cite[Table 1]{clss1992}), so is $H$, and so $N_{\bar G}(H)$ is always contained in $N_{\bar G}(J)$. This completes the proof for $p\neq 2,3$.

Since we do not consider $p=3$, it remains to consider $p=2$. From \cite[Table 6.31]{litterickmemoir} we see that $M(F_4){\downarrow_H}$ is irreducible and the action is unique up to isomorphism. Let $L$ denote a parabolic subgroup $3^2\rtimes (\SL_2(3).2)$ of $H$. By Lemma \ref{lem:normtorus}, the normal subgroup $L_0\cong 3^2$ of $L$ is toral. Let $L_1$ denote a complement to $L_0$ in $L$, which is unique up to $L$-conjugacy.

Using Magma, one can restrict $L(F_4)$ to $L_0$, and the $1$-eigenspace of $L_0$ on $L(F_4)$ is exactly $4$-dimensional. This means that $C_\mbG(L_0)^\circ$ is a maximal torus, say $\mb T$ (containing $L_0$). Since $N_\mbG(L_0)$ normalizes the connected centralizer of $L_0$, $N_\mbG(L_0)$ lies in $N_\mbG(\mb T)$.

We now prove that $L$ is unique up to $\mbG$-conjugacy in $N_\mbG(\mb T)$, given the action of $H$ on $M(F_4)$. To do so we must first understand the conjugacy classes of subgroups $\SL_2(3).2\cong L_1$ inside $N_\mbG(\mb T)$.

By a computer calculation, there are three conjugacy classes of subgroups $K$ of the Weyl group $N_\mbG(\mb T)/\mb T\cong W(F_4)$ isomorphic to $L_1$, which is a $\{2,3\}$-group. For each of these we compute the number of conjugacy classes of complements to $\mb T$ in the preimage $\hat K$ of $K$ in $N_\mbG(\mb T)$. By the Schur--Zassenhaus theorem, to compute this number we only need consider the `$\{2,3\}$-part' of $\mb T$, and since $p=2$ we only have to consider the prime $3$. In other words, the number of conjugacy classes of complements to $\mb T$ in $\hat K$ is determined by the $1$-cohomology of $K$ on the $3$-part of $\mb T$.

It is easy to see that the only simple module for $L_1\cong \SL_2(3).2$ with non-trivial $1$-cohomology in characteristic $3$ is $1_2$, the non-trivial $1$-dimensional module. However, by a computer calculation, the action of each subgroup $K$ of $W(F_4)$ on the subgroup $3^4$ of $\mb T$ splits as the sum of two non-isomorphic $2$-dimensional modules, and hence no layer $(3^n)^4/(3^{n-1})^4$ can have non-zero $1$-cohomology. This proves that all complements to $\mb T$ in $\hat K$ are conjugate.

In conclusion, we have proved that there are exactly three conjugacy classes of subgroups $\SL_2(3).2$ complementing $\mb T$ in $N_\mbG(\mb T)$, each arising from a different class of subgroups of the Weyl group $W(F_4)$. Each of these different classes of subgroups has a different Brauer character on $M(F_4)$ and $L(F_4)$, so only one of them is compatible with the action of $H$. Thus $L_1$ is determined up to conjugacy in $N_\mbG(\mb T)$.

To move from here to proving that $L$ is determined up to conjugacy in $N_\mbG(\mb T)$, we need to understand copies of $L\cong 3^2\rtimes L_1$ inside the subgroup $3^4\rtimes L_1$ in $\mb T\cdot L_1$. There are exactly two subgroups of order $9$ normalized by $L_1$, so we obtain two possible subgroups $L$. These are conjugate in $N_\mbG(\mb T)$, as confirmed by a computer check in the supplementary materials.

Thus we have proved that, given the Brauer characters on $M(F_4)$ and $L(F_4)$, $L$ is unique up to $\mbG$-conjugacy. The action of $L$ on $M(F_4)$ is easy to describe, and is the sum of three simple modules, of dimensions $2$, $8$ and $16$. We therefore use the trilinear form method from Section \ref{sec:trilinear} to prove uniqueness of $H$, which is done in the supplementary materials.

(The polynomial equations involved are much easier than those we considered in \cite{craven2019un}. If $a,b,c$ are the parameters in a generic element of the centralizer, then each of $a,b,c$ is non-zero. We obtain an equation $\omega c^2(a+c)=0$, where $\omega$ is a cube root of $1$, so $a=c$. Another equation is $c^2(b+\omega c)=0$, so $b=\omega c$. This proves that there is a unique subgroup $H$ in $\mbG$ containing $L$.) Thus for $p=2$, we still have that $H$ is contained in an exotic $3$-local subgroup.

In this case, unlike $p>3$, there is a subgroup $\PSL_3(3).2$ of $\mbG$ \cite[p.2821]{nortonwilson1989}, so $N_\mbG(H)=H.2$. This is contained in a maximal subgroup $\PSL_4(3).2$ of $F_4(2)$ (see \cite[p.2820]{nortonwilson1989}). Since both are irreducible subgroups of $\GL_{26}(k)$, there is a unique subgroup $\PSL_4(3).2$ containing any given $\PSL_3(3).2$. (The same holds for the irreducible subgroup $3^3\rtimes \PSL_3(3)$, so \cite[Theorem A]{cohenwales1995} holds without the restriction $p\neq 2$.)

Since $\PSL_3(3).2\cong\Aut(H)$ is already a subgroup of $G$, it cannot form a novelty maximal subgroup. This completes the proof.
\end{proof}

\subsubsection{\texorpdfstring{$H\cong \PSL_4(3)$}{H=PSL(4,3)}}

Here $p=2$. From \cite[Table 6.32]{litterickmemoir}, we see that $H$ acts irreducibly on $M(F_4)$ and there are two such actions up to isomorphism, which are $\Aut(H)$-conjugate.

This group lies in $F_4(2)$ by \cite{nortonwilson1989}, so it exists over all fields. In addition, Norton and Wilson in \cite{nortonwilson1989} show that $H.2$ lies in, and is maximal in, $F_4(2)$, and $H.2^2$ lies in $F_4(2).2$, so we see that $N_\mbG(H)$ and $N_G(H)$ is always $H.2$. This subgroup can only be maximal in $F_4(2)$.

There are three methods to prove that $H\cong\PSL_4(3)$ is unique up to conjugacy in $\mbG$. In all cases, let $L$ denote a copy of $3^3\cdot \PSL_3(3)$ inside $H$ (i.e., an $A_2$-parabolic). All subgroups $L$ of $\mbG$ acting irreducibly on $M(F_4)$ are $\mbG$-conjugate \cite[Table 1]{clss1992}.

First, $M(F_4){\downarrow_H}$ and $M(F_4){\downarrow_L}$ are both irreducible, there are two $H$-classes of subgroups $L$, and $C_{\GL_{26}(k)}(H)=C_{\GL_{26}(k)}(L)$. Also, the group $H.2$ contained in $\mbG$ fuses the two $H$-classes of subgroups $L$, and so there is a unique copy of $H$ containing a given $L$ in $\GL_{26}(k)$. Since there is a unique $\mbG$-class of subgroups $L$, there is a unique $\mbG$-class of subgroups $H$.

Alternatively, one notes that $S^3(M(F_4))^L$ and $S^3(M(F_4))^H$ are both $2$-dimensional, so any $F_4$-form stabilized by $L$ is also stabilized by $H$; since there is a single $\mbG$-conjugacy class of $L$ acting irreducibly on $M(F_4)$, (up to scalar) there is a unique $F_4$-form in $S^3(M(F_4))^L$, hence there is a unique $F_4$-form in $S^3(M(F_4))^H$.

Finally, a third proof does not use $M(F_4)$ at all: the $246$-dimensional module $L(\lambda_3)$ for $\mbG$ (see \cite[Appendix A.50]{luebeck2001}) that is in the exterior square of $M(F_4)$ restricts to $H$ as the sum of a $208$-dimensional and a $38$-dimensional module, by examining the factors of $\Lambda^2(M(F_4){\downarrow_H})$. The second of these restricts to $L$ as the sum of a $26$- and $12$-dimensional module, but the $208$-dimensional submodule of $L(\lambda_3){\downarrow_H}$ is irreducible on restriction to $L$. Therefore the $38$- and $208$-dimensional submodules of $L(\lambda_3){\downarrow_L}$ are preserved by $H$, and indeed by any copy of $\PSL_4(3)$ containing $L$. Since $N_\mbG(H)$ is Lie primitive, we therefore have that $N_\mbG(H)$ is the stabilizer of the $38$-dimensional subspace of $L(\lambda_3)$. Thus again there is a unique subgroup $H\cong \PSL_4(3)$ containing $L$, which lies in $\mbG$ 

\begin{prop} Let $p=2$ and let $H\cong\PSL_4(3)$ act irreducibly on $M(F_4)$. Then $H$ is Lie primitive, unique up to $\mbG$- and $G$-conjugacy (if $H\leq G$), we have $N_{\mb G}(H)=H.2_2$, and $N_{\bar G}(H)$ is maximal if and only if $G\cong F_4(2)$. In this case $N_G(H)=H.2_2$ and $N_{\bar G}(H)=H.2^2$ if $\bar G\neq G$.
\end{prop}

The proof of maximality will be given in Section \ref{sec:proofmaximal}.

\subsubsection{\texorpdfstring{$H\cong \PSU_3(3)$}{H=PSU(3,3)}}
\label{sec:f4,psu33}
Since $\PSU_3(3)\cong G_2(2)'$, we need only consider $p\neq 2,3$, so $p=7$ or $p\neq 2,3,7$. We see from \cite[Table 6.33]{litterickmemoir} that if $p\nmid|H|$ then $H$ must always stabilize a line on $L(F_4)$, hence is strongly imprimitive by Proposition \ref{prop:fix1space}.

For $p=7$, we use the Lie algebra structure, so the ideas from Section \ref{sec:liealgebra}. From \cite[Table 6.34]{litterickmemoir}, either $H$ stabilizes a line on $L(F_4)$ as it has pressure $-3$ (see Proposition \ref{prop:pressure}), or $H$ acts irreducibly on $M(F_4)$ and with factors $26,14,6^2$ on $L(F_4)$. The $14$ is projective so breaks off as a summand, and its exterior square is $14\oplus 21\oplus 28\oplus 28^*$. In particular, by Lemma \ref{lem:liesubalgebra} this means that the $14$ is a Lie subalgebra of $L(F_4)$, and is $L(G_2)$. Finally, from Proposition \ref{prop:maxg2} we see that this is the maximal $\mathfrak{g}_2$ Lie subalgebra of $\mathfrak{f}_4$, and $H$ is strongly imprimitive.

This completes the proof of the following result.

\begin{prop} If $H\cong \PSU_3(3)$ is a subgroup of $\mbG$ then $H$ is strongly imprimitive.
\end{prop}

\subsubsection{\texorpdfstring{$H\cong {}^3\!D_4(2)$}{H=3D4(2)}}

We delay maximality proofs until Section \ref{sec:proofmaximal}, but we state the full result for this group now anyway.

\begin{prop} Let $p\neq 2$ and let $H\cong {}^3\!D_4(2)$ act irreducibly on $M(F_4)^\circ$. Then $H$ is Lie primitive, unique up to $\mbG$- and $G$-conjugacy, we have $N_{\mb G}(H)=H.3$, and $N_{\bar G}(H)$ is maximal if and only if $\bar G\cong F_4(p)$.
\end{prop}
\begin{proof} For $H\cong {}^3\!D_4(2)$ and $p\neq 2$, from \cite[Tables 6.35 and 6.36]{litterickmemoir} we see that $H$ must act irreducibly on $M(F_4)^\circ$ and the module is unique up to isomorphism. By a Magma computation $S^3(M(F_4)^\circ)^H$ is $1$-dimensional, so there is a unique $H$-invariant, symmetric trilinear form on $M(F_4)^\circ$, and this form is invariant under the group $H.3$, to which the irreducible module extends. Thus if $H$ embeds in $\mbG$ then it is unique up to conjugacy by Lemma \ref{lem:uniquesymform}. (In addition, $H$ has the Ryba property for all $p\neq 2$.) This also shows that, if $H$ embeds in $\mbG$ then so does $H.3\cong \Aut(H)$.

There seem to be two proofs in the literature that ${}^3\!D_4(2)$ embeds in $\mbG$. The first is a discussion in \cite[p.489]{norton1988}, and the second is a direct construction involving generalized hexagons \cite{cohen1983}. (The Cohen--Wales paper \cite[6.2]{cohenwales1997} references \cite{norton1988} and a paper that does not contain the result because the reference should point to \cite{cohen1983}, and \cite{liebeckseitz1999} references \cite{cohenwales1997}.) Similarly to Norton's proof, but much easier, is to note that if $\chi$ is the ordinary character for $H$ of degree $52$, then $\langle \Lambda^2(\chi),\chi\rangle=1$ but $\langle \Lambda^3(\chi),\chi\rangle=0$, i.e., it has the strong Ryba property. By Lemma \ref{lem:rybaproperty}, and the fact that $\mathfrak f_4$ is the only simple complex Lie algebra of dimension $52$, $H.3$ embeds in $F_4(\C)$ and is unique up to conjugacy. Hence $H.3$ embeds in $F_4(k)$ for all primes $p$ by Theorem \ref{thm:serre}, and this embedding acts irreducibly on $M(F_4)$. We cannot deduce the strong Ryba property for $p=3,7,13$ (the odd primes dividing $|H|$) directly from it holding for $p=0$. However, we can just check the Ryba property directly in Magma for these primes, and it holds, and we do not need the strong Ryba property since existence is already known. Thus $H$ exists and is unique up to $\mbG$-conjugacy for all $p\neq 2$.

We now must establish the minimal $q$ for which $H$ (and hence $H.3$) can embed in $F_4(q)$. Since the ordinary character of $H$ on the $26$-space is rational, $H$ embeds into $\GL_{26}(q)$ for any $q$ not a power of $2,3,7,13$. In addition, the minimal fields for $p=3,7,13$ are indeed $\F_p$ as well (see \cite[p.251--253]{abc}, or note that the character is the reduction modulo $p$ of an rational character), so $\F_p$ is the minimal field for the module for all odd primes $p$. Since $S^3(M(F_4)^\circ)^H$ is $1$-dimensional, so the unique $H$-invariant symmetric trilinear form is an $F_4$-form, and therefore $H$ embeds in $F_4(p)$. (Alternatively one sees that a Frobenius endomorphism must stabilize the class of $H$ because it is unique, and then act as an inner automorphism of $H.3$ since this is $\Aut(H)$ so all automorphisms are inner. As we saw in Section \ref{sec:outerauto}, this implies that $\sigma$ centralizes a conjugate of $H$.)
\end{proof}

\subsection{Cross-characteristic subgroups \texorpdfstring{$\PSL_2(r)$}{PSL(2,r)}}

Let $H\cong \PSL_2(r)$ for $r$ one of $7,8,13,17,25,27$. We will deal with each case in turn. We will let $L$ denote a Borel subgroup of $H$, and $L_0$ its derived subgroup, a group of order $r$.

\subsubsection{$H\cong \PSL_2(7)$}

For $H\cong\PSL_2(7)\cong\PSL_3(2)$, we may exclude both $p=2$ and $p=7$, and are thus left with $p=3$ and $p\neq 2,3,7$. If $p=3$ then we know that $H$ is strongly imprimitive by \cite{litterickmemoir}, and in fact $H$ always stabilizes a line on $L(F_4)$. Indeed, we prove more.

\begin{prop}\label{prop:allpsl27p=3} Let $p=3$, let $\mb X$ be a simple, simply connected exceptional group of Lie type not of type $G_2$, and let $H\cong \PSL_2(7)$ be a subgroup of $\mb X$, or $H\cong \SL_2(7)$ with $Z(H)=Z(\mb X)$. Then $H$ stabilizes a line on $L(\mb X)^\circ$. In particular, $H$ is always strongly imprimitive in $\mb X$.
\end{prop}
\begin{proof} The non-projective simple $kH$-modules are of dimensions $1$ and $7$, with the projective cover of $7$ being $7/1/7$. Thus a $kH$-module either has a trivial submodule or quotient, or has at least twice as many $7$s as $1$s. From \cite[Tables 6.13, 6.69, 6.162, 6.163 and 6.268]{litterickmemoir}, we see that this is never the case, and therefore $H$ always stabilizes a line on $L(\mb X)^\circ$. That $H$ is strongly imprimitive follows from Proposition \ref{prop:fix1space}.
\end{proof}

If $p=0$, it is shown that the embedding of $H$ into $\mbG$ that is fixed-point-free  on both $M(F_4)$ and $L(F_4)$ lies inside $A_2A_2$ in \cite[Lemma 8.3]{cohenwales1997}. To check strong imprimitivity, we note that the composition factors of $M(F_4){\downarrow_H}$ are $3,3^*,6^2,8$ and $M(F_4){\downarrow_{A_2A_2}}$ acts as $(10,10)\oplus (01,01)\oplus(00,11)$ (so dimensions $9$, $9$ and $8$) \cite[p.246]{thomas2016}, so if $H$ is contained in $A_2A_2$ then both $H$ and $A_2A_2$ stabilize a unique $8$-dimensional subspace. Thus $H$ is strongly imprimitive for all primes $p$ (via Theorem \ref{thm:larsen}).

\begin{prop} If $H\cong \PSL_2(7)$ is a subgroup of $\mbG$ then $H$ is strongly imprimitive.
\end{prop}

\subsubsection{$H\cong \PSL_2(8)$}
\label{sec:psl28inf4}
Since $\PSL_2(8)\cong {}^2\!G_2(3)'$, we assume that $p\neq 2,3$, so the distinct cases are $p=7$ and $p\neq 2,3,7$. Let $L$ denote a Borel subgroup of $H$, of the form $2^3\rtimes 7$.

If $p=7$, then from \cite[Table 6.15]{litterickmemoir} either $H$ stabilizes a line on $M(F_4)$ or the factors are $8^3,1^2$. If $p\neq 2,3,7$ then either $H$ stabilizes a line on $M(F_4)$ or has composition factors $9_1,9_2,8$ for $9_1,9_2$ any two of the three $\Aut(H)$-conjugate $9$-dimensional simple $kH$-modules. In the case $p=7$, the only option that is compatible with an element of order $7$ acting with Jordan blocks from \cite[Table 3]{lawther1995} is a uniserial module $8/1/8/1/8$. Thus $H$ is strongly imprimitive by Proposition \ref{prop:fix1space}, or the other cases mentioned hold.

We will prove the following result, given the action above, with maximality delayed until Section \ref{sec:proofmaximal}.

\begin{prop} Let $H\cong \PSL_2(8)$ be a subgroup of $\mbG$ and let $\sigma$ be a Frobenius endomorphism of $\mbG$ with $G=\mbG^\sigma$. Suppose that $p\neq 2,3$, and that $H$ acts on $M(F_4)$ as either $8/1/8/1/8$ or $9_1\oplus 9_2\oplus 8$.
\begin{enumerate}
\item There is exactly one $\mbG$-conjugacy class of subgroups $H$. If $p\neq 2,3,7$ then $H$ is Lie primitive and $N_\mbG(H)=H$. If $p=7$ then $H$ is contained in a $G_2$ maximal subgroup of $\mbG$, but $N_\mbG(H)=H.3$ is Lie primitive.

\item There is an embedding of $H$ into $G$ if and only if the polynomial $f_5(x)=x^3-x^2-2x+1$ splits over $\F_q$ (i.e., $q\equiv 0,\pm 1\bmod 7$), where $G=F_4(q)$. If $H$ embeds in $G$ then $H$ is unique up to $G$-conjugacy.

\item If $\bar G$ is an almost simple group with socle $G=F_4(q)$, then $N_{\bar G}(H)$ is maximal in $\bar G$ if and only if either $\bar G=G$ and $\F_q$ is the minimal splitting field for $f_5(x)$, or $\bar G=F_4(p^3).3$ and $H$ does not embed in $F_4(p)$. In this second case, $N_{\bar G}(H)=H.3$.
\end{enumerate}
\end{prop}

The proof proceeds in stages.

\medskip\noindent\textbf{Determination of $N_\mbG(H)$}

\medskip\noindent Since $M(F_4){\downarrow_H}$ is not stable under the field automorphism of $H$ for $p\neq 2,3,7$, of course $H.3$ cannot embed in $\mbG$ in these cases. For $p=7$ a computer calculation shows that there are three extensions of the module $8/1/8/1/8$ to the whole of $H.3$, and exactly one of these is self-dual. We prove in the supplementary materials that this extension indeed embeds in $\mbG$. Later in this proof we show that all copies of $H$ are $\mbG$-conjugate, and this shows that $N_\mbG(H)=H.3$ for all subgroups $H$.

\medskip\noindent\textbf{Determination of $L$ up to conjugacy}

\medskip\noindent By \cite[pp.135--136]{cohenwales1997} $L$ is unique up to conjugacy in $\mbG$ if $p\neq 2,7$, although we need this in characteristic $7$ as well. Thus we assume that $p=7$ from now on. From the action of $H$ on $M(F_4)$, it is easy to check that $L$ acts on $M(F_4)$ with structure $(1/1/1/1/1)\oplus 7^{\oplus 3}$. This means that, by Lemma \ref{lem:linestabs}, $L$ lies inside either $B_4$ or a maximal parabolic subgroup. The former is impossible for $p=7$, since $B_4$ acts on $M(F_4)$ as $1\oplus 9\oplus 16$.

Note that the $C_3T_1$-parabolic has a composition factor of dimension $6$ on $M(F_4)$, so $L$ does not embed in that. Also, $L$ does not embed at all in $A_1A_2$, so $L$ is contained in a $B_3T_1$-parabolic subgroup $\mb X$. We need to compute the number of classes of subgroups $L$ inside this parabolic, which means computing the $1$-cohomology of $L$ on the unipotent radical $\mb U$ of $\mb X$. Since $p=7$, $\mb X$ acts on $M(F_4)$ as $1/8/1,7/8/1$. We see that $L$ must act irreducibly on $M(B_3)=L(100)$, and as $7\oplus 1$ on $L(001)$. Thus the dimensions of the $1$-cohomology groups of $L$ on these two modules are $0$ and $1$ respectively.

As $\mb X/[\mb X,\mb X]$ is a torus, $p=7$ and $|L/L'|=7$, we see that $L$ definitely lies in $[\mb X,\mb X]$, which is of the form $\mb U\cdot B_3$. As a module for the $B_3$, $\mb U$ has layers $L(100)$ and $L(001)$, so the $1$-cohomology of $L$ on $\mb U$ has dimension at most the sum of the dimensions of the $1$-cohomologies of the factors, so $0+1=1$. The $T_1$ part of the Levi complement normalizes $\mb U\cdot B_3$, and also acts on $H^1(L,\mb U)$. The zero of this group corresponds to a subgroup of the Levi complement; $B_3$ acts semisimply on $M(F_4)$, and this is incompatible with the structure of $M(F_4){\downarrow_L}$. Thus $L$ must be one of the complements to $\mb U$ corresponding to a non-zero point of $H^1(L,\mb U)$. All non-zero points of $H^1(L,\mb U)$ must be permuted regularly by $T_1$ as this is the only action of $k^\times$ on $k$ (see, for example, \cite[Lemma 3.2.15]{stewart2013}). (Alternatively one may see it directly in terms of the complements, by assuming that some element of $T_1$ normalizes more than one complement and proving that it then normalizes all of them.) 

Thus $L$ is unique up to conjugacy for all $p\neq 2$.
%
%
%

\medskip\noindent\textbf{Determination of $H$ up to conjugacy}

\medskip\noindent Fix some $L$, and suppose that $H$ contains $L$. First assume $p=7$. From the supplementary materials, we see that all other subgroups $J\cong H$ containing $L$ are related by $J=H^g$ for some $g\in C_\mbG(L)$. Since all $L$ are $\mbG$-conjugate, this implies that all $H$ are $\mbG$-conjugate as well.

For $p\neq 2,3,7$, we first note that, from \cite[Table 6.14]{litterickmemoir} we see there are three non-isomorphic, but $\Aut(H)$-conjugate, sets of composition factors for $H$ that have no fixed points on either $M(F_4)$ or $L(F_4)$, and each of these yields a unique embedding (up to conjugacy) into $\mbG$ by \cite[Theorem 6.10 and p.137]{cohenwales1997}. We claim that these three possible embeddings yield conjugate subgroups of $\mbG$ (they obviously yield conjugate subgroups of $\GL_{26}(k)$ because the representations are $\Aut(H)$-conjugate).

The group $N_\mbG(L)$ acts by conjugation on the overgroups of $L$ isomorphic to $H$. We see from \cite[p.136]{cohenwales1997} that there is a subgroup $2^3\cdot \GL_3(2)$, which obviously contains a subgroup $L.3$ not contained in $LC_\mbG(L)$. An element of order $3$ in $L.3\setminus L$ permutes the three classes of elements of order $7$ in $L$, and it is these classes on which the three representations for $M(F_4){\downarrow_H}$ differ. Thus this element of order $3$ must permute the three classes of embeddings by conjugation, hence they are all $\mbG$-conjugate. Thus $H$ is unique up to $\mbG$-conjugacy.

\medskip\noindent\textbf{Action of outer automorphisms}

\medskip\noindent We need only consider the field automorphism $\sigma=F_p$ of $\mbG$, since $p\neq 2$. It cannot fuse classes as $H$ is unique up to $\mbG$-conjugacy, so we may assume that $\sigma$ normalizes $H$, using Lemma \ref{lem:fixedpoints}. If $\sigma$ acts as an inner automorphism of $H$ then $H\leq \mbG^\sigma$ (up to conjugacy), as we saw in Section \ref{sec:outerauto}. Thus we assume that $\sigma$ acts as an element not in $\Aut_\mbG(H)$.

Suppose first that $p\neq 2,3,7$. Notice that the outer automorphism of $H$ permutes the three representations $M(F_4){\downarrow_H}$ with character irrationalities satisfying $f_5(x)$. Thus if they are defined over $\F_p$ then $\sigma$ cannot induce the outer automorphism, and if they are not defined over $\F_p$ then $\sigma$ must induce the outer automorphism. Finally, if $p=7$ then $\Aut(H)$ already embeds in $\mbG$, so $\sigma$ must act as an inner automorphism of $N_\mbG(H)$, and thus centralizes a conjugate of $H$, as claimed. (In particular, $F_{p^3}$ always centralizes (a conjugate of) $H$.)

\subsubsection{$H\cong\PSL_2(13)$}

Here we consider $H\cong \PSL_2(13)$, when $p\neq 13$, and let $L$ be a Borel subgroup of $H$, a group of the form $13\rtimes 6$. From \cite[Tables 6.17--6.20]{litterickmemoir}, we see that either $H$ stabilizes a line on $M(F_4)^\circ$---and hence is strongly imprimitive by Proposition \ref{prop:fix1space}---or $H$ acts on $M(F_4)^\circ$ as the sum of a $12$- and $14$-dimensional module for $p\neq 3,13$, and a $12$- and $13$-dimensional module for $p=3$. The $13$- and $14$-dimensional modules are always unique up to isomorphism. If $p\neq 7,13$ then there are three non-isomorphic, non-$\Aut(H)$-conjugate, but algebraically conjugate, representations of dimension $12$. If $p=7$ then these three representations become isomorphic.

For $p\neq 3$, the restriction of $M(F_4)$ to $L$ is the sum of two copies of each of the two $6$-dimensional modules and one of each of the two $1$-dimensional modules with $13\rtimes 4$ in the kernel. If $p=3$, the structure is
\[ 6_1^{\oplus 2}\oplus 6_2^{\oplus 2}\oplus (1/1)\]
for the $26$-dimensional module $M(F_4)^\circ/1$.

We will prove the following result, delaying maximality proofs until Section \ref{sec:proofmaximal}, as usual.

\begin{prop} Let $H\cong \PSL_2(13)$ be a subgroup of $\mbG$ and let $\sigma$ be a Frobenius endomorphism of $\mbG$ with $G=\mbG^\sigma$. Suppose that $p\neq 13$, and that $H$ acts on $M(F_4)$ as described above.
\begin{enumerate}
\item If $p\neq 7,13$ then there are exactly three $\mbG$-conjugacy classes of subgroups $H$, $H$ is Lie primitive, and $N_\mbG(H)=H.2\cong\PGL_2(13)$. If $p=7$ then there is one $\mbG$-conjugacy class of subgroups $H$, $H$ is contained in a $G_2$ maximal subgroup of $\mbG$, but $N_\mbG(H)=H.2$ is Lie primitive.

\item There is an embedding of $H$ into $G$ if and only if the polynomial $f_5(x)=x^3-x^2-2x+1$ splits over $\F_q$ (i.e., $q\equiv0,\pm1\bmod 7$), where $G=F_4(q)$. If $H$ embeds in $G$ then there are exactly three $G$-conjugacy classes of subgroups $H$ if $p\neq 7$, and exactly one if $p=7$.

\item If $\bar G$ is an almost simple group with socle $G=F_4(q)$, then $N_{\bar G}(H)$ is maximal in $\bar G$ if and only if $p\neq 2$, $\bar G=G$ and $\F_q$ is the minimal splitting field for $f_5(x)$. (For $p=2$ a subgroup ${}^2\!F_4(8)$ contains $H$.)
\end{enumerate}
\end{prop}

The proof proceeds in stages.

\medskip\noindent\textbf{Determination of $N_\mbG(H)$}

\medskip\noindent Note that either $N_\mbG(H)=H$ or $N_\mbG(H)=H.2\cong\PGL_2(13)$. By either \cite[6.9]{cohenwales1997} or \cite{serre1996}, at least one of the classes of $H$ extends to $\PGL_2(13)$, so we expect $N_\mbG(H)\cong\PGL_2(13)$ for all classes. To show this we need to determine the number of $\mbG$-classes, which we do below.

\medskip

\noindent\textbf{Determination of $L$ up to conjugacy}

\medskip\noindent If $p\neq 3,13$ then the subgroup $13\rtimes 3$ of $L$ is unique up to $\mbG$-conjugacy by \cite[Lemma 2.1]{cohenwales1993}. We could use this to show that $L=13\rtimes 6$ is unique, but we also need the case $p=3$, and to understand classes in the finite group. By Lemma \ref{lem:normtorus}, $L\leq N_\mbG(\mb T)$ for some maximal torus $\mb T$, and $L\mb T/\mb T$ has order $6$.

To apply Lemmas \ref{lem:algconjweylelements} and \ref{lem:finiteconjweylelements}, we need to find $C_{\mb T}(w)$ for an element $w$ of order $6$ in $W(F_4)$ (which is unique up to conjugacy), and in particular the $\{2,3\}$-part of it. We check with a computer that there are no elements of order $2$ or $3$ in $\mb T$ centralized by $w$ by examining the action of $w$ on elementary abelian subgroups of a torus, so $w$ is unique up to conjugacy in both $\mbG$ and any finite group $G$ containing $L$.

We now move up to the group $L$. The subgroup $L_0=13$ is unique up to $\mbG$-conjugacy, and since it is regular, the $\mbG$- and $N_\mbG(\mb T)$-classes of subgroups $L_0$ in $\mb T$ are the same. In particular, $L$ is unique up to $\mbG$-conjugacy. Moreover, by Lemma \ref{lem:fixedpoints}, all subgroups $L_0$ are $G$-conjugate, and therefore we obtain that $L$ is unique up to $G$- and $\mbG$-conjugacy, as needed.

\medskip

\noindent\textbf{Determination of $H$ up to conjugacy}

\medskip\noindent Using the method from Section \ref{sec:trilinear}, we compute the number of copies in $\mbG$ of $H$ containing $L$ under the centralizer of $L$ in $\GL_{26}(k)$. We find exactly one for each representation of $H$. This is done in the supplementary materials for $p=2,3,7$, and for $p\neq 2,3,7,13$ it is accomplished in \cite[Theorem 3.1(ii)]{cohenwales1993}. Thus if $p\neq 7,13$, there are exactly three subgroups $H$ containing a given subgroup $L$, one from each class, and for $p=7$ there is exactly one.

\medskip\noindent\textbf{Action of outer automorphisms}

\medskip\noindent We only need consider the field automorphism $\sigma=F_p$ of $G$ for $p$ odd, and the graph and field automorphisms for $p=2$. Here, since $\Aut(H)$ already embeds in $\mbG$, either $\sigma$ fuses classes---only possible if the Brauer character values of $M(F_4){\downarrow_H}$ (which has irrationalities satisfying $f_5(x)$) lie outside $\F_p$---or acts as an inner automorphism of $N_\mbG(H)$, and hence centralizes $N_\mbG(H)$ up to conjugacy. Thus $\sigma$ fuses the three classes of $H$ if $f_5(x)$ does not split over $\F_p$, and centralizes $H$ (up to conjugacy) if it does split.

For $p=2$, we also need the graph automorphism. Note that $F_8$ centralizes (a conjugate of) $H$, so we may assume that the graph automorphism acts on $G=F_4(8)$ and has order $6$. It still must permute the three $G$-classes of subgroups, so its cube must stabilize each class. However, $\Aut(H)$ embeds in $G$, so it therefore acts as an inner automorphism on $N_G(H)$, and hence up to conjugacy $H$ and $N_G(H)$ lie in ${}^2\!F_4(8)$.

\begin{rem}\label{rem:psl213} As mentioned in the introduction, the existence of $\PGL_2(13)$ in ${}^2\!F_4(8)$ was overlooked in \cite{malle1991}. The proof above shows that there are three classes, permuted transitively by the field automorphism of that group. To confirm this answer, the author has produced a direct construction of $\PGL_2(13)$ inside ${}^2\!F_4(8)$, and this is given in the supplementary materials.
\end{rem}

\subsubsection{$H\cong\PSL_2(17)$}

Here we consider $H\cong \PSL_2(17)$, when $p\neq 17$. There are three cases to consider: $p=2$, $p=3$ and $p\neq 2,3,17$. The first two of these are easy to do, and we prove that $H$ is strongly imprimitive in both cases. For the others, we proceed as in the previous section.

\medskip

If $p=2$ then the pressure of $M(F_4){\downarrow_H}$ is non-positive and there are trivial composition factors (see \cite[Table 6.23]{litterickmemoir}), so $H$ centralizes a line on $M(F_4)$ (as $M(F_4)$ is self-dual) by Proposition \ref{prop:pressure}. We therefore have that $H$ is strongly imprimitive by Proposition \ref{prop:fix1space}, agreeing with \cite[Theorem 1]{litterickmemoir}.

\medskip

Now suppose that $p=3$, and we suppose that $H\leq F_4\leq E_6$, so that $H$ acts on $M(E_6)$ as well as $M(F_4)^\circ$. Let $u$ be an element of order $9$ in $H$. The action of $H$ on $M(F_4)^\circ$ (which has dimension $25$ since $p=3$) is $16\oplus 9$ (see \cite[Table 6.22]{litterickmemoir}). The element $u$ acts on this module with Jordan blocks $9^2,7$. Thus $u$ lies in class $F_4(a_1)$ as we see from \cite[Table 6.1]{craven2015un}, and this table gives the Jordan blocks of the action of $u$ on $M(E_6)$, namely $9^2,7,1^2$.

A computer calculation shows that the module $9$ has no extension with the trivial module, but $\Ext^1_{kH}(16,1)$ is $1$-dimensional. The action of $u$ on the module $9\oplus (16/1)$ is $9^2,8$, which does not match the action of $u$ on $M(E_6)$. Thus $H$ must act on $M(E_6)$ as $16\oplus 9\oplus 1^{\oplus 2}$. Hence $H$ is strongly imprimitive by Lemma \ref{lem:nocohomologyonminimal}.

\medskip

For other primes we see from \cite[Table 6.21]{litterickmemoir} that either $H$ stabilizes a line on $M(F_4)$ or acts as a module $17\oplus 9$. The module $17$ is unique up to isomorphism, and there are two non-isomorphic but $\Aut(H)$-conjugate modules $9$. For this embedding, we have the following result, delaying maximality proofs until Section \ref{sec:proofmaximal}, as usual.

\begin{prop} Let $H\cong \PSL_2(17)$ be a subgroup of $\mbG$ and let $\sigma$ be a Frobenius endomorphism of $\mbG$ with $G=\mbG^\sigma$. Suppose that $p\neq 2,3,17$, and that $H$ acts as described above.
\begin{enumerate}
\item There is a unique $\mbG$-conjugacy class of subgroups isomorphic to $H$, $H$ is Lie primitive, and $N_\mbG(H)=H$.

\item The subgroup $H$ embeds into $G=F_4(q)$ if and only if the polynomial $f_1(x)=x^2-x-4$ splits over $\F_q$ (i.e., $q\equiv \pm 1,\pm2,\pm4,\pm8\bmod 17$). In this case there is a unique $G$-conjugacy class of subgroups $H$.

\item If $\bar G$ is an almost simple group with socle $G=F_4(q)$, then $N_{\bar G}(H)$ is maximal in $\bar G$ if and only if, either $\bar G=G$ and $\F_q$ is the minimal field over which $f_1(x)$ splits, or $f_1(x)$ does not split over $\F_p$, $q=p^2$, and $\bar G=G.2$. In this latter case, $N_{\bar G}(H)\cong \PGL_2(17)$.
\end{enumerate}
\end{prop}

As usual, the proof proceeds in stages.

\medskip\noindent\textbf{Determination of $N_\mbG(H)$}

\medskip\noindent Since the module $M(F_4){\downarrow_H}$ is not stable under the outer automorphism of $H$, we must have that $N_\mbG(H)=H$.

\medskip\noindent\textbf{Determination of $L$ up to conjugacy}

\medskip\noindent Note first that $L_0$ of order $17$ in $L$ is regular, and $C_\mbG(L_0)=\mb T$. Also, if $w\in L$ has order $8$, then $C_{\mb T}(w)$ has order $2$. This can either be seen inside $N_\mbG(\mb T)$, or can be proved directly in the supplementary materials. Either way, we see that the centralizer of $w$ in $\mb T$ is finite, so we may apply Lemmas \ref{lem:algconjweylelements} and \ref{lem:torusisgood} to see that $w$, and then $L$, is unique up to conjugacy in $\mbG$.

In the finite group, things are slightly different. Note that $|C_\mbG(L)|=2$ and if $L\leq G$ then $C_\mbG(L)\leq G$. Write $z$ for the involution in $C_\mbG(L)$. Thus by Lemma \ref{lem:finiteconjweylelements} there are at most two classes of elements $w$, and hence at most two classes of subgroups $L$. (Since $C_\mbG(L_0)$ is a torus, by Lemma \ref{lem:fixedpoints} all such subgroups $L_0$ are conjugate.) Note that $N_\mbG(L)$ in particular normalizes $L_0$; we see that $N_\mbG(L)=\gen{L,z}$.

To see that there are exactly two classes of subgroups $L$ in $G$, note that all involutions in $G$ that are $\mbG$-conjugate are $G$-conjugate, so we may assume that if $L$ and $L_1$ are two subgroups $17\rtimes 8$ that are conjugate in $\mbG$, then $C_\mbG(L)=C_\mbG(L_1)$. If $L$ and $L_1$ are conjugate then they are via an element of $C_\mbG(z)$, which is $B_4$ (as $z$ has trace $-6$ on $M(F_4)$, by a computer check). The finite version of this is $\Spin_9(q)$, which has an outer automorphism of order $2$. Since this lies in $B_4$ (as it is diagonal) it cannot normalize $L$, as $N_\mbG(L)\leq G$. Thus it must fuse two classes, so there are exactly two classes of subgroups $L$ in the finite group $G$.

\medskip\noindent\textbf{Determination of $H$ up to conjugacy}

\medskip\noindent Using the standard centralizer method from Section \ref{sec:trilinear}, in the supplementary materials we find exactly two subgroups $H$ containing a given subgroup $L$. These are conjugate via the involution in $C_\mbG(L)$. Since all subgroups $L$ are conjugate in $\mbG$, this means that all subgroups $H$ are conjugate in $\mbG$. Applying Corollary \ref{cor:centsubgroups}, we obtain that $H$ is unique up to $G$-conjugacy as well, whenever $H\leq G$.

(One might be surprised that there are two classes of subgroups $L$ in $G$ but only one of $H$. Each class of subgroups $L$ of $G$ has two overgroups $H$ in $\mbG$: for one class, both of these lie in $G$, and for the other the two classes are swapped by $\sigma$, with them lying in $\mbG^{\sigma^2}$ but not in $\mbG^\sigma$.)

\medskip\noindent\textbf{Action of outer automorphisms}

\medskip\noindent We only need consider the field automorphisms of $G$ since $p$ is odd. Since $H$ is unique up to conjugacy in $\mbG$, a Frobenius endomorphism $\sigma$ cannot fuse classes, so must normalize $H$. As $|\Out(H)|=2$, $H$ is centralized (up to conjugacy) by $F_{p^2}$, so $H\leq F_4(p^2)$. Note that the character $\chi$ of $M(F_4){\downarrow_H}$ has values in $\F_q$ if and only if $f_1(x)$ splits, as is easily checkable, and an element $g\in H$ of order $17$ has non-rational character value. If $\sigma=F_p$ and $f_1(x)$ does not split over $\F_p$, then $\sigma$ has to send $\chi(g)$ to $\chi(g^3)$, so cannot centralize $H$. In particular, we obtain the statement that $N_{\bar G}(H)\cong \PGL_2(17)$ when $f_1(x)$ does not split over $\F_p$.

If $f_1(x)$ does split over $\F_p$ however, then $\sigma$ cannot send $\chi(g)$ to $\chi(g^3)$, so $\sigma$ cannot induce the outer automorphism on $H$. Thus $\sigma$ must induce an inner automorphism on $H$, so some conjugate of $H$ is centralized by $\sigma$, as claimed.

\subsubsection{$H\cong\PSL_2(25)$}

Here we consider $H\cong \PSL_2(25)$, when $p\neq 5$, and let $L$ be the Borel subgroup of $H$, a group of the form $5^2\rtimes 12$. From \cite[Tables 6.24--6.26]{litterickmemoir}, we see that $H$ acts irreducibly on $M(F_4)^\circ$, and there is a unique such action. The restriction to $L$ is the sum of the two simple $12$-dimensional modules and the two $1$-dimensional modules with kernel $5^2\rtimes 4$ unless $p=3$, in which case the structure is
\[ 12_1\oplus 12_2\oplus (1/1)\]
on the $26$-dimensional module $M(F_4)$.

The subgroup $H$ was constructed in $F_4(\C)$ in \cite[6.6]{cohenwales1997}, and therefore occurs in all characteristics by Theorem \ref{thm:serre}. We will prove the following result, as usual delaying proofs of maximality until Section \ref{sec:proofmaximal}.

\begin{prop} Let $H\cong \PSL_2(25)$ be a subgroup of $\mbG$ and let $\sigma$ be a Frobenius endomorphism of $\mbG$ with $G=\mbG^\sigma$. Suppose that $p\neq 5$, and that $H$ acts as described above.
\begin{enumerate}
\item 
There is a unique $\mbG$-conjugacy class of subgroups isomorphic to $H$, $H$ is Lie primitive, and $N_\mbG(H)=H.2$, with the extension being non-split and with the field-diagonal automorphism ($\PSL_2(25).2_3$ in Atlas notation).

\item There is always an embedding of $H$ into $G$, and $H$ is unique up to $G$-conjugacy.

\item If $\bar G$ is an almost simple group with socle $G=F_4(q)$, then $N_{\bar G}(H)$ is maximal in $\bar G$ if and only if $\bar G=G$ and $q=p\neq 2$.
\end{enumerate}
\end{prop}

The proof proceeds in stages.

\medskip\noindent\textbf{Determination of $N_\mbG(H)$}

\medskip\noindent Since $\Out(H)\cong 2^2$, and $M(F_4){\downarrow_H}$ is $\Out(H)$-stable in all cases, we have to check traces for $p$ odd and unipotent actions for $p=2$, so we start with $p$ odd. In $\Aut(H)$ there are two classes of outer involutions, with traces $0$ and $\pm 4$ on any extension of the $26$-dimensional module to $\Aut(H)$. Since these are not traces of involutions in $\mbG$, we see that $N_\mbG(H)$ cannot be a split extension of $H$ by a non-trivial $2$-group. This means that $N_\mbG(H)$ is either $H$ or is the extension with the diagonal-field automorphism. It follows that $N_\mbG(H)>H$ in characteristic $0$, and therefore for all characteristics by Theorem \ref{thm:serre}.

(It is stated in \cite[6.6]{cohenwales1997} that $H.2$ embeds in $F_4(\C)$, but the proof in \cite{cohenwales1997} is incorrect. In \cite{cohenwales1997}, the proof is that the field-diagonal automorphism is present for $p=2$, and therefore is true for all characteristics. This only works if the automorphism is present for a characteristic not dividing $|H|$ (where one may use Theorem \ref{thm:larsen}), and is false in general; for a counterexample, see Section \ref{sec:psl28inf4}.)

In characteristic $2$, the two outer involutions act on the (unique) $26$-dimensional simple module for $H$ with Jordan blocks $2^{13}$ and $2^{11},1^4$. As neither of these appears in \cite[Table 3]{lawther1995}, we see again that $N_\mbG(H)$ is not a split extension with a non-trivial $2$-group, and the result holds again (see also \cite{nortonwilson1989}).

\medskip\noindent\textbf{Determination of $L$ up to conjugacy}

Note that the subgroup $L_0\cong 5^2$ of $L$ is toral by Lemma \ref{lem:normtorus}, and has exactly four trivial factors on $L(F_4)$, whence $C_\mbG^\circ(L_0)=\mb T$ for some maximal torus $\mb T$. Since any element of $\mbG$ that normalizes $L_0$ normalizes $C_\mbG^\circ(L_0)$, we have
\[ L\leq N_\mbG(L_0)\leq N_\mbG(\mb T).\]

We can then show that $L$ is unique up to conjugacy in $N_\mbG(\mb T)$, completing the proof that $L$ is unique up to $\mbG$-conjugacy. First, there is a unique $W(F_4)$-class of groups of order $5^2$ in $\mb T$ whose Brauer character consists only of values $26$ and $1$ (as $M(F_4){\downarrow_L}$ must). The centralizer of this group in $N_\mbG(\mb T)$ is simply $\mb T$, so the centralizer is connected. In addition, we check that $C_{\mb T}(w)=1$ for $w\in L$ of order $12$, and hence all copies of $L$ containing a fixed $L_0$ are $\mb T$-conjugate using Lemma \ref{lem:torusisgood}. Note that $w$ acts on $5^4$ to normalize a $5^2\times 5^2$ decomposition; both $5^2$s are $G$-conjugate to $L_0$.

\medskip\noindent\textbf{Determination of $H$ up to conjugacy}

\medskip\noindent Using the method from Section \ref{sec:trilinear}, we compute the number of copies of $H$ containing $L$ under the centralizer of $L$ in $\GL_{26}(k)$. We find exactly one. This is done in the supplementary materials for $p=2,3,13$, and for $p\neq 2,3,5,13$ it is accomplished in \cite[Theorem 4.2.2]{pacheraphd}. (Uniqueness was not given in \cite[Section 6.6]{cohenwales1997}.)


\medskip\noindent\textbf{Action of outer automorphisms}

\medskip\noindent We only need consider the field automorphisms of $G$ if $p$ is odd, and the graph automorphism (that powers to the field automorphism) for $p=2$. Notice that, by Corollary \ref{cor:centsubgroups}, $H$ is unique up to $G$-conjugacy if $H\leq G$. If $p=2$ then $H\leq F_4(2)$, and in fact $H\leq {}^2\!F_4(2)$
(see \cite{wilson1984a} or, for example, \cite[p.74]{atlas} for a list of maximal subgroups of ${}^2\!F_4(2)$). Thus $H$ is centralized by all outer automorphisms of $G$, but $H$ is not maximal if $p=2$. Thus $p$ is odd.

As $H$ is unique up to $\mbG$-conjugacy, $F_p$ cannot fuse classes, so must normalize $H$. Since $\Out(H)$ has exponent $2$, certainly $F_{p^2}$ always centralizes (up to conjugacy) $H$, so the question is whether $H$ embeds in $F_4(p)$ or whether $F_p$ induces an outer automorphism on $H$.

We could check if $L$ embeds in $G=F_4(p)$. If this is the case then the Frobenius endomorphism $F_p$ permutes the overgroups of $L$ isomorphic to $H$ (since it centralizes $L$), and as we have proved there is exactly one such overgroup, $F_p$ normalizes $H$. But then $F_p$ acts as an automorphism of $H$ centralizing $L$, and thus $F_p$ centralizes $H$. If $p\equiv\pm 1\bmod 5$ then there is a subgroup $5^4$ of $G$, normalized by the Weyl group, so $L$ embeds in $G$.

From, for example, \cite[Table 3]{bmm1993}, the complex reflection group $G_8$ is the automizer of a Sylow $5$-subgroup of $F_4(p)$ for $p\equiv \pm 2\bmod 5$. This contains elements of order $12$, so there is a group $5^2\rtimes 12$ in $F_4(p)$ for these primes. Indeed, as we stated above, given the element of order $12$, there is a unique (up to conjugacy) group of order $25$ that it normalizes, and so the group $5^2\rtimes 12$ must be $\mbG$-conjugate to $L$. In particular, this proves that $H$ embeds in $F_4(p)$ for all primes $p$.

\subsubsection{$H\cong\PSL_2(27)$}

Here we consider $H\cong \PSL_2(27)$, when $p\neq 3$, and let $L$ be a Borel subgroup of $H$, a group of the form $3^3\rtimes 13$.

From \cite[Tables 6.27--6.29]{litterickmemoir}, we see that $H$ acts irreducibly on $M(F_4)$, and if $p\neq 3,7$ then there are three non-isomorphic, but $\Aut(H)$-conjugate, such representations. If $p=7$ then these representations become isomorphic. The restriction to $L$ is the sum of two $13$-dimensional modules in all cases.

We will prove the following result, delaying maximality proofs until Section \ref{sec:proofmaximal}.

\begin{prop} Let $H\cong \PSL_2(27)$ be a subgroup of $\mbG$ and let $\sigma$ be a Frobenius endomorphism of $\mbG$ with $G=\mbG^\sigma$. Suppose that $p\neq 3$, and that $H$ acts as described above.
\begin{enumerate}
\item 
There is a unique $\mbG$-conjugacy class of subgroups isomorphic to $H$, $H$ is Lie primitive, and $N_\mbG(H)=H$. 

\item There is an embedding of $H$ into $G$ if and only if the polynomial $f_5(x)=x^3-x^2-2x+1$ splits over $\F_q$ (i.e., $q\equiv 0,\pm1\bmod 7$), where $G=F_4(q)$. If $H$ embeds in $G$ then $H$ is unique up to $G$-conjugacy.

\item Let $\bar G$ be an almost simple group with socle $G=F_4(q)$. If $p$ is odd, then $N_{\bar G}(H)$ is maximal in $\bar G$ if and only if $\F_q$ is the minimal splitting field for $f_5(x)$, and $N_{\bar G}(H)=H$ if $\bar G=G$, and is $H.3$ if $\bar G=G.3$. If $p=2$ then $N_{\bar G}(H)$ is always maximal in $\bar G$, if $G=F_4(8)$.
\end{enumerate}
\end{prop}

The proof proceeds in stages.

\medskip\noindent\textbf{Determination of $N_\mbG(H)$}

\medskip\noindent First, note that $\PGL_2(27)$ does not embed in $\mbG$ for $p\neq 3$. For $p\neq 2,3$, this is because an outer involution has trace $0$, hence cannot be conspicuous. For $p=2$, an outer involution acts on all $26$-dimensional modules for $\PGL_2(27)$ with blocks $2^{13}$, which does not appear in \cite[Table 3]{lawther1995}. If $p\neq 3,7$ then the $\Aut(H)$-class of simple modules of dimension $26$ has length $3$, so $H.3$ cannot embed in $\mbG$. For $p=7$, we will prove that $|N_\mbG(H):H|\geq 3$, so is equal to $3$ as it cannot be $6$. Thus $N_\mbG(H)=H$ for all $p\neq 3,7$ and $N_\mbG(H)=H.3$ for $p=7$.

\medskip

\noindent\textbf{Determination of $L$ up to conjugacy}

\medskip\noindent By, for example, \cite[Table II]{griess1991}, this group is unique up to conjugacy in $\mbG$, and even in the finite group $G$ (see \cite{clss1992}, for example, or \cite[Theorem 3.5]{cohenwales1997} for a list of several other proofs).

\medskip

\noindent\textbf{Determination of $H$ up to conjugacy}

\medskip\noindent Using the method from Section \ref{sec:trilinear}, we compute the number of copies of $H$ containing $L$ under the centralizer of $L$ in $\GL_{26}(k)$. For $p\neq 3,7$, we find exactly three subgroups $H$, permuted by a generator for $N_\mbG(L)/N_H(L)$, which has order $3$. This is done in the supplementary materials for $p=2,13$, and for $p\neq 2,3,7,13$ it is accomplished in \cite[Theorem 4.3.2]{pacheraphd}. (Uniqueness was not given in \cite[Section 6.5]{cohenwales1997}.) Hence all subgroups $H$ are $\mbG$-conjugate. For $p=7$, in fact $N_\mbG(L)$ normalizes $H$, so $N_\mbG(L)=H.3$ embeds in $\mbG$, as claimed earlier, and again all subgroups are $\mbG$-conjugate.

\medskip\noindent\textbf{Action of outer automorphisms}

\medskip\noindent For odd primes we need only consider the field automorphisms of $G$, whereas for $p=2$ we have the field and the graph automorphisms.

No automorphism can fuse classes as $H$ is unique up to $\mbG$-conjugacy. Thus any automorphism either acts as an outer automorphism on (a conjugate of) $H$, or centralizes (a conjugate of) $H$.

We start with $p=2$, and let $\tau$ denote a generator for $\Out(G)$, where $G=F_4(2^n)$. Thus $\tau$ has order $2n$. Note that $\Out(H)$ is cyclic of order $6$, and that $M(F_4){\downarrow_H}$ requires $\F_8$ to be a subfield to be realized. Thus $n$ is a multiple of $3$, so $o(\tau)$ is a multiple of $6$. Thus $\tau^6$ must centralize $H$, so $H\leq F_4(8)$ and we may assume that $n=3$. Since $H\not\leq {}^2\!F_4(8)$ (as noted in \cite{malle1991}, $H$ has $3$-rank $3$ and ${}^2\!F_4(8)$ has $3$-rank $2$), we have that $\gen\tau$ induces the full group $\Out(H)$ on $H$.

If $p=7$ then $H$ embeds in $F_4(7)$, as we see in the supplementary materials, so $\sigma$ always centralizes $H$, and we may assume that $p\neq 2,3,7$.

For other primes $p$, note that $L$ is centralized by $\sigma=F_p$, so $\sigma$ either permutes the three subgroups $H$ containing $L$ or it normalizes at least one of them. But then $\sigma$ acts as an automorphism of $H$ that centralizes $L$, and this is only $1$. Thus either $\sigma$ permutes the three subgroups or $\sigma$ centralizes $H$. In the former case, $\sigma$ must induce the field automorphism of $H$, which also permutes the three representations of $H$. This can only occur if $f_5(x)$ does not split over $\F_p$, for then the Galois automorphism of $\F_{p^3}$ permutes the roots of $f_5(x)$, hence the Brauer character values of the three representations.

This completes the proof of the proposition.

\begin{rem} Chris Parker and Kay Magaard have produced a proof that $H$ is unique up to conjugacy in $F_4(8)$ that is theoretical in nature, using the symmetric trilinear form directly. As of the time of writing, this proof is not publicly available.
\end{rem}

\section{Subgroups of \texorpdfstring{$\mbG=E_6$}{G=E6}}
\label{sec:e6}

Now let $\mbG=E_6$, so that $M(\mbG)$ has dimension $27$ and $L(\mbG)^\circ$ has dimension $78-\delta_{p,3}$. Table \ref{t:e6tocheck} lists the groups and primes that we need to consider. We start from the list given in \cite{liebeckseitz1999}, and then subtract those that have been dealt with in \cite{craven2017} and \cite{litterickmemoir}, where explicit theorems were tabulated at the start. For other papers, particularly \cite{cohenwales1997}, which includes lots of results about Lie primitive subgroups, we use them as and when the results apply. (This paper is not as useful to us as one might think because many of the proofs are existence proofs, which we will need to reprove on the way to counting numbers of classes.)

\begin{table}
\begin{center}
\begin{tabular}{ll}
\hline Prime & Group
\\\hline $p\nmid |H|$ & $\Alt(7)$, $M_{11}$, $\PSL_2(q)$, $q=4,7,8,9,11,13,17,19,25,27$,
\\ & $\PSL_3(3)$, $\PSU_3(3)$, $\PSU_4(2)$, ${}^3\!D_4(2)$, ${}^2\!F_4(2)'$
\\ $p=2$ & $\Alt(n)$, $n=7,9,10,11,12$, $M_{11}$, $M_{12}$, $M_{22}$,
\\ & $J_2$, $J_3$, $Fi_{22}$, $\PSL_2(q)$, $q=11,13,17,19,25,27$,
\\ & $\PSL_3(3)$, $\PSL_4(3)$, $\PSU_4(3)$, $\Omega_7(3)$, $G_2(3)$
\\ $p=3$ & $\Alt(7)$, $M_{11}$, $\PSL_2(q)$, $q=4,7,11,13,17,19,25$, ${}^3\!D_4(2)$, ${}^2\!F_4(2)'$
\\ $p=5$ & $\Alt(7)$, $M_{11}$, $M_{12}$, $\PSL_2(q)$, $q=9,11,19$, $\PSU_4(2)$, ${}^2\!F_4(2)'$
\\ $p=7$ & $\Alt(7)$, $\PSL_2(q)$, $q=8,13,27$, $\PSU_3(3)$,  ${}^3\!D_4(2)$
\\ $p=11$ & $M_{11}$, $J_1$
\\ $p=13$ & $\PSL_2(q)$, $q=25,27$, $\PSL_3(3)$, ${}^3\!D_4(2)$, ${}^2\!F_4(2)'$
\\ \hline
\end{tabular}
\end{center}
\caption{Simple groups inside $E_6$ that are not of Lie type in defining characteristic.}
\label{t:e6tocheck}
\end{table}

In the general setup, $H$ is a quasisimple subgroup of $\mbG$ such that $\bar H$ is a simple group in $\mbG/Z(\mbG)$, and in $G$ whenever $H$ is centralized by $\sigma$. If $Z(H)=1$ then we may identify $H$ and $\bar H$, and we do so, so $\bar H$ will only make an appearance when $Z(H)=Z(\mbG)>1$. By Theorem \ref{thm:larsen}, if $H$ is a quasisimple group, then we need only consider primes dividing $|H|$, and then all other primes and $p=0$ are a single case. We will not mention this again in this section, as it will be used in almost every case for $H$.

The aim of this section is to prove that, if $N_{\bar G}(\bar H)$ is maximal in $\bar G$ then $\bar H$ appears in Tables \ref{tab:e6curlyS} and \ref{tab:2e6curlyS}. We delay until Section \ref{sec:proofmaximal} the proof that the subgroups here are indeed maximal, by showing that there are no overgroups inside $G$.

\subsection{Alternating groups}

Most alternating groups were dealt with in \cite{craven2017}. This left over a couple of cases, for $\Alt(6)$ and $p=0$. One case was proved to be Lie imprimitive in \cite{craven2017}, but strong imprimitivity was not proved there. The other case was settled in \cite{pacheraphd}. In all cases we can prove that $H$ is strongly imprimitive.

\begin{prop} Let $H\cong \Alt(n)$ be a subgroup of $\mbG$ for some $n\geq 5$, or $H\cong 3\cdot \Alt(n)$. Then $H$ is strongly imprimitive.
\end{prop}
\begin{proof} As in the proof of Proposition \ref{prop:f4;alt}, if $H$ stabilizes a $1$- or $2$-space on $M(E_6)$, or a line on $L(E_6)^\circ$, then $H$ is strongly imprimitive. If $H\cong \Alt(n)\leq \mbG$ then from \cite{craven2017} this is always the case except for $n=6$ and $p=0$, and for $n=7$ and $H$ contained in a maximal $A_2$ subgroup $\mb X$ acting as $19/8$ on $M(E_6)$ (up to duality). In the former case, $H$ is shown to lie in a subgroup of type $C_4$ and be strongly imprimitive in \cite[Theorem 5.3.1]{pacheraphd}. In the latter case, since $H$ acts on $M(E_6)$ with composition factors $6,8,13$ (see \cite[Table 6.45]{litterickmemoir}), and as $6,13/8$, we see that $H$ and $\mb X$ stabilize a unique $8$-space on $M(E_6)\oplus M(E_6)^*$. Hence $H$ is strongly imprimitive by Proposition \ref{prop:intorbit}. If $H\cong 3\cdot \Alt(7)$ then $H$ always stabilizes a line on $L(E_6)$ by \cite{craven2017}.

The rest of the proof deals with the case where $H\cong 3\cdot\Alt(6)$, so assume this from now on. If $p=2$ we use \cite[Proposition 10.3]{craven2019un}, which states that $H$ is always strongly imprimitive. For $p=5$ this result is proved in \cite{craven2017}, but we offer a shorter alternative proof below. Finally, for $p=0$ in \cite{craven2017} we showed that $H$ is Lie imprimitive but did not show strong imprimitivity.

Thus suppose that $H\leq \mb X$ for some maximal, connected, positive-dimensional subgroup $\mb X$, which are listed, for example, in \cite{liebeckseitz2004} or \cite[Theorem 3.1]{thomas2016}. Note that $\mb X$ cannot split over $Z(\mbG)$, which severely curtails the possibilities for $\mb X$, to $A_5A_1$, an $A_5$-parabolic subgroup, $A_2A_2A_2$, and $A_2G_2$.

If $p=0$ then from \cite[Table 6.51]{litterickmemoir} we are left with the case of dimensions $9^2,6,3$. Since $H$ cannot map onto $A_1$, if $\mb X=A_5A_1$ then $H$ lies in $A_5$. This acts on $M(E_6)$ with factors $15,6,6$, which is not compatible with the factors of $H$. Since $\Alt(6)$ does not embed in $G_2$ (see \cite{kleidman1988}), we cannot have a copy of $H\cong 3\cdot \Alt(6)$ in $A_2G_2$ with $Z(H)=Z(\mbG)$. This leaves $A_2A_2A_2$.

Suppose that $\mb X=A_2A_2A_2$, and note that $M(E_6){\downarrow_{\mb X}}$ is the sum of the three modules
\[ (10,01,00),\quad (00,10,01),\quad (01,00,10).\]
In order for $H$ to lie in $\mb X$, $Z(H)$ must act on the three different modules $M(A_2)$ as (up to permutation) the scalar matrices $1$, $\zeta^2$ and $\zeta$, where $\zeta$ is a primitive cube root of unity. Since there is no (non-trivial) $3$-dimensional representation of $\Alt(6)$, $H$ must act trivially on one of the factors; hence there are at least six $3$-dimensional factors in the action of such a subgroup $H$ on $M(E_6)$.

Thus this case cannot occur in $\mbG$. (This agrees with \cite[Table 44]{frey2016} and \cite{cohenwales1997}.)

If $p=5$ then we see from \cite[6.53]{litterickmemoir} that $H$ could act on $L(E_6)$ with factors $10^2,8^7,1^2$ and still not stabilize a line on $L(E_6)$. This case was resolved in \cite{craven2017}, but we present a shorter proof here. If $H$ does not stabilize a line on $L(E_6)$, then $H$ acts on $L(E_6)$ as
\[ 10^{\oplus 2}\oplus P(8)^{\oplus 2}\oplus 8,\]
where $P(8)$ is the projective cover of $8$, of the form $8/1,8/8$. If $u$ denotes an element of order $5$ in $H$, then $u$ acts on this module with Jordan blocks $5^{15},3$, whence we see that $u$ belongs to class $A_4+A_1$ from \cite[Table 6]{lawther1995}. This class acts on $M(E_6)$ with blocks $5^5,2$. The composition factors of $H$ on $M(E_6)$ are either $15,6^2$ or $6^3,3^3$. The former must be semisimple as there are no extensions between factors (an easy computer check), but then $u$ acts with blocks $5^5,1^2$, which is a contradiction. In the second case, to produce a single non-projective summand (as the action of $u$ requires), the action must be up to duality $P(6)\oplus (3/3,6)$. Thus $H$ is Lie imprimitive as it stabilizes a $3$-space on $M(E_6)$, and by simple dimension counting such stabilizers are positive-dimensional.

So far we have followed the proof of \cite[Proposition 6.1]{craven2017}. But now we use our above options for $\mb X$ that $H$ must lie in. The $15+12$ decomposition above rules out $A_2A_2A_2$, $A_2G_2$ and the $A_5$-Levi, hence $A_5A_1$, so $H$ lies in the $A_5$-parabolic. But if $H$ acts on $M(A_5)$ as $6$ then the factors of $M(E_6){\downarrow_H}$ are $15,6^2$, and if it acts as $3^2$ then the factors are $6,3^3$. Neither of these is as required, so $H$ cannot embed in this way.
\end{proof}

\subsection{Sporadic groups}

Here we have the group $M_{11}$ for all primes, and then $M_{12}$ for $p=2,5$, and also $M_{22}$, $J_2$, $J_3$ and $Fi_{22}$ for $p=2$ and $J_1$ for $p=11$.

\subsubsection{$H\cong M_{11}$}
\label{sec:e6,m11}

From Table \ref{t:e6tocheck} we see that $p=2$, $p=3$, $p=5$, $p=11$, or $p\neq 2,3,5,11$. Unless $p=3,5$, $H$ is shown to be strongly imprimitive in \cite[Theorem 1]{litterickmemoir}. Thus we assume that $p=3,5$. If $p=5$ then from \cite[Table 6.58]{litterickmemoir} we have that the composition factors of $H$ on $M(E_6)$ have dimensions $11$ and $16$. If $p=3$ then from \cite[Table 6.59]{litterickmemoir} we find multiple conspicuous sets of factors that have positive pressure but have trivial factors as well. We will show that these always yield stabilized lines.

We prove the following result.

\begin{prop}\label{prop:e6,m11} Let $H\cong M_{11}$ be a subgroup of $\mbG$.
\begin{enumerate}
\item If $p\neq 5$ then $H$ is strongly imprimitive.
\item If $p=5$ then $H$ is either strongly imprimitive or Lie primitive, and in the latter case is contained in a copy of $M_{12}$. The group $N_{\bar G}(H)$ is not maximal in any almost simple group $\bar G$.
\end{enumerate}
\end{prop}
\begin{proof} When $p=3$ we prove that $H$ stabilizes a line or hyperplane on $M(E_6)$, and hence is strongly imprimitive by Proposition \ref{prop:fix1space} again. First, note that if $H$ is Lie imprimitive then it stabilizes a line: since $H$ does not embed in $F_4$ (see Table \ref{t:f4tocheck}), $G_2$ (as it is not in $F_4$), $A_3$ (no non-trivial module of dimension $4$) or $C_4$ (no non-trivial self-dual  module of dimension $8$), we must have that $H$ lies in either a $D_5$-parabolic subgroup, hence stabilizes a line or hyperplane on $M(E_6)$, or inside an $A_5$-parabolic subgroup. But since $H$ has simple modules of dimensions $1$, $5$ and $10$, $H$ must stabilize a line or hyperplane on $M(A_5)$, hence lie in an $A_4$-parabolic subgroup, thus in a $D_5$-parabolic subgroup. Since a $D_5$-parabolic subgroup stabilizes a line or hyperplane, we are done.

Thus we must prove that $H$ is Lie imprimitive. If $H$ stabilizes a line or hyperplane on $M(E_6)$ then $H$ is strongly imprimitive by Proposition \ref{prop:fix1space}, so we may assume that this is not the case. There are two $kH$-modules with non-trivial $1$-cohomology (see \cite[Table A.3]{litterickmemoir}): in the notation of \cite{litterickmemoir}, they are $5^*$ and $10_b$. in \cite[Table 6.59]{litterickmemoir}, we see that either $H$ has negative pressure on either $M(E_6)$ or its dual, hence stabilizes a line on one of them, or has composition factors $10_b^*,5,(5^*)^2,1^2$. Of these factors, only $5^*$ has non-zero $1$-cohomology: by quotienting out by any submodules $10_b^*$ and $5$, we obtain a module $W$ whose socle is either $5^*$ or $5^*\oplus 5^*$, and with two trivial composition factors, but no trivial quotient (by assumption).

A short Magma computation shows that the largest submodule of $P(5^*)$ with composition factors in $\{1,5,5^*,10_b^*\}$ is
\[ 1/5^*/5,10_b^*/1/5^*.\]
This does not have two trivial composition factors that are not quotients, so $5^*$ cannot be the socle of $W$. Thus $W$ must have the socle structure $5,10_b/1,1/5^*,5^*$, and in particular $W$ is the whole of $M(E_6){\downarrow_H}$. Up to isomorphism there is a unique such module that has no trivial quotient, namely
\[ (5/1/5^*)\oplus (10_b^*/1/5^*).\]
Now let $L$ denote a subgroup $\PSL_2(11)$ of $H$. The restriction of this module to $L$ is
\[ 1\oplus (10/1)\oplus 5^*\oplus (5/5^*).\]
Since $L$ stabilizes a hyperplane and does not lie in $F_4$ (as $M(E_6){\downarrow_H}$ is not self-dual, or see Table \ref{t:f4tocheck}) we see that $L$ lies inside a $D_5$-parabolic subgroup acting uniserially with layers $1/16/10$ on $M(E_6)$. In particular, $L$ must act on $M(D_5)$ as $5/5^*$ and hence stabilize a $5$-space on $M(D_5)$, whose stabilizer is an $A_4$-parabolic subgroup $\mb X$ of $D_5$. Of course, $\mb X$ stabilizes this $5$-dimensional subspace of $M(E_6)$, which $H$ also stabilizes. Thus $H$ is Lie imprimitive, as claimed.

\medskip

We finally consider $p=5$, where the factors of $M(E_6){\downarrow_H}$ have dimensions $11$ and $16$. We restrict to $L\cong \PSL_2(11)$, where the simple module $11$ becomes $1\oplus 10_1$, and the $16$ becomes $5\oplus 11$. (Note that $5$ is not self dual.) Since $M(E_6){\downarrow_L}$ has pressure zero, $L$ stabilizes a line or hyperplane of $M(E_6)$, thus $L$ is contained in $F_4$ or a $D_5$-parabolic. As with $p=3$, the former is impossible.

We now classify copies of $\PSL_2(11)$ inside the $D_5$-parabolic that act as $10_1$ on $M(E_6)$, and therefore act as $5\oplus 11$ on the $16$. Clearly there is a unique class in the $D_5$-Levi subgroup, and so we need to understand $H^1(L,16)$, which is $1$-dimensional (see, for example, \cite[Table A.4]{litterickmemoir}). Let $J$ be a copy of $M_{11}$ inside the $D_5$-Levi subgroup acting on $M(E_6)$ with composition factors $10,16,1$. The restriction map $H^1(J,16)\to H^1(L,11\oplus 5)$ is an isomorphism. So every copy of $\PSL_2(11)$ in the $D_5$-parabolic subgroup whose image in the $D_5$-Levi subgroup is correct is contained in a copy of $M_{11}$, a complement with image $J$ in the Levi. Therefore we see that the stabilizer of the unique $11$-space stabilized by $L$ contains a copy of $M_{11}$ inside the $D_5$-parabolic.

Thus if there is another copy $H$ of $M_{11}$ acting on $M(E_6)$ as $16/11$, it is contained in this $11$-space stabilizer, which must be larger than $H$. Thus while $H$ is Lie primitive, it is not maximal. Indeed, we will see in Section \ref{sec:e6,m12} below that there is such a class, and the subspace stabilizer is $M_{12}$.

Since $\Out(H)=1$, there can be no novelty maximals, and this completes the proof.
\end{proof}

\subsubsection{$H\cong M_{12}$}
\label{sec:e6,m12}

The group $M_{12}$ embeds in $\mbG$ for $p=2,5$. If $p=2$ then $H$ is strongly imprimitive by \cite[Theorem 1]{litterickmemoir}, so we assume that $p=5$. In this case, from \cite[Table 6.61]{litterickmemoir} we see that $H$ must act on $M(E_6)$ with composition factors $11,16$. There are two non-isomorphic, $\Aut(H)$-conjugate $11$-dimensional modules, both self-dual, and two dual $16$-dimensional modules, again $\Aut(H)$-conjugate. In order to comply with the unipotent action from \cite[Table 5]{lawther1995}, the module $M(E_6){\downarrow_H}$ cannot be semisimple, and the projective cover of $11_1$ is, for some choice of labelling,
\[ 11_1/16^*/11_2/16/11_1.\]
So there are four non-isomorphic and non-$\Aut(H)$-conjugate indecomposable modules of dimension $27$ with factors of dimension $11$ and $16$, yielding exactly four potential types of embeddings of $M_{12}$ inside $\mbG$. These subgroups were constructed in \cite{kleidmanwilson1990} but maximality was not proved there. We delay our maximality proofs until Section \ref{sec:proofmaximal}.

\begin{prop} Let $p=5$, and let $H\cong M_{12}$ be a subgroup of $\mbG$, and suppose that $H$ acts as described above.
\begin{enumerate}
\item There are exactly four $\mbG$-conjugacy classes of subgroups $H$, one for each representation, and $N_\mbG(H)=H\times  Z(\mbG)$ in all cases. They are swapped in pairs by the graph automorphism.
\item Each class has a representative in $G=E_6(5)$ and there are $4\cdot\gcd(3,q-1)$ distinct $E_6(q)$-conjugacy classes. The group $H$ does not embed in ${}^2\!E_6(q)$ for any $q$.

\item The subgroup $H$ is maximal in $\bar G$ if and only if $\bar G=E_6(5)$, and there are four classes in this group.
\end{enumerate}
\end{prop}
\begin{proof}
Since $M(E_6){\downarrow_H}$ is not stable under the outer automorphism of $H$, we have that $N_\mbG(H)=H\cdot Z(\mbG)$. Since it is not self-dual either we obtain the second statement from (i).

In \cite{kleidmanwilson1990} it is proved that there are exactly four classes of subgroups $M_{12}$ in $E_6(5)$, that are fused into two by the graph automorphism. We thus want to prove that every copy of $M_{12}$ lying in $E_6(5^n)$ is conjugate to one in $E_6(5)$, for all $n\geq 1$. Hence if we fix an isomorphism type of $M(E_6){\downarrow_H}$, we wish to show that all such embeddings are $E_6(5^n)$-conjugate to one in $E_6(5)$. For definiteness, set it to be $16/11_1$.

Let $L$ denote a copy of $M_{11}$ in $M_{12}$ such that $L$ stabilizes a line on $11_1$. (There are two classes $L_1$ and $L_2$ of $M_{11}$ inside $M_{12}$, and $L_i$ fixes a point on $11_i$ but not on $11_{3-i}$. Thus the other $M_{11}$ acts on $M(E_6)$ as $16/11$, as suggested in the previous section.) The restriction of $16/11_1$ to $L$ is
\[ 10\oplus (16/1),\]
with of course the $11_1$ restricting to the $10\oplus 1$ and the $16$ restricting irreducibly. In particular, $L$ stabilizes a unique $11$-dimensional subspace $W$ of $M(E_6)$.

We now count such copies of $M_{11}$ inside $E_6$, proving that there is a unique class in the finite group $G_{\mathrm{ad}}$ as well. Since $L$ stabilizes a line on $M(E_6)$, and $M(E_6){\downarrow_L}$ is not self-dual (as $5$ is not self-dual), we have that $L$ is contained in a $D_5$-parabolic subgroup $\mb X$, which acts on $M(E_6)$ uniserially with factors $10/16/1$. We see that the image $\bar L$ of $L$ inside the $D_5$-Levi subgroup acts as $10$ on the natural and $16$ on the $16$. Let $\mb Y$ denote the preimage of $\bar L$ in $\mb X$. In particular, $H^1(L,16)$ is $1$-dimensional, so there are exactly two conjugacy classes of complements to the unipotent radical in $\mb Y$ by, for example, \cite[Lemma 3.2.15]{stewart2013}. One class of complements is contained in the $D_5$-Levi subgroup and this acts semisimply on $M(E_6)$ as $1\oplus 10\oplus 16$, so $L$ must belong to the other class of complements in $\mb Y$. Thus unique up to $\mbG$-conjugacy subject to that action on $M(E_6)$, as needed.

Since $L$ may be embedded in $E_6(5)$, one may conjugate $H$ so that $L$ is contained in $E_6(5)$. As $L$ stabilizes a unique $11$-space $W$, this subspace is $\sigma$-stable, whence its stabilizer is $\sigma$-stable. But the stabilizer of $W$ contains $M_{12}$, and as $M_{12}$ is Lie primitive and maximal in $E_6(5)$ (see Section \ref{sec:proofmaximal} later), $M_{12}$ \emph{is} the stabilizer of $W$. Thus the Frobenius endomorphism $F_5$ must normalize $H$ and centralize $L\leq H$. There is no non-trivial such automorphism, so $H\leq E_6(5)$.

To see the actions of outer automorphisms, note that diagonal automorphisms must fuse classes together by Corollary \ref{cor:diagaut}, since $N_\mbG(H)=H\cdot Z(\mbG)$. The graph automorphism must fuse the two classes that act as duals to one another on $M(E_6)$, as in \cite{kleidmanwilson1990}. Since $H$ embeds in $E_6(5)$, the field automorphism centralizes $H$.

We have now proved all of this except the statement about ${}^2\!E_6(q)$; but the field automorphisms stabilize the classes of subgroups $H$ and the graph automorphism fuses classes of subgroups $H$, so $H$ can never be $\sigma$-stable if $\sigma$ is not a standard Frobenius endomorphism.
\end{proof}

\subsubsection{$H\cong 3\cdot M_{22}$}

If $\bar H\cong M_{22}$ then $H\cong 3\cdot M_{22}$ and $p=2$ in this case, and by \cite[Theorem 1]{litterickmemoir} $H$ is strongly imprimitive. In fact, $H$ must stabilize a line on $L(E_6)$, as we see from \cite[Table 6.63]{litterickmemoir}.

\begin{prop} Any copy of $H\cong 3\cdot M_{22}$ in $\mbG$ for $p=2$ is strongly imprimitive.
\end{prop}

\subsubsection{$H\cong J_1$}

From Table \ref{t:e6tocheck} we see that $p=11$. The minimal dimension of a representation of $H$ is $7$. There are two conspicuous sets of composition factors for $M(E_6){\downarrow_H}$ in this case: $7^3,1^6$ and an irreducible action $27$. For $7^3,1^6$, there are no extensions between the composition factors (by a computer check, for example, or using the Brauer tree of $H$) so $M(E_6){\downarrow_H}$ is semisimple. Since $H$ centralizes a $6$-space it lies in a $D_5$-parabolic subgroup, then inside a $D_5$-Levi subgroup as it stabilizes a complement to the line, and then inside a $G_2$ subgroup as $H$ must act as $7\oplus 1^{\oplus 3}$ on $M(D_5)$. This is the copy of $G_2$ that lies inside the $A_2G_2$ maximal subgroup of $\mbG$ (see \cite[p.246]{thomas2016}, for example, for the composition factors of maximal reductive subgroups of $\mbG$ on $M(E_6)$).

If $H$ acts irreducibly on $M(E_6)$ then $H$ also lies in a $G_2$ subgroup, but the irreducible one. There are two ways to prove this. The first is to note that the space of $H$-invariant symmetric trilinear forms on the $27$-dimensional simple module is $2$-dimensional, and the same is true for $G_2(11)$. Since there is a unique copy of $G_2(11)$ containing $J_1$ in $\GL_{27}(11)$ (as $J_1$ is unique up to conjugacy in $G_2(11)$, they are both irreducible on the $27$-space and have no outer automorphisms\footnote{If there were more than one copy of $J\cong G_2(11)$ containing $H\cong J_1$ in $\GL_{27}(11)$, then there would exist $g\in \GL_{27}(11)$ such that $H\leq J,J^g$. Then $H^{g^{-1}}\leq J$ so there exists $n\in J$ such that $H^{ng^{-1}}=H$ and $J^{ng^{-1}}=J^g$. Then $ng^{-1}$ normalizes $H$ but not $J$, but our conditions force the normalizer of $H$ inside the normalizer of $J$.}) we see that any $E_6$-form for $J_1$ extends to an $E_6$-form for $G_2(11)$. (Since $G_2$ contains $H$, the space of forms is at most $2$-dimensional, and since there are two classes of $G_2$ in $E_6$, the space of forms cannot be $1$-dimensional. Thus one does not need to check the claim with a computer.)

Alternatively, we note that the action of $H$ on $L(E_6)$ is $14\oplus 64$, and $\Hom_{kH}(\Lambda^2(14),14\oplus 64)$ is $1$-dimensional, with image contained in the $14$. Thus the $14$ is a subalgebra of the Lie algebra of $E_6$, and we may proceed as in the proof for $\PSU_3(3)$ in $F_4$ in Section \ref{sec:f4,psu33} to see that $H$ must lie in a maximal $G_2$ (for exactly the same reason). Notice that $G_2$ also acts on $L(E_6)$ as $14\oplus 64$, so $H$ is strongly imprimitive by Proposition \ref{prop:intorbit}.

\begin{prop} Every copy of $H\cong J_1$ inside $\mbG$ for $p=11$ lies inside a $\sigma$-stable $G_2$ subgroup of $\mbG$. In particular, $H$ is strongly imprimitive.
\end{prop}

\begin{rem}
This case is left unresolved in \cite{aschbacherE6Vun}, being subgroup 14 in UNK. The condition on the field $F$ in that statement, that $x^2+7$ splits over $F$, is unnecessary, as over $\F_{11}$ this factorizes as $(x-2)(x+2)$.
\end{rem}

\subsubsection{$H\cong J_2$}

From Table \ref{t:e6tocheck}, this only occurs for $p=2$. Let $L\cong \PSU_3(3)$ be a subgroup of $H$. From \cite[Proposition 10.4]{craven2019un}, either $L$ stabilizes a line on $M(E_6)$ or $L(E_6)$, or $L$ acts as $14\oplus 32_1\oplus 32_2$ on $L(E_6)$, and $L$ is contained in a $G_2$-subgroup acting as $14\oplus 64$.

From \cite[p.102]{abc}, we see that the dimensions of simple $kH$-modules are $1,6,14,36,64,84,160$. If $L$ acts on $L(E_6)$ as $14\oplus 32_1\oplus 32_2$ then clearly $H$ must act as $14\oplus 64$ for some modules $14$ and $64$, and we see that $H$ is a blueprint for $L(E_6)$.

The permutation module $P_L$ of $H$ on the cosets of $L$ has structure
\[ 1,36/6_1,6_2/1,1/6_1,6_2/1,36,\]
and the only quotient of this not involving a copy of $36$ is $1$. Thus if $L$ stabilizes a line on $M(E_6)$ then so does $H$. (Of course, a copy of $36$ cannot appear in $M(E_6){\downarrow_H}$ since it has dimension $27$.)

The final option is that $L$ stabilizes a line on $L(E_6)$, in which case from the proof of \cite[Proposition 10.4]{craven2019un} the composition factors of $L(E_6){\downarrow_L}$ are $14^2,6^7,1^8$. Since the $kH$-module $36$ restricts to $L$ with composition factors $14^2,6,1^2$, and only $6$-dimensional simple $kH$-modules have non-zero $1$-cohomology by \cite[Table A.2]{litterickmemoir}, we see that $H$ always has non-positive pressure on $L(E_6)$. As $L(E_6)$ is self-dual, this means that $H$ stabilizes a line on $L(E_6)$ by Proposition \ref{prop:pressure}.

In all cases, we have the following result by Propositions \ref{prop:fix1space} and \ref{prop:intorbit}.

\begin{prop} Any subgroup $H\cong J_2$ of $\mbG$ for $p=2$ is strongly imprimitive.
\end{prop}

\subsubsection{$H\cong 3\cdot J_3$}

From Table \ref{t:e6tocheck}, we only need consider $p=2$ here, with an embedding of $H\cong 3\cdot J_3$ inside $\mbG$ with centres coinciding, acting irreducibly on $L(E_6)$ (see \cite[Table 6.66]{litterickmemoir}). A short computer calculation shows that $\Hom_{kH}(\Lambda^2(L(E_6)),L(E_6))$ is $1$-dimensional, hence $J_3$ has the Ryba property and so is unique up to conjugacy inside $\Aut(\mbG)$. There are two $78$-dimensional simple $kH$-modules though, swapped by the outer automorphism of $H$, so $H.2$ cannot embed in $\Aut(\mbG)$. In particular, this means that the graph automorphism of $\mbG$ fuses two $\mbG$-classes of subgroups $H$. This tallies with, and extends to all fields, \cite[Theorem 1]{kleidmanwilson1990}, which states that there are exactly six classes of $J_3$ over $\F_4$, fused by the outer automorphism group. Aschbacher has also provided a computer-free construction, and uniqueness proof, of $J_3$ in \cite{aschbacher1990}. As $|\Out(H)|=2$ and the outer automorphism of $H$ does not stabilize the $78$-dimensional module, $N_\mbG(H)=H$.

In \cite{kleidmanwilson1990} it is shown that the field automorphism of $E_6(4)$ normalizes (but does not centralize) $H$ and all other automorphisms fuse classes. This means that, if $G=E_6(4)$, then $N_{\bar G}(\bar H)$ is maximal in $\bar G$ if and only if $\bar G=G$ or $\bar G$ is $G$ extended by a field automorphism. Furthermore, because the field automorphism normalizes $H$ and the graph automorphism does not, $H$ cannot be centralized by any field-graph automorphism. Thus $\bar H$ does not lie in ${}^2\!E_6(2^n)$ for any $n$. Diagonal automorphisms fuse classes by Corollary \ref{cor:diagaut}.

\begin{prop} If $p=2$, there are exactly two $\mbG$-conjugacy classes of subgroups isomorphic to $H\cong 3\cdot J_3$, each Lie primitive. The graph automorphism swaps the two classes.

In the finite groups $G=E_6(2^{2n})$, there are six classes of subgroups isomorphic to $J_3$. The group $\bar G$ possesses a maximal subgroup $N_{\bar G}(\bar H)$ if and only if $n=1$ and $\bar G$ is either $G$ or $G$ extended by the field automorphism, and there are six such classes in both cases.

The group $J_3$ does not embed in ${}^2\!E_6(2^n)$ for any $n$.
\end{prop}

Maximality of this subgroup was not proved in either \cite{kleidmanwilson1990} or \cite{aschbacher1990}. We delay our maximality proofs until Section \ref{sec:proofmaximal}.

\subsubsection{$H\cong 3\cdot Fi_{22}$}
\label{sec:e6,fi22}

As with $J_3$, $\bar H\cong Fi_{22}$ only embeds in $G$ for $p=2$, and $H\cong 3\cdot Fi_{22}$ acts irreducibly on both $M(E_6)$ and $L(E_6)$ (see \cite[6.2.67]{litterickmemoir}). We again have that $H$ has the Ryba property and so $H$ is unique up to $\Aut(\mbG)$-conjugacy, but this time there is a unique $78$-dimensional simple module, so this module extends to $H.2$ and must also have the Ryba property. (This is forced: since $p=2$ and the index is $2$, the trivial $kH$-submodule of $\Hom_k(\Lambda^2(L),L)$ (where $L$ is the $78$-dimensional module) must extend to a trivial module for $H.2$, since there is only one $1$-dimensional module for $H.2$ when $p=2$.) Thus $H.2$ embeds in $\mbG.2$.

Since the outer automorphism of $H$ inverts $Z(H)$, $N_\mbG(H)=H$. Thus since $H.2$ embeds in $\mbG.2$, we must have a single $\mbG$-class of subgroups $H$, normalized by the graph automorphism.

Since the $\mbG$-class of $H$ contains all subgroups isomorphic to $H$, it is stable under any Frobenius endomorphism, in particular $F_2$. Since clearly $H\not\leq E_6(2)$ (as $Z(\mbG^{F_2})=1$), we must have that $F_2$ normalizes but does not act as an inner automorphism. (We can also see this as the Brauer character values of the $27$-dimensional module lie in $\F_4$ but not in $\F_2$). Hence $H$ lies in $E_6(4)_{\mathrm{sc}}$ but not in $E_6(2)$. Furthermore, since the graph automorphism also normalizes, but clearly cannot centralize, $H$, the graph-field automorphism must act as an inner automorphism, hence $\bar H\leq G={}^2\!E_6(2)$. This agrees with \cite[p.191]{atlas}.

The group of diagonal automorphisms of $G$ must fuse $G$-classes by Corollary \ref{cor:diagaut}. Since there is a single $\Aut(G)$-class of subgroups $\bar H$, there must be three in $G$, permuted by the $\Sym(3)$ of automorphisms, with the graph fixing one class and swapping the other two. (This again agrees with \cite[p.191]{atlas}.)

\begin{prop}\label{prop:e6,fi22} Let $p=2$ and let $H$ be $3\cdot Fi_{22}$. Then $H$ is Lie primitive, $N_\mbG(H)=H$, and there is exactly one $\mbG$-conjugacy class of subgroups $H$, normalized by the graph automorphism of $\mbG$. Furthermore, $\bar H$ embeds in $G$ if and only if $G$ is either $E_6(4^n)$ and ${}^2\!E_6(2^{2n+1})$ for some $n$, and is unique up to $\Aut(G)$-conjugacy in this case.

If $N_{\bar G}(\bar H)$ is maximal in $\bar G$ then $G={}^2\!E_6(2)$, and either $\bar G=G$ and there are three $G$-classes of subgroups $\bar H$, or $\bar G=G.2$ and there is exactly one class of maximal subgroups $\bar H.2$.
\end{prop}

We will prove maximality in Section \ref{sec:proofmaximal}, or see \cite{wilson2018un}.

\subsection{Cross-characteristic subgroups not \texorpdfstring{$\PSL_2(r)$}{PSL(2,r)}}

According to Table \ref{t:e6tocheck}, the groups in this section are $\PSL_3(3)$, $\PSU_3(3)$, $\PSU_4(2)$, ${}^3\!D_4(2)$ and ${}^2\!F_4(2)'$ for all $p$, and $\Omega_7(3)$ and $G_2(3)$ for $p=2$. Maximality proofs can be found in Section \ref{sec:proofmaximal}.

\subsubsection{$H\cong \PSL_3(3)$}

We assume $p\neq 3$, so the cases to consider are $p=2$, $p=13$, and $p\neq 2,3,13$. From \cite[Tables 6.96--6.98]{litterickmemoir} we can find the composition factors of $H\cong \PSL_3(3)$ on $M(E_6)$ in all characteristics (other than $3$). If $p=2$ then $H$ acts with factors $26,1$, so clearly $H$ stabilizes a line or hyperplane on $M(E_6)$ and lies in $F_4$ by Lemma \ref{lem:linestabs}, hence is strongly imprimitive by Proposition \ref{prop:fix1space}. If $p>3$ then one possibility is again that $H$ acts on $M(E_6)$ with factors $26,1$, stabilizes a line or hyperplane on $M(E_6)$ and then lies in $F_4$. The other option is that $H$ acts irreducibly on $M(E_6)$ for $p\neq 13$, and $H$ acts with factors $16,11$ when $p=13$.

\begin{prop} Let $H\cong \PSL_3(3)$ be a subgroup of $\mbG$ for $p\neq 2,3$, and suppose that $H$ acts irreducibly on $M(E_6)$. Then $H$ is Lie primitive, and $N_\mbG(H)=Z(\mbG)\times (H.2)$. There are two $\mbG$-conjugacy classes of subgroups $H$, swapped by the graph automorphism of $\mbG$. Each such subgroup $H$ is contained in a single subgroup $J\cong {}^2\!F_4(2)$.

Each subgroup $H$ of the finite group $G$ is contained in a copy of $J$ inside $G$, and $N_{\bar G}(H)$ is never maximal in any almost simple group $\bar G$.
\end{prop}
\begin{proof}
Let $L$ denote a maximal parabolic of $H$, which has the form $3^2\rtimes \SL_2(3).2$. The action of this on $M(E_6)$ is the sum of $3$-, $8$- and $16$-dimensional simple modules. Such a subgroup must lie in a proper, positive-dimensional subgroup of $\mbG$.

The group $L''\cong 3^2\rtimes Q_8$ acts on $M(E_6)$ as
\[ 1_2\oplus 1_3\oplus 1_4\oplus 8^{\oplus 3},\]
and so stabilizes three distinct lines on $M(E_6)$. Since it acts non-trivially on these lines (and $p\nmid |L|$ so it lies inside a Levi subgroup whenever it lies in a parabolic subgroup), $L''$ lies first in $D_5T_1$ (see Lemma \ref{lem:linestabs} for the line stabilizers on $M(E_6)$), then in a line stabilizer on $M(D_5)$, so $B_4T_1$, then in a line stabilizer on $M(B_4)$, hence $D_4T_2$. Furthermore, the component group of $N_\mbG(D_4T_2)$ is $\Sym(3)$, and this is also the quotient $L/L''$, so the projection of $L''$ on $D_4$ is normalized by all graph automorphisms of $D_4$. Thus the projection of $L''$ is determined up to conjugacy in $D_4$, by for example \cite[Lemma 1.8.10(ii)]{bhrd}. The projection of $L''$ onto the $T_2$ factor is unique (not even just up to conjugacy), as it is simply the Klein four group. Thus $L''$ is determined up to conjugacy in each of $D_4$ and $T_2$.

However, because there is no element of $D_4$ acting as an outer element of order $2$ on the subgroup $L''$, there are \emph{two} $D_4T_2$-conjugacy classes of subgroups $L''$ diagonally embedded in $D_4T_2$: one is obtained by twisting one of the two factors by the graph automorphism of order $2$. Of course, both are stable under the $\Sym(3)$ of graph automorphisms acting simultaneously on both $D_4$ and $T_2$, so both extend to subgroups $L$ in $\mbG$. Notice that the graph automorphism of $\mbG$ induces a graph automorphism of order $2$ on one of the factors $D_4$ and $T_2$, hence swaps these two copies of $L$. (There is an element of order $3$ in $D_4$ acting like the graph automorphism of order $3$ on $L''$, so we see only two classes, not six.)

In the supplementary materials we show that there is a unique copy of $H$ in $\mbG$ above a given copy of $L$. Hence there are exactly two $\mbG$-conjugacy classes of subgroups $H$, swapped by the graph automorphism. We then apply Corollary \ref{cor:centsubgroups} to obtain either six or two $G_{\mathrm{sc}}$-conjugacy classes of subgroups $H$, depending on whether $\Out(G)$ contains diagonal automorphisms or not.

\medskip

Now note that there are two $\mbG$-classes of subgroups ${}^2\!F_4(2)$, each containing a class of subgroups $\PSL_3(3).2$, swapped by the graph automorphism. (See Section \ref{sec:e6,2f42} below, which does not depend on this section.) Thus each copy of $H$ is contained in a copy of $J\cong {}^2\!F_4(2)$, and exactly one since there is a unique copy of $J$ containing $H$ in $\GL_{26}(k)$. Also, since every copy of $H$ is contained in a copy of $\Aut(H)\leq \mbG$, any Frobenius endomorphism $\sigma$ of $\mbG$ either centralizes $H$ or fuses classes, as with $J$. Furthermore, $\sigma$ centralizes $H$ if and only if it centralizes $J$, so $N_{\bar G}(H)$ cannot be maximal for any $\bar G$. In Section \ref{sec:e6,2f42} below we determine that $J$ embeds in $G=E_6(q)$ if $q\equiv 1\bmod 4$ and in ${}^2\!E_6(q)$ for $q\equiv -1\bmod 4$.
\end{proof}

\subsubsection{$H\cong \PSU_3(3)$}

Let $H\cong \PSU_3(3)\cong G_2(2)'$. Here $p=7$ or $p\neq 2,3,7$. In the first case, $H$ stabilizes a line on either $M(E_6)$ or $L(E_6)$, as we see from \cite[Table 6.101]{litterickmemoir}. Hence $H$ is strongly imprimitive by Proposition \ref{prop:fix1space}, as stated in \cite{litterickmemoir}. For $p\neq 2,3,7$, we note that $H$ either stabilizes a line on $M(E_6)$ or $L(E_6)$, or acts irreducibly on $M(E_6)$ by \cite[Table 6.100]{litterickmemoir}. In this latter case, $H$ stabilizes a $\mathfrak{g}_2$ subalgebra and thus is strongly imprimitive by Proposition \ref{prop:maxg2}, as seen in the proof of \cite[Lemma 8.2]{cohenwales1997}. Thus we obtain the following result.

\begin{prop} Any copy of $\PSU_3(3)$ in $\mbG$ is strongly imprimitive.
\end{prop}

\subsubsection{$H\cong \PSL_4(3),\PSU_4(2),3\cdot \PSU_4(3),{}^3\!D_4(2)$}

Each of these was proved in \cite{litterickmemoir} to be strongly imprimitive. Indeed, in all cases $H$ stabilizes a line on either $M(E_6)$ or $L(E_6)$.

\begin{prop} Any copy of $\PSL_4(3)$, $\PSU_4(2)$, $3\cdot \PSU_4(3)$ or ${}^3\!D_4(2)$ in $\mbG$ is strongly imprimitive.
\end{prop}

\subsubsection{$H\cong 3\cdot \Omega_7(3)$, $3\cdot G_2(3)$}

These two cases only occur for $p=2$, and if $\bar H$ is simple then $H\cong 3\cdot \bar H$. In both cases, $H$ must act irreducibly on $L(E_6)$. Also in both cases $\Hom(\Lambda^2(78),78)$ is $1$-dimensional, so $H$ has the Ryba property, and the $78$-dimensional module is unique up to isomorphism. Thus each of these is unique up to $\Aut(\mbG)$-conjugacy, and unique up to $\Aut(G_{\mathrm{ad}})$-conjugacy whenever $\bar H$ embeds in the finite group $G_{\mathrm{ad}}$. Furthermore, it shows that $H$ extends to $H.2$ in $\mbG.2$ (as with $3\cdot Fi_{22}$). As with $3\cdot Fi_{22}$, this means that there is a unique $\mbG$-class of both groups, normalized by the graph automorphism.

The proof for minimal field and action of automorphisms is the same as for $3\cdot Fi_{22}$ in Section \ref{sec:e6,fi22}, so we are a bit less detailed in this case. Again, $\bar H\leq E_6(4)$ but $\bar H\not\leq E_6(2)$ because of the fact that $Z(E_6(2)_{\mathrm{sc}})=1$. Consequently $\bar H\leq {}^2\!E_6(2)$.

Note that 
\[G_2(3)\leq \Omega_7(3)\leq Fi_{22}\leq {}^2\!E_6(2),\]
with the final inclusion appearing in Proposition \ref{prop:e6,fi22} (see also \cite[p.191]{atlas}). Thus we see that $\bar H$ embeds in $E_6(4^n)$ and ${}^2\!E_6(2^{2n+1})$ for all $n$; in particular, such $G$ always have a non-trivial diagonal automorphism, which cannot normalize $\bar H$ by Corollary \ref{cor:diagaut}. Thus if $N_{\bar G}(\bar H)$ is maximal in $\bar G$ then $\bar G$ cannot induce this diagonal automorphism on $G$.

If $N_{\bar G}(\bar H)$ is maximal in $\bar G$, then $G={}^2\!E_6(2)$ and either $\bar G=G$ or $\bar G=G.2$. We can just read off the answer from \cite{atlas} now, and find no classes of maximal $G_2(3)$s, and one class of (novelty) maximal $\Omega_7(3).2$ in the group $\bar G=G.2$. To prove it independently of \cite{atlas}, $\Out(G)$ acts as $\Sym(3)$ on the three conjugacy classes of $\Omega_7(3)$ and $G_2(3)$, and therefore exactly one is normalized in $G.2$, which yields $\Omega_7(3).2$ and $G_2(3).2$. Since the graph automorphism of $Fi_{22}$ interchanges two classes of $\Omega_7(3)$ \cite{kleidmanwilson1987}, we see that this subgroup $\Omega_7(3).2$ cannot lie in $Fi_{22}.2$, so yields a novelty maximal. On the other hand, the group $G_2(3).2$ \emph{is} contained as a novelty maximal subgroup of $Fi_{22}.2$ (see \cite{kleidmanwilson1987} again), hence is not a novelty maximal of $\bar G$.

\begin{prop} Let $p=2$ and let $H$ be one of $3\cdot G_2(3)$ or $3\cdot \Omega_7(3)$. Then $H$ is Lie primitive, $N_\mbG(H)=H$, and there is exactly one $\mbG$-conjugacy class of subgroups $H$, normalized by the graph automorphism of $\mbG$. Furthermore, $\bar H$ embeds in $G$ if and only if $G$ is one of $E_6(4^n)$ and ${}^2\!E_6(2^{2n+1})$ for some $n$, there are exactly three $G$-conjugacy classes of subgroups in this case, permuted transitively by the diagonal automorphism of $G$.

If $N_{\bar G}(\bar H)$ is maximal in $\bar G$ then $\bar H\cong \Omega_7(3)$, $G={}^2\!E_6(2)$, and $\bar G=G.2$, in which case there is one class of novelty maximal subgroups.
\end{prop}

We will prove maximality in Section \ref{sec:proofmaximal}, or see \cite{wilson2018un}.

\subsubsection{$H\cong {}^2\!F_4(2)'$}
\label{sec:e6,2f42}
For all characteristics $p\neq 2,3$, the action of $H$ on $L(E_6)$ is irreducible, and again possesses the Ryba property, so we are done as for the previous two cases. In characteristic $3$ we must be very slightly more careful, and note that $\Hom_{kH}(\Lambda^2(L(E_6)^\circ),L(E_6)^\circ)$ is $1$-dimensional instead. (In characteristic $0$ this was already noted by Cohen and Wales \cite[6.1]{cohenwales1997}.) This $kH$-module homomorphism extends to one for $H.2$, so this group lies in $\mbG.2$ and is also unique up to conjugacy. (In fact, $H$ has the strong Ryba property in characteristic $0$, so existence is also easily proved via Theorem \ref{thm:serre}.)

Note that the ordinary character of $M(E_6){\downarrow_H}$ has field of values $\Q(\I)$ \cite[p.75]{atlas}. The same holds for the Brauer characters when $p=3,5,13$, so we always require that $f(x)=x^2+1$ splits over $\F_q$ for $H$ to be defined over $\F_q$.

As $M(E_6){\downarrow_H}$ is not self-dual, and the $27$-dimensional modules for $H$ are not swapped by the outer automorphism of $H$, we see that $H.2$ in fact is contained in $\mbG$, not just $\mbG.2$, and so there must be two $\mbG$-conjugacy classes of subgroups $H$ (and $H.2$), swapped by the graph automorphism.

Suppose that $p\equiv 3\bmod 4$. In this case the Brauer character values of $M(E_6){\downarrow_H}$ are not fixed by $F_p$, and so $F_p$ must swap the two $\mbG$-classes. Since the graph automorphism also swaps them, the product $\sigma$ normalizes a representative $H.2$, and since $\Aut(H.2)=\Aut(H)$, must act as an inner automorphism. Hence $H.2\leq \mbG^\sigma=G={}^2\!E_6(p)$, and $N_G(H)=H.2$ must be maximal in this group.

The opposite holds if $p\equiv 1\bmod 4$. In this case $F_p$ fixes the Brauer character values and $\sigma$ does not, so $H\leq E_6(p)$.

Diagonal automorphisms always fuse classes by Corollary \ref{cor:diagaut}. We will prove maximality in Section \ref{sec:proofmaximal}.

\begin{prop} Let $H\cong {}^2\!F_4(2)'$ and $p\neq 2$. Then $H$ is Lie primitive, $N_\mbG(H)=Z(\mbG)\times (H.2)$ and there are exactly two $\mbG$-conjugacy classes of subgroups $H$, swapped by the graph automorphism of $\mbG$. If $H\leq G$ then $H$ is unique up to $\Aut(G)$-conjugacy.

If $p\equiv 1\bmod 4$ then for all $n$, $H\leq E_6(p^n)$, and if $p\equiv 3\bmod 4$ then for all $n$, $H\leq {}^2\!E_6(p^{2n+1}),E_6(p^{2n})$. If $N_{\bar G}(H)$ is maximal in $\bar G$ then $\bar G=G$, $N_G(H)=H.2$, $G={}^\varepsilon E_6(p)$ for the appropriate $\varepsilon$ and there are exactly $2\cdot\gcd(3,q-\varepsilon)$ many $G$-conjugacy classes.
\end{prop}

\subsection{Cross-characteristic subgroups \texorpdfstring{$\PSL_2(r)$}{PSL(2,r)}}

Let $H\cong \PSL_2(r)$ for $r$ one of $7,8,11,13,17,19,25,27$, and let $L$ denote a Borel subgroup of $H$. We will deal with each case in turn. Proofs of maximality, as with the other groups before, are completed in Section \ref{sec:proofmaximal}.

\subsubsection{$H\cong\PSL_2(7)$}

Here $p=3$ or $p\neq 2,3,7$. If $p=3$ then $H$ is strongly imprimitive by Proposition \ref{prop:allpsl27p=3}. This leaves $p\neq2,3,7$. From \cite[Table 6.68]{litterickmemoir}, there is a unique conspicuous set of composition factors for which $H$ acts fixed point freely on both $M(E_6)$ and $L(E_6)$; it acts on $M(E_6)$ with factors $8,7,6^2$ and on $L(E_6)$ with factors $8^5,7^2,6^2,(3,3^*)^2$. A copy of $H$ with these composition factors on $M(\mbG)$ and $L(\mbG)$ is known to exist as it lies in $\Sp_8$ or $G_2$. (Note that there is a copy of $H$ in $A_5$ with the same action on $M(E_6)$, but it has fixed points on $L(E_6)$.)

\begin{prop} Suppose that $p\neq 2,3,7$ and $H\cong \PSL_2(7)$ acts fixed-point freely on both $M(E_6)$ and $L(E_6)$. Then $H$ is unique up to $\mbG$-conjugacy, $N_\mbG(H)=Z(\mbG)\times \PGL_2(7)$, and $H$ is strongly imprimitive.
\end{prop}

The proof of this proposition has both theoretical and computational aspects. It is among the more challenging of the proofs here.

\medskip\noindent\textbf{Determination of $N_\mbG(H)$}

\medskip\noindent This will become clear once we have proved that $H$ is unique up to $\mbG$-conjugacy. Since $\PGL_2(7)\leq \PSU_3(3)\leq G_2(p)$, and $C_\mbG(H)=Z(H)$, we obtain the result.

\medskip\noindent\textbf{Determination of $L$ up to conjugacy}

\medskip\noindent Let $w\in L$ have order $3$. The action of $L$ on $M(E_6)$ has four each of the two $3$-dimensional simple $kL$-modules, and a single copy of each $1$-dimensional simple module. In particular, $L$ centralizes a line, so lies in $F_4$ or $D_5T_1$. However, if $L\leq D_5T_1$ then $L$ must act as $3+3+3+1$ on $M(D_5)$, but such a module cannot be self-dual as the two $3$-dimensional modules are dual to one another (and the $1$-dimensional module is not self-dual). Thus $L\leq F_4$.

Since $L$ is a supersoluble $p'$-group it lies in $N_\mbG(\mb T)$ by Lemma \ref{lem:normtorus}, where $\mb T$ is a maximal torus of $F_4$, not of $E_6$. Now we note that $w$ lies in the class of elements of order $3$ that act on $M(F_4)$ and $L(F_4)$ with traces $-1$ and $-2$ respectively. Such a class does not have eigenvalue $1$ on a torus, so the centralizer of $w$ on $\mb T$ is finite. Thus by Lemma \ref{lem:algconjweylelements}, all elements of order $3$ in $\mb Tw$ are conjugate, so we may choose one. There are sixteen classes of subgroups $7\rtimes 3$ in $\mb T\rtimes \gen w$, which fall into two $N_\mbG(\mb T)$-classes, one for each rational $F_4$-class of elements of order $7$. Thus $L$ is determined up to conjugacy.

\medskip\noindent\textbf{Determination of $H$ up to conjugacy}

\medskip\noindent This is done in the supplementary materials using the trilinear form method from Section \ref{sec:trilinear}, as we do in the other cases, but this is much more difficult and one cannot simply solve the equations. The main obstacle here is the determination of the centralizer $C_\mbG(L)$, which is more complicated than most other cases. Since this centralizer acts on the solutions, a complicated centralizer yields a complicated set of solutions.

One may determine $C_\mbG(L)$ theoretically, but since we need explicit elements from it, we may as well use the computer. This determines that $C_\mbG(L)$ is the semidirect product of a  rank-$2$ torus by a group $3\times 3$. We also require the centralizer in $\GL_{27}(k)$ of $H$, which is clearly isomorphic to $\GL_2(k)\times (k^\times)^2$, so of dimension $6$. Notice that the centralizer $\mb C$ in $\GL_{27}(k)$ of $L$ has dimension $35$, and is of the form
\[ \GL_4(k)\times \GL_4(k)\times k^\times\times k^\times\times k^\times.\]
One may therefore write an arbitrary element of $\mb C$ in terms of 35 parameters, which in the supplementary materials are $m_1,\dots,m_{16}$ (corresponding to matrix entries for the first subgroup), $n_1,\dots,n_{16}$ (for the second), and $a,b,c$ (which are parameters for the $3$-dimensional torus).

In order to be able to do the computations, we need to take representatives for the orbits on $\mb C$ of $C_\mbG(L)$ on the right and the centralizer of $H$ in $\GL_{27}(k)$ on the left. In fact, we will move the centre of $\GL_{27}(k)$ to the right to make the descriptions easier to understand.

The rest of the proof that is here is to show that the orbit representatives for the left and right actions, that are used in the supplementary materials, are correct.

We have a torus $\mb S$ of rank $3$ acting on the right. The elements of $\mb S$, written as triples $(t_1,t_2,t_3)$ in $k^3$, act as polynomials on the 35 variables in $\mb C$. This action preserves the subset $\{m_9,m_{11},m_{12}\}$ of the variables (in the labelling of the supplementary materials), in the sense that the action of the triple on $m_9$ only involves the $t_i$, $m_9$, $m_{11}$ and $m_{12}$, and similarly for $m_{11}$ and $m_{12}$.

For $t_1,t_2,t_3\in k$, write \[f(t_1,t_2,t_3)=t_1^3+t_2^3+t_3^3-3t_1t_2t_3=(t_1+t_2+t_3)(t_1+\omega t_2+\omega^2 t_3)(t_1+\omega^2t_2+\omega t_3),\]
where $\omega$ is a primitive cube root of unity. With respect to the specific action of the $t_i$ given in the supplementary materials, it is not true that a triple $(t_1,t_2,t_3)\in k^3$ lies in $\mb S$ if and only if $t_1t_2t_3\neq 0$, but rather that $f(t_1,t_2,t_3)\neq 0$.

The action of $(t_1,t_2,t_3)$ on $(m_9,m_{11},m_{12})$ is
\[ (m_9,m_{11},m_{12})^{(t_1,t_2,t_3)}=(m_9t_3+m_{11}t_1+m_{12}t_2,m_9t_2+m_{11}t_3+m_{12}t_1,m_9t_1+m_{11}t_2+m_{12}t_3).\]
(For example, the centre of $\GL_{27}(k)$ is given by $(0,0,t_3)$.) Of course it would be possible to reparametrize $\mb S$ so that $(t_1,t_2,t_3)$ lies in $\mb S$ if and only if $t_1t_2t_3\neq 0$, but then the action on the $m_i$ would be much more complicated.

We claim that there are exactly eight orbits of $\mb S$ on all triples $(m_9,m_{11},m_{12})\in k^3$:
\begin{itemize}
\item One orbit containing all points such that $f(m_9,m_{11},m_{12})\neq 0$;
\item Six orbits of non-zero points such that $f(m_9,m_{11},m_{12})=0$, with representatives $(1,-\omega^i,0)$ for $i=0,1,2$, (these have stabilizer a rank-$1$ torus) and $(1,\omega,\omega^2)$, for $i=0,1,2$ (these have stabilizer a rank-$2$ torus);
\item One orbit consisting of $(0,0,0)$.
\end{itemize}

We prove this now. A fast way to prove this is to determine the stabilizer of a point in each orbit over $\F_q$, then note that
\[ q^3=(q-1)^3+3(q-1)^2+3(q-1)+1.\]
We provide a direct proof, both to show where the orbits come from and also because it then works over any field of characteristic $p$ with a cube root of unity.

For $\zeta$ any cube root of unity, write $f_\zeta(t_1,t_2,t_3)=t_1+\zeta t_2+\zeta^2 t_3$. We first note that, if $(a,b,c)\in k^3$ and $(t_1,t_2,t_3)\in \mb S$, then $f_\zeta(a,b,c)=0$ if and only if $f_\zeta\left((a,b,c)^{(t_1,t_2,t_3)}\right)=0$. Our proposed orbits are exactly the collections of triples $(a,b,c)$ for which some fixed subset of the $f_\zeta$ is zero.

This shows first that each of our eight elements lie in different orbits, so we do not need to prove that. If $\Omega$ denotes the set of all non-zero triples satisfying $f_\zeta=0$ and $f_{\zeta'}=0$ for $\zeta\neq \zeta'$, then since $f_\zeta$ and $f_{\zeta'}$ are distinct linear conditions, $\Omega$ is simply all scalar multiples of a single non-zero element. Furthermore, since we may scale by any non-zero element of $k$, this forms a single orbit.

The orbit containing $(1,0,0)$ is also clear: it is given by
\[ (1,0,0)^{(t_1,t_2,t_3)}=(t_3,t_2,t_1),\]
such that $f(t_1,t_2,t_3)\neq 0$. Thus all triples $(m_9,m_{11},m_{12})$ such that $f(m_9,m_{11},m_{12})\neq 0$ lie in the same orbit.

The only case left is the orbit of $(1,-1,0)$, which is $(t_3-t_1,t_2-t_3,t_1-t_2)$ for all $t_1,t_2,t_3$ such that $f(t_1,t_2,t_3)\neq 0$. If this is $(\alpha,\beta,-\alpha-\beta)$, then $t_3=t_1+\alpha$ and $t_2=t_3+\beta=t_1+\alpha+\beta$. Thus for any $\alpha,\beta$ we have \[(1,-1,0)^{(t_1,t_1+\alpha+\beta,t_1+\alpha)}=(\alpha.\beta,-\alpha-\beta).\]
Furthermore,
\[f(t_1,t_1+\alpha+\beta,t_1+\alpha)=(\alpha^2+\alpha\beta+\beta^2)(t_1+2\alpha+\beta)/3.\]
As $t_1$ can be chosen arbitrarily, this is non-zero unless $\alpha^2+\alpha\beta+\beta^2=0$, i.e., $\alpha=\omega\beta$ or $\alpha=\omega^2\beta$. The element is then $(1,\omega,-1-\omega)=(1,\omega,\omega^2)$ or $(1,\omega^2,\omega)$, which is indeed in a different orbit, as we have already seen.

One of the elements in the centralizer but not inside $\mb S$, denoted $w$ in the supplementary materials, has action on $(m_9,m_{11},m_{12})$ given by
\[ (m_9,m_{11},m_{12})\mapsto (m_9,\omega m_{11},\omega^2m_{12}).\]
This clearly permutes the three orbits with representatives $(1,1,1)$, $(1,\omega,\omega^2)$ and $(1,\omega^2,\omega)$, and with representatives $(1,-1,0)$, $(1,-\omega,0)$ and $(1,-\omega^2,0)$.

This completes the proof of the claim.

\medskip

Having computed the orbits, this means that we have exactly four options for $(m_9,m_{11},m_{12})$. We can use still more of the centralizer of $H$ in $\GL_{27}(k)$. The copy of $\GL_2(k)$ in that centralizer acts as $2\times 2$ matrices on the variables in $\mb C$ labelled $m_1,m_2,m_5,m_6$. Thus if $s_1,s_2,s_3,s_4$ are the variables in a $(2\times 2)$-matrix element
\[ \begin{pmatrix} s_1&s_2\\s_3&s_4\end{pmatrix}\]
of the $\GL_2(k)$ subgroup of the centralizer, we have
\[ (m_1,m_2,m_5,m_6)^{(s_1,s_2,s_3,s_4)}=(m_1s_1+m_5s_2,m_2s_1+m_6s_2,m_1s_3+m_5s_4,m_2s_3+m_6s_4).\]
The orbits of $2\times 2$ matrices under left multiplication by elements of $\GL_2(k)$ are much clearer than for the triples above, and they are of course
\[ \begin{pmatrix}1&0\\0&1\end{pmatrix},\;\;\begin{pmatrix}1&m_2\\0&0\end{pmatrix},\;\;\begin{pmatrix}0&1\\0&0\end{pmatrix}\;\;\begin{pmatrix}0&0\\0&0\end{pmatrix}.\]

For each pair $(m_9,m_{11},m_{12})$ and $(m_1,m_2,m_5,m_6)$ of orbit representatives we can compute the number of solutions to the equations given by imposing the trilinear form as in Section \ref{sec:trilinear}, and find a single solution up to action of the centralizer (which occurs when the representatives are $(1,0,0)$ and $(1,0,0,1)$).

Thus $H$ is unique up to $\mbG$-conjugacy.

\medskip

\noindent \textbf{Strong imprimitivity}

\medskip\noindent To prove this is easy, once we know that $H$ is unique up to $\mbG$-conjugacy. Since $N_\mbG(H)=Z(\mbG)\times \PGL_2(7)$, we see that any Frobenius endomorphism induces an inner automorphism on some conjugate of $H$. Thus $H$ is contained in all groups $G$. In any adjoint group $G_{\mathrm{ad}}$ there exists a unique class of subgroups $\PSp_8(q).2$ (see Tables \ref{tab:e6othermaximals} and \ref{tab:2e6othermaximals}), which contains a copy of $H.2$. Also, there exists a unique class of subgroups $H$, which is necessarily contained in $\PSp_8(q)$. Thus $H$, and indeed $H.2$, is always contained in $\PSp_8(q).2$, and thus $N_{\bar G}(H)$ can never be maximal in $\bar G$. We even obtain strong imprimitivity because the group $\PSp_8(q).2$ is the fixed points of a $C_4$ subgroup.


\subsubsection{$H\cong\PSL_2(8)$}
\label{sec:e6;sl28}

Let $H\cong \PSL_2(8)\cong {}^2\!G_2(3)'$. The only cases are $p=7$ and $p\neq 2,3,7$. For $p=7$ we always have that $H$ is strongly imprimitive (see \cite[Table 6.71]{litterickmemoir}), and for $p\neq 2,3,7$ if $H$ is not strongly imprimitive then $M(E_6){\downarrow_H}$ is the sum of the three non-isomorphic but $\Aut(H)$-conjugate $9$-dimensional modules (see \cite[Table 6.70]{litterickmemoir}).

For $p\neq 2,3,7$, \cite[Theorem 29.3]{aschbacherE6Vun} proves the result, but since that work is unpublished, we will reprove the result here using similar, but not quite the same, means. In \cite{ryba2007}, it is proved that there are two $\mbG$-conjugacy classes of subgroups $H.3$, swapped by the graph automorphism. However, this says nothing about the group $H$ itself, which could have classes that are self-normalizing.

We will not determine, in this section, if $N_G(H)=H$ or $N_G(H)=H.3$ in the simple group $G$. The proof of this is complicated, and will be delayed until Section \ref{sec:sl28}. We will prove maximality in Section \ref{sec:proofmaximal}.

\begin{prop} \label{prop:sl28ine6} Let $H\cong \PSL_2(8)$, suppose that $p\neq 2,3,7$, and that $H$ acts as described above.
\begin{enumerate}
\item We have $N_\mbG(H)=Z(\mbG)\times (H.3)$.
\item There are exactly two $\mbG$-conjugacy classes of subgroups $H$, swapped by the graph automorphism of $\mbG$.
\item The group $H$ embeds in $G=E_6(q)$ if $q\equiv 1,2,4\bmod 7$, and in ${}^2\!E_6(q)$ if $q\equiv 3,5,6\bmod 7$.
\end{enumerate}
\end{prop}

The proof proceeds in stages.

\medskip\noindent\textbf{Determination of $N_\mbG(H)$}

\medskip\noindent We show in the supplementary materials, or see also \cite{aschbacherE6Vun} and \cite{ryba2007}, that in fact $H.3$ embeds in $\mbG$, rather than just $H$. In the course of the proof we will see that all copies of $H$ are $\Aut(\mbG)$-conjugate, and so $N_\mbG(H)=Z(\mbG)\times (H.3)$.

\medskip\noindent\textbf{Determination of $L$ up to conjugacy}

\medskip\noindent The action of $L$ on $M(E_6)$ is the sum of the six non-trivial $1$-dimensional modules and three copies of the $7$-dimensional module. This places $H$ inside $D_5T_1$, then inside $B_3T_2$ since $H$ must stabilize three lines on $M(D_5)$. The copy of $L$ inside $B_3$ acts irreducibly on $M(B_3)$ and lies in $G_2$, whence $L$ lies in $G_2T_2\leq G_2A_2$. From the structure above, we see that $L$ lies diagonally in $G_2A_2$, projecting as $L$ on $G_2$ and $7$ on $A_2$. There is a unique class of $7$s in $A_2$ that act on $S^2(M(A_2))$ as the sum of all six non-trivial modules, and $L$ is unique up to conjugacy in $G_2$. Furthermore, each side is normalized by a $3$, so since $|\Out(L)|=3$ and $|\Out(7)|=6$, there are two diagonal classes of $L$ in $G_2A_2$. Since the graph automorphism of $E_6$ normalizes (but does not act as an inner automorphism on) the $A_2$ factor, it swaps these two classes of $L$.

Thus there are two $\mbG$-classes of subgroups $L$ in $\mbG$, swapped by the graph automorphism. Note that $L$ can be chosen to lie in an irreducible $G_2$ subgroup, and inside $G_2(p)$. (See \cite{kleidman1988} for a list of the maximal subgroups of $G_2(q)$, which includes $2^3\cdot \SL_3(2)\leq G_2(p)$ for $p\geq 5$.)

\medskip\noindent\textbf{Determination of $H$ up to conjugacy}

\medskip\noindent Using the method from Section \ref{sec:trilinear}, we compute the number of copies of $H$ containing $L$ under the centralizer of $L$ in $\GL_{27}(k)$. Since $C_\mbG(L)$ is a $2$-dimensional torus, we must find many overgroups $H$ of a given subgroup $L$. However, we prove in the supplementary materials, by a counting argument, that all overgroups $H$ are $C_\mbG(L)$-conjugate.

\medskip\noindent\textbf{Action of some outer automorphisms}

\medskip\noindent As $H.3$ has no outer automorphisms, an automorphism $\sigma$ of $\mbG$ either fuses the two $\mbG$-classes or it acts as an inner automorphism on them, hence centralizes a representative from each class by Lemma \ref{lem:fixedpoints}.

By \cite{kleidman1988}, $H$ lies in $G_2(q)$ if and only if $x^3-3x+1$ splits over $\F_q$ (i.e., $q\equiv\pm1,\pm3,\pm4\bmod 13$), so in particular $H\leq G_2(p^3)$ always. We see from Tables \ref{tab:e6othermaximals} and \ref{tab:2e6othermaximals} that $H\leq G_2(p^3)\leq E_6(p^3)$ if $p\equiv 1,2,4\bmod 7$, and $H\leq G_2(p^3)\leq {}^2\!E_6(p^3)$ if $p\equiv 3,5,6\bmod 7$. Notice that $F_p$ swaps the two classes if and only if $F_{p^3}$ does, so $F_p$ centralizes $H$ (up to conjugacy) if and only if $p\equiv 1,2,4\bmod 7$. Thus $H\leq E_6(p)$ if $p\equiv 1,2,4\bmod 7$ and $H\leq {}^2\!E_6(p)$ if $p\equiv 3,5,6\bmod 7$.

The remaining question is whether a non-trivial diagonal automorphism of $G$ (if it exists) induces the outer automorphism of $H$, i.e., if $N_G(H)=H.3$ when $H\leq G$. Of course, if $p\not\equiv \varepsilon\bmod 3$ then $G={}^\varepsilon\!E_6(p)$ does not have a non-trivial (outer) diagonal automorphism, so $N_G(H)=H.3$ in this case. We complete this in Section \ref{sec:sl28}.

\subsubsection{$H\cong\PSL_2(11)$}
\label{sec:e6;psl211}
Let $H\cong \PSL_2(11)$. The possible cases are $p=2$, $p=3$, $p=5$ and $p\neq 2,3,5,11$. Contrary to the assertions in \cite{aschbacherE6Vun} and \cite{cohenwales1997}, we will find a Lie primitive copy of $H$ inside $\mbG$ in all characteristics not equal to $5$ and $11$.

We start with a preliminary lemma, which might be of independent interest. It classifies non-abelian subgroups of order $55$ in $\mbG$, for $p\neq 5,11$. It gives a case where Lemma \ref{lem:torusisgood} does not apply, and there are two $\mbG$-classes of subgroup $11\rtimes 5$ where the subgroups $11$ are both $\mbG$-conjugate and the subgroups $5$ are both $\mbG$-conjugate, but there are two non-conjugate ways to put them together.

\begin{prop}\label{prop:frobenius55} Let $p\neq 5,11$. There are three conjugacy classes of non-abelian subgroups $11\rtimes 5$ of order $55$ in $\mbG$, with representatives $L_1,L_2,L_3$. We may choose labellings and representatives so that $L_1$ has a rational element of order $5$, and $L_2$ and $L_3$ both have semirational elements of order $5$ with trace $\zeta+\zeta^{-1}$ on $M(E_6)$. The subgroups of orders $5$ and $11$ in $L_2$ and $L_3$ are $\mbG$-conjugate, but $L_2$ and $L_3$ are not $\mbG$-conjugate.
\end{prop}
\begin{proof} Using Theorem \ref{thm:larsen}, to count the number of $\mbG$-classes of subgroups $11\rtimes 5$ we may assume that $\mbG$ has characteristic $2$. Embed a subgroup $L\cong 11\rtimes 5$ into $N_\mbG(\mb T)$ via Lemma \ref{lem:normtorus}. Tori are complemented in their normalizer in characteristic $2$ (but not in odd characteristics) so we may assume that $N_\mbG(\mb T)=\mb T\rtimes W(E_6)$, which simplifies matters. (We can proceed without this assumption, but things are more complicated.)

Note that $W(E_6)$ contains a unique class of elements of order $5$, so let $w$ be a representative of this class. Any element of order $5$ in $N_\mbG(\mb T)$ is conjugate to an element of the form $tw$ for $t\in \mb T$. We first count the number of classes of elements of order $5$ in $\mb T.\gen w$ whose image modulo $\mb T$ is $w$. For this we need the action of $W(E_6)$ on $\mb T$.

Since $W(E_6)$ has a unique class of elements of order $5$, and $\Sym(5)\leq W(E_6)$, we may choose a subgroup $\Sym(5)$ containing $w$. (There are actually four conjugacy classes of subgroups $\Sym(5)$ in $W(E_6)$.) By choosing the correct subgroup $\Sym(5)$, it acts on $\mb T$ as the direct sum of the permutation representation and the trivial representation, i.e., there is a spanning set $\mb T_i$ of six $1$-dimensional tori with $\mb T_i^w=\mb T_{i+1}$ for $1\leq i\leq 4$, $\mb T_5^w=\mb T_1$ and $w$ centralizes $\mb T_6$.

For explicit computations we will move to finite groups. Choosing $\mb T$ to be maximally split with respect to a Frobenius endomorphism $F_q$, any element of order $5$ lies in some $N_\mbG(\mb T)^{F_q}$ for some $q$. In the finite group $T\rtimes \gen w\cong (q-1)^6\rtimes \gen w$, we count elements of order $5$ not in the homocyclic subgroup $(q-1)^6$. By the Schur--Zassenhaus theorem, to count conjugacy classes we may replace $T\rtimes \gen w$ with a group $X\rtimes \gen w$, where $X$ is the Sylow $5$-subgroup of $T$, and we do this.

Let $t\in X$. Notice that $tw$ has order $5$ if and only the product of $t^{w^i}$ for $i=0,\dots,4$ is $1$. If $x_1,\dots,x_6$ is a basis for $X$, then we may choose that $x_i$ so that $x_i^w=x_{i+1}$ for $1\leq i\leq 4$, $x_5^w=x_1$ and $x_6^w=x_6$.

Write $t=x_1^{a_1}x_2^{a_2}\ldots x_6^{a_6}$, and let $5^n$ be the order of each $x_i$. We have
\[ \prod_{i=0}^4 t^{w^i}=(x_1x_2x_3x_4x_5)^{a_1+a_2+a_3+a_4+a_5}\cdot x_6^{5a_6}.\]
Thus $5^n\mid (a_1+a_2+a_3+a_4+a_5)$ and $5^{n-1}\mid a_6$. This yields $5^{4n+1}$ options for the $a_i$, so there are that number of elements $t$ such that $o(tw)=5$. On the other hand, $|C_X(tw)|=|C_X(w)|=5^{2n}$, and $|X|=5^{6n}$. Hence there are exactly five conjugacy classes of elements of order $5$ in $\mb T\gen w$ with image $w$ modulo $\mb T$. (They are $x_6^iw$ for the appropriate $i$.)

The normalizer in $W(E_6)$ of $\gen w$ has order $40=8\cdot 5$, and a Sylow $2$-subgroup $Q$ of this acts on the five conjugacy classes of complements above. Of course it stabilizes $\gen w$. We check with a computer (by constructing them in a subgroup $5^6\rtimes 5$ of $\mbG$) that the four classes other than $w$ have Brauer character $\zeta_5+\zeta_5^{-1}$ ($\zeta_5$ a primitive fifth root of unity) on $M(E_6)$, so cannot be rational. Thus the stabilizer in $Q$ of any conjugacy class other than the one containing $\gen w$ has order $2$, and $Q$ acts transitively on the classes. Thus in $N_\mbG(\mb T)$ there are exactly two classes of subgroups of order $5$ outside $\mb T$, $\gen w$ and (say) $\gen{w_1}$. (One may also simple check in $5^6\rtimes W(E_6)$ that there are two classes of subgroups of order $5$ outside the subgroup $5^6$.)

To produce a subgroup $L$, we need to consider the action of $\gen w$ and $\gen{w_1}$ on the subgroup $11^6$ of $\mb T$. We know from above that, as a $6$-dimensional module over $\F_{11}$, $\gen w$ and $\gen{w_1}$ act the same, and as a sum of two trivial modules and one copy each of the four non-trivial modules. For $\gen w$, the normalizer of $\gen w$ acts to permute all non-trivial modules, and hence there is up to $N_\mbG(\mb T)$-conjugacy a unique subgroup $11\rtimes 5$ containing $w$. For $\gen{w_1}$, the normalizer has order $2$, so the four non-trivial modules are swapped in pairs, and we obtain two non-conjugate subgroups $L_2$ and $L_3$ containing $w_1$.

To see that $L_2$ and $L_3$ are non-conjugate in $\mbG$, not just in $N_\mbG(\mb T)$, we actually note that they are non-conjugate even in $\GL_{27}(k)$, because the modules are not $\Aut(L_i)$-conjugate. Let $1_1,1_2,1_2^*,1_3,1_3^*$ denote the five $1$-dimensional $kL$-modules. The module action of $L_2$ on $M(E_6)$ is
\[ 5^{\oplus 3}\oplus (5^*)^{\oplus 2}\oplus 1_2\oplus 1_2^*,\]
and this is not $\Aut(L_2)$-conjugate to the sum of those $5$-dimensional modules and $1_3\oplus 1_3^*$. Thus $L_2$ and $L_3$ are non-conjugate in $\mbG$. The three subgroups $L_1$, $L_2$ and $L_3$ are not $\GL_{27}(k)$-conjugate, so are certainly not $\mbG$-conjugate. 
\end{proof}

Now we come back to $H\cong \PSL_2(11)$. In \cite[Tables 6.73--6.76]{litterickmemoir}, we see that if $p=5$ then $H$ is strongly imprimitive, so we assume that $p\neq 5,11$. If $p=2$ then there are two rows of \cite[Table 6.76]{litterickmemoir} labelled `\textbf P'. One of these, row 4, has factors $10,5^2,5^*,1^2$ on $M(E_6)$. Both $5$ and $5^*$ have non-zero $1$-cohomology (but $10$ has zero $1$-cohomology) \cite[Table A.4]{litterickmemoir}, so the pressure is $1$. However, we show that $H$ nevertheless must stabilize a line or hyperplane on $M(E_6)$, and hence is strongly imprimitive. As the pressure is $1$, we may assume that the socle is $5$ or $5^*$: a computer check on Magma shows that the largest submodule of $P(5)$ with composition factors from the set $\{10,5,5^*,1\}$ has structure
\[ 5/1,5^*/5.\]
Since there is only one trivial module in this module, we cannot produce a module with factors $10,5^2,5^*,1^2$ with two trivial composition factors and no trivial submodule or quotient. Thus this row may be excluded, and there is a single row to consider for each of $p=2$, $p=3$ and $p\nmid |H|$.

If $p=2$ or $p\nmid |H|$ then the dimensions of the composition factors of $M(E_6){\downarrow_H}$ are $5$, $10$ and $12$. There are two (dual) possible simple modules of dimension $5$ and two of dimension $12$, yielding four possibilities for $M(E_6){\downarrow_H}$ (two up to outer automorphism of $H$, as it swaps the two $12$-dimensional modules). If $p=3$ then again there are two possible simple modules of dimension $12$ and two (dual) of dimension $5$, and the factors of $M(E_6){\downarrow_H}$ are $12,5,5,5^*$. Thus for $p\neq 5,11$, we obtain four possibilities for the composition factors of $M(E_6){\downarrow_H}$, two up to $\Aut(H)$-conjugacy.

For $p\neq 3$, $M(E_6){\downarrow_H}$ must be semisimple. For $p=3$, the only action compatible with the unipotent action from \cite[Table 5]{lawther1995} is $(5/5^*/5)\oplus 12$, so this must be $M(E_6){\downarrow_H}$.

For this action, we have the following result, delaying the proof of maximality until Section \ref{sec:proofmaximal}.

\begin{prop} Let $H\cong \PSL_2(11)$ and let $p\neq 5,11$. Suppose that $H$ acts as described above.
\begin{enumerate}
\item The group $H$ is Lie primitive and $N_\mbG(H)=Z(\mbG)\times H$.
\item There are two $\mbG$-classes of subgroups $H$, both normalized by the graph automorphism, so that $H.2\cong \PGL_2(11)$ embeds in $\mbG.2$.
\item The group $H$ always embeds in $E_6(p^2)$, and embeds in $E_6(p)$ if and only if $f_2(x)=x^2+x+3$ and $f_3(x)=x^2+x-1$ split over $\F_p$ (i.e., $p\equiv \pm 1\bmod 5$ and $p\equiv 1,3,4,5,9\bmod 11$). If $f_3(x)$ splits but $f_2(x)$ does not (i.e., $p\equiv \pm 1\bmod 5$ and $p\equiv 2,6,7,8,10\bmod 11$) then $H$ embeds in ${}^2\!E_6(p)$.
\item The subgroup $H$ is maximal in $G$ if and only if $G={}^\varepsilon\!E_6(p^a)$ is the minimal group into which $H$ embeds, as given in the previous part. There are $2\cdot \gcd(3,p-\varepsilon)$ classes of subgroups $H$. If $\bar G\neq G$ and $N_{\bar G}(H)$ is maximal in $\bar G$ then $\bar G=G.2$, $H$ is maximal in $G$, $\bar G$ induces a graph automorphism on $G$, and there are two $\bar G$-classes of subgroups $H.2$.
\end{enumerate}
\end{prop}

The proof proceeds in stages.

\medskip\noindent\textbf{Determination of $N_\mbG(H)$}

\medskip\noindent Since $M(E_6){\downarrow_H}$ is not stable under the outer automorphism of $H$, we must have $N_\mbG(H)=Z(\mbG)\times H$.

\medskip\noindent\textbf{Determination of $L$ up to conjugacy}

\medskip\noindent The Brauer character of $H$ implies that an element of order $5$ acts on $M(E_6)$ with character $\zeta_5+\zeta_5^{-1}$. Hence there are two options for $L$ up to $\mbG$-conjugacy, by Proposition \ref{prop:frobenius55}.

\medskip

\noindent\textbf{Determination of $H$ up to conjugacy}

\medskip \noindent In the supplementary materials we pick one of the two classes of subgroups $L$, and prove that for all primes, all subgroups $H$ containing $L$ are $C_\mbG(L)$-conjugate. Thus there are exactly two $\mbG$-conjugacy classes of subgroups $H$.

(Note that we only construct one of the classes of subgroups $H$ in the supplementary materials. We will prove the other case exists in the next part.)

\medskip\noindent\textbf{Action of outer automorphisms}

\medskip\noindent Note that an outer automorphism of $H$ swaps the $kH$-modules $5$ and $5^*$ and fixes $12_1$ and $12_2$. Thus a graph automorphism of $\mbG$ cannot fuse classes, and cannot centralize $H$ as $H$ is not contained in $F_4$ or $C_4$, the centralizers of graph automorphisms of $\mbG$. Thus it normalizes both classes, and $\PGL_2(11)$ embeds in $\mbG.2$.

The diagonal automorphisms of course merge classes in the finite group by Corollary \ref{cor:diagaut}, so it remains to consider the field automorphisms. Note that the module $5$ exists over $\F_q$ if and only if the polynomial $f_2(x)$ splits (and this irrationality appears for an element of order $11$), and the module $12_1$ exists over $\F_q$ if and only if the polynomial $f_3(x)$ splits (and this irrationality appears for an element of order $5$). Of course, $F_{p^2}$ fixes all of these irrationalities, so $F_{p^2}$ centralizes $H$. Whether $F_p$ centralizes $H$, normalizes $H$, or fuses the two $\mbG$-classes, depends on the splitting of $f_2(x)$ and $f_3(x)$. This also proves that the second class over $\mbG$ exists, since it must do for certain primes (as the field automorphism does not stabilize a class) and hence for all primes by Theorem \ref{thm:serre}.

The map $F_p$ swaps the two $\mbG$-classes if and only if $f_3(x)$ does not split over $\F_p$, so we assume that $f_3(x)$ splits over $\F_p$, and therefore $F_p$ normalizes $H$. If $f_2(x)$ does not split over $\F_p$ then $H$ cannot embed in $E_6(p)$, so that $F_p$ acts as an outer automorphism on $H$ and the product with the graph centralizes $H$. Thus $H\leq {}^2\!E_6(p)$. Conversely if $f_2(x)$ splits over $\F_p$ then $F_p$ centralizes all four representations of $H$, so cannot act as the outer automorphism. Thus $H\leq E_6(p)$. In both cases for $G$, $H.2$ embeds in $G.2$.

This completes the proof.

\begin{rem} The error in \cite{aschbacherE6Vun} appears in the proof of (17.1), where it is assumed that (in that notation) $h$ is equal to $k$. Since $C_{\mb T}(w)$ is positive-dimensional (where $w$ is of order $5$ in $W(E_6)$) there can be, and are, multiple classes of complements. As no proof is given in \cite{cohenwales1997}, the source of the error there cannot be found.
\end{rem}

\subsubsection{$H\cong\PSL_2(13)$}

Let $H\cong \PSL_2(13)$, and let $L$ denote the subgroup $13\rtimes 6$ of $H$. Let $L_1$ denote the subgroup $13\rtimes 3$ of index $2$ in $L$. The cases to consider are $p=2$, $p=3$, $p=7$ and $p\neq 2,3,7,13$. If $p=7$ then $H$ is strongly imprimitive by \cite[Theorem 1]{litterickmemoir}, so we ignore this case. From \cite[Tables 6.77, 6.79 and 6.80]{litterickmemoir}, we see that either $H$ is strongly imprimitive or the composition factors of $M(E_6){\downarrow_H}$ are: $13+14$ for $p\neq 2,3,7,13$ (a unique such module with rational character values); $13^2,1$ if $p=3$ (and such a module must be uniserial $13/1/13$, else $H$ stabilizes a line on $M(E_6)$); $14,6_1,6_2,1$ if $p=2$. In the final case, the $14$ splits off as a summand, since it lies in a different block, and if $H$ does not stabilize a line or hyperplane on $M(E_6)$ the remaining summand must be (up to automorphism) $6_1/1/6_2$. So in all cases the module $M(E_6){\downarrow_H}$ is determined completely, up to automorphism for $p=2$.

The character of $M(E_6){\downarrow_H}$ is always rational-valued, so there is no obstacle to $H$ lying in $G=E_6(p)$ for any prime $p$.

We will prove the following, delaying the proof of maximality until Section \ref{sec:proofmaximal}.

\begin{prop} Let $p\neq 7,13$ and $H\cong \PSL_2(13)$ be a subgroup of $\mbG$, and let $\sigma$ be a Frobenius endomorphism of $\mbG$ with $G=\mbG^\sigma$. Suppose that $H$ acts as described above.
\begin{enumerate}
\item There is a unique $\mbG$-conjugacy class of subgroups isomorphic to $H$, $H$ is contained in a maximal $G_2$ subgroup of $\mbG$, and $N_\mbG(H)=Z(\mbG)\times H$. The group $H.2\cong \PGL_2(13)$ embeds in $\mbG.2$ and is Lie primitive.
\item $H$ embeds into exactly one of $E_6(p)$ and ${}^2\!E_6(p)$. It embeds in ${}^2\!E_6(p)$ if and only if exactly one of $p\equiv1,2,4\bmod 7$ and $p\equiv \pm1,\pm3,\pm4\bmod 13$ holds.
\item If $p=2$ then $N_{\bar G}(H)$ is never maximal in an almost simple group $\bar G$. If $p$ is odd, then $H$ is maximal in $G={}^\varepsilon\!E_6(p)$ if and only if $p\equiv\pm2,\pm5,\pm6 \bmod 13$. If $p\equiv \pm1,\pm3,\pm4\bmod 13$ then $H.2$ is a novelty maximal subgroup of ${}^\varepsilon\!E_6(p).2$. There are $\gcd(3,p-\varepsilon)$ classes of subgroups $N_G(H)$ and one class in $N_{\bar G}(H)$ for $\bar G=G.2$.
\end{enumerate}
\end{prop}

The proof proceeds in stages.

\medskip\noindent\textbf{Determination of $N_\mbG(H)$}

\medskip\noindent Although $M(E_6){\downarrow_H}$ is stable under the outer automorphism of $H$, any extension of it to $\PGL_2(13)$ has trace $\pm 1$ for an outer involution. Thus $N_\mbG(H)=Z(\mbG)\times H$ for $p\neq 2$. For $p=2$, as the $6$-dimensional simple modules are swapped by the outer automorphism, $M(E_6){\downarrow_H}$ is not stable under the outer automorphism. Thus $N_\mbG(H)=Z(\mbG)\times H$ in all cases.

\medskip\noindent\textbf{Determination of $L$ up to conjugacy}

\medskip\noindent By \cite[Lemma 2.1]{cohenwales1993}, $L_1$ is unique up to $\mbG$-conjugacy for $p\neq 3,13$, and $L$ is unique up to $\mbG$-conjugacy for $p\neq 2,3,13$. Thus we first assume that $p=3$, where we will show that $L_1$ is again unique up to conjugacy. As an element $x$ of order $13$ in $L$ is regular, $C^\circ_\mbG(x)$ is a torus. Since any element normalizing $\gen x$ normalizes $C^\circ_\mbG(x)$, $L\leq N_\mbG(\mb T)$. Finally, we need to understand elements of order $3$ in $W(E_6)$. There are three conjugacy classes, each belonging to a different unipotent class of $\mbG$ (as demonstrated by their Jordan block action on $M(E_6)$ and $L(E_6)$). Thus $w\in L_1$ of order $3$ is determined, and subgroups $13\rtimes 3$ of $N_\mbG(\mb T)$ containing $w$ are $N_\mbG(\mb T)$-conjugate. Thus $L_1$ is unique up to conjugacy in characteristic $3$.

To obtain uniqueness of $L$, we just need to know that $C_\mbG(L_1)$ has odd order, for then $13\rtimes 6$ must be uniquely determined in $N_\mbG(L_1)$. For $p=2$, we have $C_\mbG(x)=\mb T$ since $\mbG$ is simply connected, so $C_\mbG(L_1)$ is a group of odd order, and thus $L$ is unique up to conjugacy. For $p=3$, we need to check that $L_1$ does not centralize an involution in $\mb T$, and this is also the case. Thus $L$ is unique up to $\mbG$-conjugacy in all characteristics.

\medskip

\noindent\textbf{Determination of $H$ up to conjugacy}

\medskip\noindent Using the method from Section \ref{sec:trilinear}, we compute the number of copies of $H$ containing $L$ under the centralizer of $L$ in $\GL_{26}(k)$. If $p\neq 2,3,7,13$ then in \cite[Theorem 3.1]{cohenwales1993} it was proved that there are exactly two subgroups $H$ in $\mbG$ containing a fixed $L$, and these can be conjugated by an element of $N_\mbG(L)$ (in which $N_H(L)\times Z(\mbG)$ has index $2$). Thus $H$ is unique up to $\mbG$-conjugacy. In fact, $H$ is contained in a maximal $G_2$ subgroup, as proved in \cite{cohenwales1993}.

For $p=2,3$, in the supplementary materials we also find exactly two subgroups $H$ containing $L$, swapped by $N_\mbG(L)$. In these cases we also find a copy of $G_2$ containing $H$.

\medskip\noindent\textbf{Action of outer automorphisms}

\medskip\noindent The graph automorphism cannot fuse classes, and cannot centralize $H$ as $H$ is not contained in $F_4$ or $C_4$. Thus it normalizes $H$, and $\PGL_2(13)$ embeds in $\mbG.2$.

Again, $F_p$ cannot fuse classes, so it either centralizes $H$---in which case $H\leq E_6(p)$---or it acts as an outer automorphism on $H$---in which $H\leq {}^2\!E_6(p)$. In both cases $H$ lies in $E_6(p^2)$, in fact in $G_2(p^2)\leq E_6(p^2)$.

Since $H\leq G_2(p)$ if and only if $p\equiv \pm1,\pm3,\pm4\bmod 13$ (see \cite{kleidman1988} or \cite[Table 8.41]{bhrd}), if $p$ is congruent to one of those numbers then $H\leq E_6(p)$ if and only if $G_2(p)\leq E_6(p)$, i.e., $p\equiv 1,2,4\bmod 7$. Thus we may assume that $p\equiv \pm2,\pm5,\pm6\bmod 13$ and thus the field automorphism of $G_2(p^2)$ acts as an outer automorphism on $H$.

If $p\equiv 1,2,4\bmod 7$ then the field automorphism $F_p$ does not swap the two classes of $G_2$ subgroup, so induces a field automorphism on them. Thus $G_2(p^2).2$ is a subgroup of $E_6(p^2).2$ (field automorphism). Thus the field automorphism $F_p$ of $E_6(p^2)$ induces an outer automorphism on $H$. Hence $H\leq {}^2\!E_6(p)$.

On the other hand, if $p\equiv 3,5,6\bmod 7$ then the field-graph automorphism $F_p\cdot \tau$ of $E_6$ normalizes the two $G_2$ classes, and so $G_2(p^2).2\leq E_6(p^2).2$ (field-graph automorphism). As in the previous paragraph, we obtain $H\leq E_6(p)$.

Thus $H\leq E_6(p)$ if and only if $p\equiv \pm1,\pm3,\pm4\bmod 13$ and $p\equiv 1,2,4\bmod 7$, or $p\equiv \pm2,\pm5,\pm6\bmod 13$ and $p\equiv3,5,6\bmod 7$. On the other hand, $H\leq {}^2\!E_6(p)$ if and only if $p\equiv \pm1,\pm3,\pm4\bmod 13$ and $p\equiv 3,5,6\bmod 7$, or $p\equiv \pm2,\pm5,\pm6\bmod 13$ and $p\equiv1,2,4\bmod 7$.

Note that $H.2$ is a novelty maximal subgroup of ${}^\varepsilon\!E_6(p).2$ if and only if $H\leq G_2(p)$, so $p\equiv\pm1,\pm3,\pm4\bmod 13$, and in the other cases $H$ is a maximal subgroup of $G$. (Proof of maximality is in Section \ref{sec:proofmaximal}.) The only case where this is not true is for $p=2$, since $H$ is contained in $G_2(3)\leq \Omega_7(3)\leq {}^2\!E_6(2)$. To see that $H$ cannot form a novelty maximal subgroup, from \cite[Table 8.42]{bhrd} there is a single class of subgroups $\PSL_2(13)$ in $G_2(3)$, and that the graph automorphism of $G_2(3)$ normalizes $\PSL_2(13)$. Thus $H$ cannot form either a type I or a type II novelty. (This tallies with \cite[p.191]{atlas}, where $H$ does not appear as a novelty maximal subgroup of ${}^2\!E_6(2)$.)

%

\subsubsection{$H\cong\PSL_2(17)$}

Let $H\cong \PSL_2(17)$. The possible primes here are $p=2$, $p=3$ and $p\neq 2,3,17$. This case appears to have been erroneously excluded in both \cite[Section 21]{aschbacherE6Vun} and \cite{cohenwales1997}, particularly as in both cases it is proved that a subgroup acting with composition factors of degrees $9$ and $18$ on $M(E_6)$ does not exist. One exists in $C_4$. Because of this we shall be particularly careful. If $p=2$ then $H$ is always strongly imprimitive by \cite{litterickmemoir}.

Let $L\cong 17\rtimes 8$ denote a Borel subgroup of $H$. Notice that $M(E_6){\downarrow_L}$ is self-dual, but not stable under the outer automorphism of $L$ (that forms $17\rtimes 16$). Hence $L$ is centralized by a graph automorphism, thus lies inside $F_4$ or $C_4$. The former is impossible since $L$ does not centralize a line on $M(E_6)$, so $L\leq C_4$. Then $L$ is unique up to conjugacy in $\mbG$ by \cite[Lemma 1.8.10(ii)]{bhrd}.

We see in the supplementary materials that in characteristic $p\neq 2,17$ there is a unique copy of $H$ containing a given $L$, up to $C_\mbG(L)$-conjugacy. Since there is a copy of $H$ inside $C_4$ we must have that $H$ is Lie imprimitive.

To prove strong imprimitivity, we use the Lie algebra. In characteristic $p\neq 2,3,17$, the composition factors of $H$ on $L(E_6)$ have dimensions $9$, $16$, $17$, $18$ and $18$, and the dimensions of the factors of $C_4$ on $L(E_6)$ are $36$ and $42$. Thus $H$ is strongly imprimitive by Proposition \ref{prop:intorbit}, since clearly the $36$ is the sum of the two $18$s. In characteristic $3$, the factors for $H$ are $1,9,16,16,18,18$ (as $L(E_6)$ has a trivial factor for $p=3$), and the factors for $C_4$ are $1,36,41$. Again, the $36$ is the sum of the two $18$s, and so $H$ is strongly imprimitive again.

\begin{prop} Any subgroup $H\cong \PSL_2(17)$ in $\mbG$ is strongly imprimitive.
\end{prop}

\subsubsection{$H\cong\PSL_2(19)$}

Let $H\cong \PSL_2(19)$, so we may assume $p\neq 19$. The cases to consider are $p=2$, $p=3$, $p=5$ and $p\neq 2,3,5,19$. For $p=2$ there is a copy of $\PSL_2(19)$ inside $J_3$, which lies in $E_6(4)$, so we might expect to see a different answer for $p=2$ to other primes.

Consulting \cite[Tables 6.84--6.87]{litterickmemoir}, the composition factors of $M(E_6){\downarrow_H}$ are of dimensions $9$ and $18$ unless $p=5$, in which case they are $9,9,9^*$. Furthermore, there are no extensions between $9$ and $18$ so the module is always semisimple for $p\neq 5$. For $p=5$ the $9$-dimensional modules have no self-extensions, but there is a non-split extension between $9$ and $9^*$. An element $u\in H$ of order $5$ acts on each $9$ with Jordan blocks $5,4$, and on $9/9^*$ with blocks $5^3,3$, and on $9/9^*/9$ with Jordan blocks $5^5,2$. We see from \cite[Table 5]{lawther1995} that the only option is $9/9^*/9$.

There are two dual irreducible $kH$-modules of dimension $9$, permuted by the outer automorphism of $H$, and for $p\neq 5$ two $18$-dimensional $kH$-modules (with trace of involution $+2$) that are stabilized by the outer automorphism of $H$. Thus one obtains four conjugacy classes of representations of $H$ on $M(E_6)$, two up to duality for $p\neq 5$, and two classes of representations and one up to duality for $p=5$.

We will delay the proof of maximality in the next proposition until Section \ref{sec:proofmaximal}.

\begin{prop} Let $H\cong \PSL_2(19)$ be a subgroup of $\mbG$ and let $\sigma$ be a Frobenius endomorphism of $\mbG$ with $G=\mbG^\sigma$. Suppose that $p\neq 19$, and that $H$ acts as described above.
\begin{enumerate}
\item 
If $p\neq 5$ there are exactly two $\mbG$-conjugacy classes of subgroups isomorphic to $H$, and if $p=5$ there is a single $\mbG$-conjugacy class of subgroups isomorphic to $H$. In both cases, $H$ is always Lie primitive, and $N_\mbG(H)=Z(\mbG)\times H$. Furthermore, $H$ extends to $H.2\cong \PGL_2(19)$ in $\mbG.2$.

\item There is an embedding of $H$ into $G=E_6(q)$ if and only if the polynomial $(x^2-x-1)(x^2+x+5)$ splits over $\F_q$ (i.e., $q\equiv 0,\pm 1\bmod 5$ and $q\equiv 1,4,5,6,7,9,11,16,17\bmod 19$). If $H$ embeds in $G$ then there are $2\cdot\gcd(3,q-1)$ conjugacy classes. Furthermore, $H$ embeds in $G={}^2\!E_6(p)$ if and only if $f_3(x)=x^2-x-1$ splits but $f_4(x)=x^2+x+5$ does not split (i.e., $p\equiv 0,\pm 1\bmod 5$ and $p\equiv 2,3,8,10,12,13,14,15,18\bmod 19$). In this case there are $2\cdot \gcd(3,p+1)$ conjugacy classes.

\item Suppose that $\bar G$ is an almost simple group such that $N_{\bar G}(H)$ is maximal in $\bar G$. Suppose that $p$ is odd. If $G=E_6(q)$ then either $q=p^2$ and $p\equiv \pm2\bmod 5$, or $q=p$ and both $f_3(x)$ and $f_4(x)$ split over $\F_p$. In both cases, either $\bar G=G$ or $\bar G=G.2$, where the outer automorphism of $G$ is the graph automorphism.

If $G={}^2\!E_6(q)$, then $q=p$ and $f_3(x)$ splits over $\F_p$ but $f_4(x)$ does not. The group $\bar G$ is either $G$ or $G.2$.

If $p=2$ then $\bar G=E_6(4).2$ with the extension being the graph automorphism, and $N_{\bar G}(H)\cong \PGL_2(19)$ is a novelty maximal subgroup. There are two such classes, swapped by the subgroup consisting of field automorphisms.
\end{enumerate}
\end{prop}

The proof proceeds in stages.

\medskip\noindent\textbf{Determination of $N_\mbG(H)$}

\medskip\noindent Since $M(E_6){\downarrow_H}$ is not stable under the outer automorphism of $H$, we must have $N_\mbG(H)=Z(\mbG)\times H$.

\medskip\noindent\textbf{Determination of $L$ up to conjugacy}

\medskip\noindent Since $L_0\leq L$ of order $19$ is regular, the normalizer of $L_0$ is contained in $N_\mbG(\mb T)$. Hence $L$ is contained in $N_\mbG(\mb T)$ in all cases.

Let $w\in L$ have order $9$. An easy computer check shows that the centralizer of $w$ on $\mb T$ is finite, so all elements of $\mb Tw$ are $\mb T$-conjugate by Lemma \ref{lem:algconjweylelements}. (See also \cite[Proof of 6.7]{cohenwales1997}.) If $p=3$ then $\mb T$ possesses no $3$-elements, so as there is a single class of elements of order $9$ in $W(E_6)$, we see that $L$ is unique up to conjugacy in $N_\mbG(\mb T)$. If $p\neq 3,19$ then $L$ is a supersoluble $p'$-group, and is unique up to $\mbG$-conjugacy. This can be seen since the $19$ is regular and the centralizer of the $9$ on $\mb T$ is finite (and we may apply Lemma \ref{lem:torusisgood}), and also appears in \cite{cohenwales1997}. In fact, one can see that the centralizer of $L$ in $\mbG$ is $Z(\mbG)$ so all copies of $L$ in any finite $G$ are conjugate in $G_{\mathrm{ad}}$ by Corollary \ref{cor:centsubgroups}.

\medskip

\noindent\textbf{Determination of $H$ up to conjugacy}

\medskip\noindent Using the method from Section \ref{sec:trilinear}, we compute the number of copies of $H$ containing $L$ under the centralizer of $L$ in $\GL_{27}(k)$. We find exactly one for each representation of $H$, so two copies of $H$ containing $L$ for $p\neq 5,19$, and one for $p=5$. This is done in the supplementary materials for $p=2$, $p=3$, $p=5$ and $p\neq 2,3,5,19$.

\medskip\noindent\textbf{Action of outer automorphisms}

\medskip\noindent Note that the graph automorphism must act as the outer $H$-automorphism on each $\mbG$-class, since $M(E_6)^*{\downarrow_H}$ is the image of $M(E_6){\downarrow_H}$ under the automorphism. The character of $M(E_6){\downarrow_H}$ has two irrationalities: elements of order $5$ require the polynomial $f_3(x)=x^2-x-1$ to split; elements of order $19$ require the polynomial $f_4(x)=x^2+x+5$ to split. So in order for $M(E_6){\downarrow_H}$ to exist over $\F_q$, the polynomial $(x^2-x-1)(x^2+x+5)$ must split over $\F_q$.

Since $L$ is unique up to $\mbG$-conjugacy, we see that $F_{p^2}$ must centralize $L$. Since one cannot normalize $H$ while centralizing $L$ without centralizing $H$, and $F_{p^2}$ cannot fuse $\mbG$-classes of subgroups, we see that $F_{p^2}$ also centralizes $H$. Since the graph automorphism acts as the outer automorphism, we have that $F_p$ normalizes $H$ if and only if its product with the graph automorphism does. In this case, exactly one centralizes $H$. Diagonal automorphisms of $G$ fuse classes by Corollary \ref{cor:diagaut}.

Note that the two $18$-dimensional modules are not swapped by the outer automorphism of $H$. Thus if the character values of these modules do not lie in $\F_p$ then $F_p$ swaps the two characters, hence swaps the two $\mbG$-classes of subgroups $H$. Thus we may suppose that $f_3(x)$ (which is the minimal polynomial for an irrationality in the character of an $18$) splits over $\F_p$. If $f_4(x)$ splits over $\F_p$ then the two dual $9$-dimensional modules, which are swapped by the outer automorphism of $H$, are stabilized by $F_p$, hence $F_p$ centralizes $H$ and the product of $F_p$ with the graph automorphism does not centralize $H$. The converse holds: if $f_4(x)$ does not split over $\F_p$ then $F_p$ normalizes but doesn't centralize $H$, and the product with the graph automorphism does centralize $H$. This yields the proposition above.

\medskip

We just need to check that $H\leq J\cong J_3$ for $p=2$. First, $H\leq E_6(4)$ as $(x^2+1)(x^2+x+5)$ splits over $\F_4$, but $H$ does not lie in either $E_6(2)$ or ${}^2\!E_6(2)$. Note that there are two subgroups $3\cdot J_3$ contained in $G$, and since $\PSL_2(19)\leq J_3$ we cannot have that $H$ is maximal in $G$. Note that $H$ is normalized by a graph automorphism, and so is $J$, so any potential novelty maximal subgroup lies in this group $\bar G=G.2$. The outer automorphism of $J_3$ fuses the two classes of $\PSL_2(19)$ subgroups, and so $\PGL_2(19)$ cannot be contained in $N_{\bar G}(H)$. This completes the proof.

\subsubsection{$H\cong\PSL_2(25)$ or $H\cong\PSL_2(27)$}

In both of these cases, $H$ was proved to be strongly imprimitive for all primes in \cite{litterickmemoir}, so we have the following result.

\begin{prop} Let $H$ be isomorphic to either $\PSL_2(25)$ or $\PSL_2(27)$, and suppose that $H$ is a subgroup of $\mbG$. Then $H$ is strongly imprimitive.
\end{prop}

\section{\texorpdfstring{$\PSL_2(8)$}{PSL(2,8)}}\label{sec:sl28}

This section considers the question of whether, in ${}^\varepsilon\!E_6(q)$, the subgroup $\PSL_2(8)$ is self-normalizing, where $q$ is a power of a prime $p\neq 2,3,7$. It turns out that this is a subtle question, depending on whether a particular number is a cube in $\F_q$.

The proof here follows closely that of \cite[Section 29]{aschbacherE6Vun}. Indeed, it relies heavily upon Aschbacher's ideas, particularly when constructing the $E_6$-forms on the $27$-dimensional module for $H=\PSL_2(8).3$. However, our proof is not identical to that of \cite{aschbacherE6Vun}; the author had difficulty with the final paragraph of the proof of \cite[(29.18)(4)]{aschbacherE6Vun}, so we proceed in a different way, more in keeping with the philosophy of this paper. This means that we determine the maximal subgroups of the finite groups by considering the action of field automorphisms on the subgroups of the algebraic group, rather than work directly in the finite group. Because we work with $\PSL_2(8).3$ rather than $\PSL_2(8)$, the calculations are significantly easier.

Let $p\neq 2,3,7$ be a prime, and let $\mbG$ be the simply connected $E_6$ in characteristic $p$. We begin with an easy lemma that gives some elementary facts about $H$ in $\mbG$.

\begin{lem}\label{lem:sl28prelim} Let $\sigma$ be a Frobenius endomorphism on $\mbG$ such that $\sigma$ centralizes $H'\cong \PSL_2(8)$.
\begin{enumerate}
\item There are exactly six $\mbG$-conjugacy classes of subgroups $H$ acting irreducibly on $M(E_6)$.
\item We have that $\sigma$ centralizes a subgroup of $N_\mbG(H')$ isomorphic to $H$ (but not necessarily $H$ itself) if and only if $\sigma^2$ centralizes $H$.
\item If $\sigma$ inverts $Z(\mbG)$ (i.e., $\mbG^\sigma$ has no diagonal automorphisms) then $\sigma$ centralizes a subgroup of $N_\mbG(H')$ isomorphic to $H$.
\item It is always true that $\sigma^3$ centralizes a subgroup of $N_\mbG(H')$ isomorphic to $H$.
\end{enumerate}
\end{lem}
\begin{proof} From Proposition \ref{prop:sl28ine6} we see that there are two $\mbG$-conjugacy classes of the subgroup $H'$. Certainly $H$ lies inside $N_\mbG(H')=Z(\mbG)\times H$ by Proposition \ref{prop:sl28ine6}. Inside here there are four subgroups of index $3$, one of which is $3\times H'$, and the other three are isomorphic to $H$. Let $H_1$ and $H_2$ be two of them. If $g\in\mbG$ conjugates $H_1$ to $H_2$, then $g$ normalizes $Z(\mbG)\times H_1=Z(\mbG)\times H_2$. But then $g$ normalizes $H'$, and this is a clear contradiction. Thus each class of $H'$ contributes three classes of $H$, so six classes in total. This proves (i).

For the rest of the proof, note first that if $\sigma$ normalizes $H$ (or one of the subgroups of $H\times Z(\mbG)$ isomorphic to $H$) then it centralizes it, since $\sigma$ centralizes $H'$. Since $\sigma$ centralizes $H'$ it normalizes its normalizer, which is $Z(\mbG)\times H$. Note that $\Out(Z(\mbG)\times H)\cong \Sym(3)$, acting faithfully on the three subgroups isomorphic with $H$, and with elements of order $2$ inverting $Z(\mbG)$. If $\sigma$ does not centralize $Z(\mbG)$ then it inverts it, and so has order $2$. Thus $\sigma$ normalizes one of the three subgroups isomorphic to $H$, hence centralizes it. This proves (iii).

Certainly $\sigma^3$ cannot induce an outer automorphism of order $3$ on $Z(\mbG)\times H$, so it is either outer of order $2$---and so (iv) holds as (iii) held above---or is inner and centralizes $Z(\mbG)\times H'$, thus is trivial and (iv) holds again. It remains to show (ii). But in this case $\sigma$ cannot be outer of order $3$ as then $\sigma^2$ would also, so it has order $1$ or $2$, and thus (ii) holds as it held for (iii) and (iv).
\end{proof}

Assume that $\sigma$ centralizes $H'\cong \PSL_2(8)$, and write $G_\mathrm{sc}=\mbG^\sigma$. We are interested in whether $N_G(H')=H'$. The previous lemma shows that if $\mbG^\sigma$ is ${}^2\!E_6(p^3)$ or $Z(G_\mathrm{sc})=1$ then $H'.3\leq N_G(H')$, and also if $G={}^2\!E_6(p)$ then $N_G(H')$ is equal to its normalizer in $E_6(p^2)$. Thus we may assume that $\sigma$ is a power of the standard Frobenius endomorphism that centralizes $Z(\mbG)$ in what follows.

We will construct the six $E_6$-forms on $M(E_6)$ (there are six by Lemma \ref{lem:sl28prelim}(i)) and then consider the action of $\sigma$ on these. We see that $\sigma$ centralizes all six forms if and only if all six $\mbG$-conjugacy classes have a $\sigma$-fixed representative. Since $\sigma$ centralizes the centre of $\mbG$, the only other alternative is that it permutes the six forms in two $3$-cycles, and then $H$ does not embed in $\mbG^\sigma$. Thus we can determine if $H$ is contained in $\mbG^\sigma$ by considering the coefficients of the $E_6$-forms on $M(E_6)$.

With this approach we can then use the results from \cite[Section 29]{aschbacherE6Vun}. Since this work is unpublished, we reproduce (with permission) what we need from it, but we can simplify the calculations somewhat because we are considering $H$ rather than $H'$, as was done in \cite{aschbacherE6Vun}.

\medskip

We begin by defining some elements of $H$, and then construct the $27$-dimensional module $V$ that is isomorphic to $M(E_6){\downarrow_H}$. Let $p\neq 2,3,7$, and let $k$ be an algebraically closed field of characteristic $p$. Let $\xi$ be a primitive $7$th root of unity in $\F_8$ and let $\zeta$ be a primitive $7$th root of unity in $k$. Over $\F_2$, $x^7-1$ factorizes as $(x-1)(x^3+x+1)(x^3+x^2+1)$, so choose $\xi$ to satisfy $x^3+x+1$. If $i\in \F_8^\times$, write $\log(i)$ for the quantity determined by $i=\xi^{\log(i)}$. For $i\in \F_8$ and $j\in \F_8^\times$, write
\[ g_i=\begin{pmatrix} 1&i\\0&1\end{pmatrix},\quad h_j=\begin{pmatrix} j^{-1}&0\\0&j\end{pmatrix},\quad t=\begin{pmatrix}0&1\\1&0\end{pmatrix}.\]
Let $s$ be the element of $H$ that centralizes $\gen{g_1,t}\cong \SL_2(2)$ and maps $g_\xi$ to $g_{\xi^2}$. Let $B=\gen{g_1,h_\xi}$ be a Borel subgroup of $H'$.

Let $U_\zeta$ be the $9$-dimensional $kB$-module with basis $\{x_\infty\}\cup \{x_i:i \in \F_8\}$, so a projective line. The actions of $g_i$, $h_j$ and $t$ are as follows.

\begin{center}\begin{tabular}{llll}
\hline Basis element & Image under $g_i$ & Image under $h_j$ & Image under $t$
\\ \hline $x_\infty$ & $x_\infty$ & $\zeta^{\log(j)}x_\infty$ & $x_0$
\\ $x_0$ & $x_1$ & $\zeta^{-\log(j)}x_0$ & $x_\infty$
\\ $x_l$, $l\neq 0$ & $x_{l+i}$ & $\zeta^{-\log(j)}x_{lj^2}$ & $\zeta^{\log(l)}x_{l^{-1}}$
\\ \hline
\end{tabular}\end{center}
Define $Y$ and $Z$ similarly, with bases $\{y_\infty\}\cup \{y_i:i \in \F_8\}$ and $\{z_\infty\}\cup \{z_i:i \in \F_8\}$ respectively, and with actions given as for $X$ but with $\zeta$ replaced by $\zeta^4$ and $\zeta^2$ respectively. Let $V=X\oplus Y\oplus Z$, a vector space of dimension $27$.

The action of $s$ is defined by $s^3=1$, and $x_\infty s=y_\infty$, $y_\infty s=z_\infty$. From this and the fact that $sg_1=g_1s$ and $st=ts$, we see that $x_0s=y_0$, $x_1s=y_1$, and similarly $y_0s=z_0$ and $y_1s=z_1$. The rest of the basis elements are permuted in a more complicated way, but it is not necessary to know what it is in what follows. (In the supplementary materials we give this representation and prove that it is indeed a representation of $H$.)

The subgroup $B$ acts on $V$ as the sum of all six non-trivial $kB$-modules and three copies of the $7$-dimensional $kB$-module. The $1$-dimensional modules are spanned by $x_\infty$, $y_\infty$, $z_\infty$, 
\[x_\Sigma=\sum_{i\in \F_8} x_i,\quad y_\Sigma=\sum_{i\in \F_8} y_i,\quad z_\Sigma=\sum_{i\in \F_8} z_i.\]
Notice that
\[ x_\Sigma g_i=x_\Sigma,\quad x_\Sigma h_j=\zeta^{-\log(j)} x_\Sigma.\]

\medskip

Let $f$ be a symmetric trilinear form on $V$, and suppose that $f$ has symmetry group $E_6$. Recall that a point $x\in V$ is \emph{singular} if, for all $y\in V$, $f(x,x,y)=0$. (This is the same as a white or $D_5$-parabolic point of $V$.) By \cite[(3.16)(2)]{aschbacher1987} $\mbG$ is transitive on singular points, and by \cite[(3.6)(1)]{aschbacher1987}, the dimension of $x\Delta$, the radical of the bilinear form $f(x,-,-)$, is $17$-dimensional whenever $x$ is singular. Notice that the space $x\Theta$ of elements $y$ such that $f(x,x,y)=0$ is either $26$- or $27$-dimensional, since it is the kernel of the linear transformation $y\mapsto f(x,x,y)$. If $\gen x$ is a $kB$-submodule of $V$ then so are $x\Theta$ and $x\Sigma$.

In Section \ref{sec:e6;sl28}, we proved that the Borel subgroup $B$ lies in $D_5T_1$ and then inside $B_4T_1$. Thus $B$ stabilizes both $D_5$-parabolic and $B_4$-type lines on $M(E_6)$ (white and grey lines in the language of, for example, \cite{cohencooperstein1988}). Since $x_\infty$, $y_\infty$ and $z_\infty$ lie in an $s$-orbit, and $x_\Sigma$, $y_\Sigma$ and $z_\Sigma$ lie in an $s$-orbit, one must be white and one must be grey.

Since $H'$ is $3$-transitive on the projective line we determine $f$ on all triples of basis elements if we determine $f$ on all triples of basis elements with subscripts from $\{\infty,0,1\}$. Thus we define:
\[ c_x=f(x_\infty,x_0,x_1),\quad c_{xy}=f(x_\infty,y_0,y_1),\quad c_{yx}=f(y_\infty,x_0,x_1),\]
\[ c_\infty=f(x_\infty,y_\infty,z_\infty),\quad c=f(x_\infty,y_0,z_1).\]
We claim that $f(x_\infty,x_\infty,v)=0$ if $v$ is one of $x_\infty,x_0,y_\infty,y_0,z_\infty$, so the only one for which the form could possibly be non-zero is $v=z_0$. To see this, we use $h_\xi$-equivariance: for example,
\[ f(x_\infty,x_\infty,x_0)=f(x_\infty h_\xi,x_\infty h_\xi,x_0 h_\xi)=f(\zeta x_\infty,\zeta x_\infty,\zeta^{-1}x_0)=\zeta f(x_\infty,x_\infty,x_0),\]
so $f(x_\infty,x_\infty,x_0)=0$. The same holds for $v$ one of $x_\infty$, $y_0$, $y_\infty$, $z_0$ or $z_\infty$ unless $v h_\xi=\zeta^{-2} v$, which is only true for $z_0$. Since $x_\infty\Theta$ is a $kB$-submodule of dimension least $26$, it is either $V$ itself or a complement to $\gen{z_\Sigma}$.

Similarly, $f(x_\Sigma,x_\Sigma,v)=0$ for $v$ one of $x_\infty,x_\Sigma,y_\infty,y_\Sigma,z_\Sigma$, and only $f(x_\Sigma,x_\Sigma,z_\infty)$ can be non-zero. For $f(x_\infty,y_\infty,v)$, we again see that this is zero for $v=x_\infty,x_\Sigma,y_\infty,y_\Sigma,z_\Sigma$, so can be non-zero only for $z_\infty$ (where it takes value $c_\infty$).

We now can compute the value of the form when each of the basis elements lies in $X$. We do the same for $Y$ and $Z$.

\begin{lem}[{{\cite[(29.6)]{aschbacherE6Vun}}}]\label{lem:xxx} We have \[ f(x_\infty,x_i,x_j)=c_x\zeta^{4\log(i+j)},\qquad f(y_\infty,y_i,y_j)=c_x\zeta^{2\log(i+j)},\qquad f(z_\infty,z_i,z_j)=c_x\zeta^{\log(i+j)}.\]
\end{lem}
\begin{proof} Since $B$ is $2$-transitive on $\{x_i:i \in \F_8\}$, we need merely show that this formula is equivariant with respect to $g_1$ and $h_\xi$. Thus:
\[ f(x_\infty g_1,x_ig_1,x_jg_1)=f(x_\infty,x_{i+1},x_{j+1})=c_x\zeta^{4\log(i+j)},\]
and
\[ f(x_\infty h_\xi,x_ih_\xi,x_jh_\xi)=f(\zeta x_\infty,\zeta^{-1}x_{i\xi^2},\zeta^{-1}x_{j\xi^2})=\zeta^{-1}c_x\zeta^{4\log((i+j)\xi^2)}=c_x\zeta^{4\log(i+j)},\]
proving the result. To obtain the other two formulae, note that $Y$ has the same action as $X$ except for $\zeta$ replaced by $\zeta^4$ (so $\zeta^{4\log(i+j)}$ becomes $\zeta^{16\log(i+j)}=\zeta^{2\log(i+j)}$, and $Z$ has the same action as $X$ except for $\zeta$ replaced by $\zeta^2$.
\end{proof}

We have a very similar lemma when the basis elements do not all come from the same $kH'$-submodule.

\begin{lem}[{{\cite[(29.8)(2)]{aschbacherE6Vun}}}]\label{lem:xyy} We have
\[ f(x_\infty,y_i,y_j)=c_{xy},\quad f(y_\infty,x_i,x_j)=\zeta^{-\log(i+j)}c_{yx},\quad f(x_\infty,y_i,z_j)=\zeta^{-\log(i+j)}c.\]
\end{lem}
\begin{proof} The proof is the same as for Lemma \ref{lem:xxx}. We have
\[ f(x_\infty g_1,y_0 g_1,y_1 g_1)=f(x_\infty,y_1,y_0)=c_{xy},\]
and
\[ f(x_\infty h_\xi,y_ih_\xi,y_jh_\xi)=f(\zeta x_\infty,\zeta^{-4}y_{i\xi^2},\zeta^{-4}y_{j\xi^2})=\zeta^{1-8}c_{xy}=c_{xy}.\]

The other two cases are proved in the same way, and are omitted.\end{proof}

Notice that $B$ lies in $D_5T_1$, and this latter subgroup is invariant under the graph automorphism. The graph automorphism inverts the torus, so in effect sends the $1$-dimensional modules in $V{\downarrow_B}$ to their duals. This swaps $x_\infty$ and $x_\Sigma$; thus we may assume that $x_\infty$ is singular, and we will do so. We aim to produce three distinct trilinear forms on $V$, under this assumption.

As $x_\infty$ is singular, $x_\Sigma$ is non-singular, and therefore $z_\infty$ lies outside $x_\Sigma \Theta$. Hence $0\neq f(z_\infty,x_i,x_j)=f(x_\infty,y_i,y_j)$ for some $i,j\in \F_8$, and in particular this means that $c_{xy}\neq 0$.

\begin{lem} $x_\infty\Delta$ contains $x_\infty$, $x_\Sigma$ and $z_\Sigma$, and does not contain $y_\infty$, $z_\infty$ and $y_\Sigma$.
\end{lem}
\begin{proof} Since $x_\infty$ is singular, certainly $f(x_\infty,x_\infty,v)=0$ for all $v\in V$, so $x_\infty\in x_\infty \Delta$. If $y_\infty\in x_\infty\Delta$ then by applying $s^{-1}$ we see that $z_\infty\in x_\infty\Delta$, as $f(x_\infty,y_\infty,v)=0$ for all $v\in V$ if and only if $f(z_\infty,x_\infty,v)=0$ for all $v\in V$. Since $x_\infty\Delta$ is a $kB$-module of dimension $17$, it contains two copies of the $7$-dimensional $kB$-module and three $1$-dimensional modules.

We now show that $x_\Sigma$ lies in $x_\infty\Delta$, thus proving that $y_\infty$ and $z_\infty$ do not lie in it. For this we show that $f(x_\infty,x_\Sigma,v)=0$ for $v$ each basis element of $V$. Note that
\[ f(x_\infty h_\xi,x_\Sigma h_\xi,vh_\xi)=f(x_\infty,x_\Sigma,vh_\xi),\]
and if $v\in\{x_\infty,x_0,y_\infty,y_0,z_\infty,z_0\}$ then $v h_\xi$ is a non-unity multiple of $v$. Hence $f(x_\infty,x_\Sigma,v)=0$ for those elements. We then use the fact that $B$ is transitive on $x_i$ for $i\in \F_8$, and stabilizes $x_\infty$ and $x_\Sigma$, to see that $f(x_\infty,x_\Sigma,v)=0$ for all $v\in V$.

Thus it remains to check whether $y_\Sigma\in x_\infty\Delta$ or $z_\Sigma\in x_\infty\Delta$. Notice that
\[ f(x_\infty h_\xi,z_\Sigma h_\xi,vh_\xi)=\zeta^{-1}f(x_\infty,x_\Sigma,vh_\xi),\]
and if $v\in\{x_0,y_\infty,y_0,z_\infty,z_0\}$ then $v h_\xi\neq \zeta v$, but is some other scalar multiple of $v$. The last case, $v=x_\infty$, does work, but $x_\infty$ is singular so the form is zero for this option for $v$ as well. Thus, as the previous case, $z_\Sigma\in x_\infty\Delta$.

Alternatively, one sees that $f(x_\infty,y_\Sigma,y_\Sigma)=56c_{xy}\neq 0$ by Lemma \ref{lem:xyy}.
\end{proof}

From this we see that $f(x_\infty,y_\infty,v)\neq 0$ for some basis element $v$. As with the previous proof, we see that the form vanishes on this triple for $v$ one of $x_\infty,x_0,y_\infty,y_0,z_0$, and so the only option is that $c_\infty=f(x_\infty,y_\infty,z_\infty)\neq 0$.

%

Thus we have so far proved that $c_{xy}$ and $c_\infty$ are non-zero. In fact, $c_x$, $c_{yx}$ and $c$ are also non-zero, but we need to determine the two $7$-dimensional submodules in $x_\infty\Delta$ to prove that.

\medskip

Inside $X$, the $7$-dimensional submodule is simply all elements
\[ \sum_{i\in \F_8} a_ix_i\]
such that $\sum a_i=0$. The same holds for $Y$ and $Z$. However the isomorphism from this summand to $X$ to this summand of $Y$ is more delicate to write down. In \cite[pp.141--142]{aschbacherE6Vun} a `fixed point' is found for this map, and we show this directly in the supplementary materials. Write
\[ x_e=(x_0+x_\xi+x_{\xi^2}+x_{\xi^4})-(x_1-x_{\xi^3}-x_{\xi^5}-x_{\xi^6}).\]
(Note that $\xi$, $\xi^2$ and $\xi^4$ are those powers of $\xi$ that satisfy the polynomial $x^3+x+1$, and the other roots satisfy $x^3+x^2+1$.) We define $y_e$ and $z_e$ similarly. In any isomorphism from the $7$-dimensional summand of $X$ to that of $Y$, $x_e$ is mapped to (a scalar multiple of) $y_e$, and similarly for $z_e$. Thus a generic $7$-dimensional $kB$-submodule of $V$ is generated by $u=\alpha x_e+\beta y_e+\gamma z_e$, where $\alpha,\beta,\gamma\in k$.

Suppose that $u$ lies in $x_\infty \Delta$. By evaluating $f(x_\infty,u,v)=0$ for various basis elements $v$ we obtain equations in $\alpha$, $\beta$ and $\gamma$ that must be satisfied. We use these to produce relations between the coefficients $c_x$, $c_{xy}$, $c_{yx}$ and $c$. (Clearly $c_\infty$ will not appear.)

To evaluate $f(x_\infty,u,x_0)$, we need an easy lemma. First, write
\[ \omega=\zeta+\zeta^2+\zeta^4-(\zeta^{-1}+\zeta^{-2}-\zeta^{-4});\]
note that $\omega^2=-7$.

\begin{lem}\label{lem:firstsimul} We have
\[f(x_\infty,x_0,x_j)=\zeta^{4\log(j)}c_x,\quad  f(x_\infty,x_0,y_j)=\zeta^{2\log(j)}c_{yx},\quad f(x_\infty,x_0,z_j)=\zeta^{\log(j)}c_{xy}.\]
Consequently, 
\[ f(x_\infty,x_e,x_0)=c_x(\omega-1),\quad f(x_\infty,y_e,x_0)=c_{yx}(\omega-1),\quad f(x_\infty,z_e,x_0)=c_{xy}(\omega-1).\]
\end{lem}
\begin{proof} The first case follows from Lemma \ref{lem:xxx}. For the second, apply $g_j$, then $t$, then evaluate using Lemma \ref{lem:xyy}. For the third, apply $g_j$, then $t$, then $s$, then evaluate using Lemma \ref{lem:xyy}.

To see the consequence, note that $f(x_\infty,x_0,x_0)=f(x_\infty,y_0,x_0)=f(x_\infty,z_0,x_0)=0$. Then the conclusion follows easily, once one notices that $\omega$ is invariant under replacing $\zeta$ by either $\zeta^2$ or $\zeta^4$.
\end{proof}

From this we see that if $u\in x_\infty\Delta$ then 
\begin{equation}
\alpha c_x+\beta c_{yx}+\gamma c_{xy}=0.\label{eq:firstsimul}
\end{equation}

We now do the same thing, but with $y_0$ instead of $x_0$. 

\begin{lem} We have
\[f(x_\infty,y_0,x_j)=\zeta^{2\log(j)}c_{yx},\quad  f(x_\infty,y_0,y_j)=c_{xy},\quad f(x_\infty,y_0,z_j)=\zeta^{-\log(j)}c.\]
Consequently, 
\[ f(x_\infty,x_e,y_0)=c_{yx}(\omega-1),\quad f(x_\infty,y_e,y_0)=-c_{xy},\quad f(x_\infty,z_e,y_0)=-c(\omega+1).\]
\end{lem}
\begin{proof} The first case follows by applying $t$ then evaluating using Lemma \ref{lem:xyy}. The second and third cases follow directly from Lemma \ref{lem:xyy}. The consequences are proved as with Lemma \ref{lem:firstsimul}.
\end{proof}

From this we see that if $u\in x_\infty\Delta$ then 
\begin{equation}
\alpha c_{yx}(\omega-1)-\beta c_{xy}-\gamma c(\omega+1)=0.\label{eq:secondsimul}
\end{equation}

Finally we use $z_0$.

\begin{lem} We have
\[f(x_\infty,z_0,x_j)=\zeta^{\log(j)}c_{xy},\quad  f(x_\infty,z_0,y_j)=\zeta^{-\log(j)}c,\quad f(x_\infty,z_0,z_j)=\zeta^{-\log(j)}c_{yx}.\]
Consequently, 
\[ f(x_\infty,z_0,x_e)=c_{xy}(\omega-1),\quad f(x_\infty,z_0,y_e)=-c(\omega+1),\quad f(x_\infty,z_0,z_e)=-c_{yx}(\omega+1).\]
\end{lem}
\begin{proof} The first case follows by applying $t$ then $s$, then evaluating using Lemma \ref{lem:xyy}. The second and third cases follow from applying $g_j$ and $s$ respectively, then applying Lemma \ref{lem:xyy}. The consequences are proved as with Lemma \ref{lem:firstsimul}.
\end{proof}

From this we see that if $u\in x_\infty\Delta$ then 
\begin{equation}
\alpha c_{xy}(\omega-1)-\beta c(\omega+1)-\gamma c_{yx}(\omega+1)=0.\label{eq:thirdsimul}
\end{equation}

We now consider simultaneous solutions $(\alpha,\beta,\gamma)$ to (\ref{eq:firstsimul}), (\ref{eq:secondsimul}) and (\ref{eq:thirdsimul}). Note that, as $c_{xy}\neq 0$, $(1,0,0)$ is not a solution to (\ref{eq:thirdsimul}), $(0,1,0)$ is not a solution to (\ref{eq:secondsimul}) and $(0,0,1)$ is not a solution to (\ref{eq:firstsimul}).

Since there is a $2$-space of solutions, there exists a unique solution $(0,1,\gamma)$ for some $\gamma\in k$. The three equations become
\[ c_{yx}+\gamma c_{xy}=c_{xy}+\gamma(\omega+1)c=c+\gamma c_{yx}=0.\]
Since $c_{xy}\neq 0$, we see that both $c_{yx}$ and $c$ are non-zero. Solving for $\gamma$ yields
\begin{equation} \gamma=-c_{yx}/c_{xy}=c_{xy}/(\omega+1)c=-c/c_{yx}.\label{eq:firstsolution}\end{equation}

We do the same with a solution of the form $(1,0,\gamma)$ to yield
\[ c_x+\gamma c_{xy}=c_{yx}(\omega-1)-\gamma c(\omega+1)=c_{xy}(\omega-1)-\gamma c_{yx}(\omega+1).\]
This time we see that $c_x$ is also non-zero. Solving for $\gamma$ yields
\begin{equation}
\gamma=-c_x/c_{xy}=c_{yx}(\omega-1)/c(\omega+1)=c_{xy}(\omega-1)/c_{yx}(\omega+1).\label{eq:thirdsolution}
\end{equation}

We may scale the trilinear form so that $c_{xy}=1$. Then from (\ref{eq:firstsolution}) we obtain $c=c_{yx}^2$ and $c_{yx}^3=-1/(\omega+1)$, and from (\ref{eq:thirdsolution}) we obtain $c_x=(1-\omega)/c_{yx}(\omega+1)=c_{yx}^2(\omega-1)$. The only parameter we have not fixed yet is $c_\infty$; it is proved in \cite[p.144]{aschbacherE6Vun} that $c_\infty$ is uniquely determined using special planes and hyperbolic subspaces. We can be much more naive, as we have Proposition \ref{prop:sl28ine6} already. (In \cite{aschbacherE6Vun} this statement was used to prove Proposition \ref{prop:sl28ine6}.) Thus we know that there is a unique (up to $C_{\GL(V)}(H')$-conjugacy) $E_6$-form for which $x_\infty$ is singular. Thus some element of $C_{\GL(V)}(H')$ conjugates any such $E_6$-form to any other, in particular, between two with the same values for $c_x$, $c_{xy}$, $c_{yx}$ and $c$, but with different values of $c_\infty$.

But $C_{\GL(V)}(H')$ is simply a scalar matrix acting on each of $X$, $Y$ and $Z$. Since $c_x$ is fixed, the scalar on each of $X$, $Y$ and $Z$ is a cube root of unity, and since $c_{xy}$ is fixed, it is the same cube root of unity for each factor. But then this is simply the centre of $E_6$, so this does not affect $c_\infty$, and the result is proved.


\medskip

We thus have six trilinear forms: the three given by the three cube roots of $-1/(\omega+1)$, and the three given by their images under the graph automorphism (where $x_\infty$ is non-singular). Let $\sigma$ denote a $q$-power map $x\mapsto x^q$ on $k$, for $q$ a power of $p$. Suppose first that $\sigma$ fixes $\zeta$ ($q\equiv 1\bmod 7$), so that $\sigma$ centralizes the $H$-action on $V$, and hence permutes the six trilinear forms on $V$. If $\sigma$ fixes the cube roots of $\omega+1$ then $\sigma$ fixes all six trilinear forms, and in particular $H$ embeds in the simple group $E_6(q)$. On the other hand, if $\sigma$ does not fix the cube roots of $\omega+1$ then $\sigma$ does not fix any trilinear form, and so $H$ cannot embed in $E_6(q)$.

If $\sigma$ does not fix $\zeta$ ($q\not\equiv 1\bmod 7$) then $\sigma$ maps the $H$-action on $V$ to a slightly different action. In order to obtain a permutation of the trilinear forms, we must replace the $q$-power map with the product of that and an element of $\GL(V)$ that maps this twisted $H$-action back to the original one.

However, the effect of the $q$-power map is easy to see: it cycles bodily the subspaces $X$, $Y$ and $Z$ (which way depends on whether $\zeta$ maps to $\zeta^2$ or $\zeta^4$), and so multiplying by a permutation matrix
\[ \begin{pmatrix}0&I&0\\0&0&I\\I&0&0\end{pmatrix},\]
where $I$ is a $9\times 9$ identity matrix, is enough to centralize $H$ again. Notice that the product of a $q$-power map and a permutation matrix does not affect the value of $f(x_\infty,y_0,y_1)$ (as the permutation matrix cycles the factors), so the effect on the six forms is the same as the $q$-power map itself. Thus the same conclusion holds, and we have the following theorem.

\begin{thm}\label{thm:sl28} Let $H\cong \SL_2(8).3$, and let $q$ be a power of $p\neq 2,3,7$. Let $\omega$ denote a square root of $-7$ in $\bar\F_p$.
\begin{enumerate}
\item The group $H'$ embeds in the simple group $E_6(q)$, acting irreducibly on $M(E_6)$, if and only if $q\equiv 1,2,4\bmod 7$, and in ${}^2\!E_6(q)$, acting irreducibly on $M(E_6)$, if and only if $q\equiv -1,-2,-4\bmod 7$.
\item The group $H$ embeds in the simple group $E_6(q)$ if and only if $q\equiv 1,2,4\bmod 7$, and one of the following holds:
\begin{enumerate}
\item $q\equiv 2\bmod 3$;
\item $q=p^{3n}$ for some integer $n$;
\item $\omega+1$ has a cube root in $\F_q$.
\end{enumerate}
\item The group $H$ embeds in the simple group ${}^2\!E_6(q)$ if and only if it embeds in $E_6(q^2)$.
\item If $H$ embeds in ${}^\varepsilon\!E_6(p)$ then it is always maximal in ${}^\varepsilon\!E_6(p)$. If $H$ embeds in ${}^\varepsilon\!E_6(p).3$ but not ${}^\varepsilon\!E_6(p)$ then it is always maximal in ${}^\varepsilon\!E_6(p).3$. If $H$ does not embed in ${}^\varepsilon\!E_6(p)$ then $H'$ is maximal in ${}^\varepsilon\!E_6(p)$ if and only if $x^3-3x+1$ does not split over $\F_p$, i.e., $p\not\equiv \pm 1\bmod 9$.
\end{enumerate}
\end{thm}

\section{The remaining maximal subgroups of \texorpdfstring{$F_4(q)$}{F4(q)} and \texorpdfstring{${}^\varepsilon\!E_6(q)$}{±E6(q)}}
\label{sec:remainingmaximals}

Let $\mbG$ be of type $F_4$ or $E_6$. Let $H$ be a maximal subgroup of $G_{\mathrm{sc}}=\mbG^\sigma$ that is not a member of $\mathcal{S}$. As mentioned in Section \ref{sec:prelim}, $H$ is one of: the fixed points $\mb X^\sigma$ for $\mb X$ a maximal positive-dimensional subgroup of $\mbG$; a subgroup of the same type as $G$; an exotic $r$-local subgroup for some $r\neq p$. The exotic $r$-local subgroups are given in \cite{clss1992}. The maximal-rank subgroups that are $\mb X^\sigma$ are given in \cite[Tables 5.1 and 5.2]{liebecksaxlseitz1992}, and the other maximal positive-dimensional subgroups appear in \cite{liebeckseitz2004}.

For $F_4(q)$, since there are no diagonal outer automorphisms the tables are fairly easy to write down, and we have Table \ref{tab:f4othermaximalsqodd} for $q$ odd and Table \ref{tab:f4othermaximalsqeven} for $q$ even. The tables are split into two because of the presence of a large number of novelty maximal subgroups whenever $\bar G$ induces a graph automorphism on $G$.

For a simple group of type $E_6$, the tables become quite complicated as there are field, diagonal and graph automorphisms. All of the subgroups are stabilized by field automorphisms, but stability under diagonal and graph automorphisms are more complicated.

Let $d=\gcd(2,q-1)$, $e=\gcd(3,q-1)$ and $e'=\gcd(3,q+1)$, $f=\gcd(4,q-1)$ and $f'=\gcd(4,q+1)$, and $g=\gcd(5,q-1)$. We write $\delta$, of order $e$ for $E_6(q)$ and $e'$ for ${}^2\!E_6(q)$, for a generator of the group of diagonal automorphisms, $\gamma$ for a graph automorphism, and $\phi$ for a generator of the group of field automorphisms. If $G={}^2\!E_6(q)$ then $\gamma$ lies in $\gen\phi$, so we simply use $\phi$ in this case.

The tables include the subgroup $\bar H$ of $G$ for which $N_{\bar G}(\bar H)$ is maximal in $\bar G$, the conditions on $p$ and $q$, and the number of classes of $\bar H$ in $G$. Note that we have written $\bar H$ so that $N_G(\bar H)=\bar H$, but of course $N_{\bar G}(\bar H)$ will be larger. We include the stabilizer in $\Out(G)$ of $\bar H$ so that the reader may compute the maximal subgroups of the almost simple groups $\bar G$.

We now show why the particular groups in these tables appear for $E_6$ and ${}^2\!E_6$. As we said earlier, we obtain from \cite{liebecksaxlseitz1992} the maximal-rank subgroups, the parabolics are simple to understand (if not describe), and the local maximal subgroups appear in \cite{clss1992}; so we must deal with reductive subgroups that are not maximal rank. These are: $F_4$, $C_4$, $A_2G_2$, $G_2$ (two classes) and $A_2$ (two classes).

For $\mb X=F_4$ there is little to say. It is centralized by the graph automorphism of $\mbG$ so appears in both tables. The fixed points must be $F_4(q)$, and the diagonal automorphism cannot normalize them, so there are $e$ or $e'$ classes.

For $\mb X=C_4$ the composition factors on the minimal module show that $\mb X$ has adjoint type in $\mbG$ (and it is not the centralizer of an involution). As it is centralized by a graph automorphism, it appears in both tables. Since the diagonal automorphism of $G$ has order $3$ it cannot induce the diagonal automorphism of $\PSp_8(q)$, so $\delta$ cannot normalize $\mb X^\sigma$. Thus $\mb X^\sigma$ is $\PSp_8(q).2$ and there are $e$ or $e'$ classes.

For $\mb X=A_2G_2$, the group is self-normalizing in $\mbG$, contrary to the statement in \cite[p.3]{liebeckseitz2004}; this references \cite[(3.15)]{seitz1991}, which does not prove this fact, and indeed the composition factors of $M(E_6){\downarrow_{\mb X}}$ do not allow this. In addition, $\mb X$ is normalized by the graph automorphism, so if $|N_\mbG(\mb X)/\mb X|=2$ then $\mb X$ would be centralized by the graph, so lie in $F_4$ or $C_4$. In the simply connected group the composition factors on the minimal module show the $A_2$ is simply connected. Thus the form of the group in $E_6(q)$ is $\PSL_3(q)\times G_2(q)$ or $\PSU_3(q)\times G_2(q)$, depending on the action of $\sigma$ on the $A_2$ factor. In addition, the diagonal automorphism must induce the diagonal automorphism on the $A_2$-factor, so there is a single class.

Since the graph automorphism normalizes $\mb X$, exactly one of $E_6(q)$ and ${}^2\!E_6(q)$ contains $\PSL_3(q)\times G_2(q)$ and the other contains $\PSU_3(q)\times G_2(q)$. From \cite[(5.7.6)]{aschbacher1990a}, we see that $\PSL_3(q)\times G_2(q)$ embeds in $E_6(q)$, completing the proof.

For $\mb X=G_2$, the graph automorphism swaps the two classes of $\mb X$. Of course $N_\mbG(\mb X)=\mb X\cdot Z(\mbG)$. We must decide if the standard Frobenius morphism of $\mbG$ stabilizes the classes or swaps them. In \cite[Theorem G.2]{testerman1989} it is shown that if $q\equiv 1,2,4\bmod 7$ then there are two (non-conjugate) subgroups $G_2(q)$ in $E_6(q)$, swapped by the graph automorphism. If $q\equiv 3,5,6\bmod 7$ then the standard Frobenius map swaps the two square roots of $-7$, so maps one class from that theorem to the other. In particular, the product with the graph automorphism does stabilize the class, so we find two classes of $G_2(q)$ in ${}^2\!E_6(q)$ in this case. Since the diagonal automorphism cannot normalize $G_2(q)$, we see that there are $2e$ classes in $E_6(q)$ and $2e'$ classes in ${}^2\!E_6(q)$. (This tallies with \cite[p.5]{aschbacherE6Vun}.)

Finally, for $\mb X=A_2$, the graph automorphism swaps the two classes of $\mb X$. Note that $\mbG$ induces a graph automorphism on $\mb X$ by \cite[Claim on p.314]{testerman1989}. The composition factors of $M(E_6){\downarrow_{\mb X}}$ show that $\mb X$ has adjoint type (and it cannot centralize an element of order $3$, of course). Thus the subgroup of $E_6(q)$ is either $\PGL_3(q).2$ or $\PGU_3(q).2$, depending on the action of $\sigma$. As with $G_2$, in \cite[Theorem A.2]{testerman1989} we find that if the standard Frobenius map fixes $\sqrt{-1}$ then it stabilizes the two $\mbG$-classes, and if it does not fix $\sqrt{-1}$ then it does not stabilize the two classes. Thus if $q\equiv 1\bmod 4$ then we should see classes in $E_6(q)$ and if $q\equiv 3\bmod 4$ then we should see classes in ${}^2\!E_6(q)$.

Unlike in the $G_2$ case though, $N_{\mbG}(\mb X)/Z(\mbG)\cdot\mb X$ has order $2$, so $N_{\mbG}(\mb X)$ is disconnected in the adjoint version of $E_6$. Lemma \ref{lem:fixedpoints} shows that there are exactly two classes in the adjoint version of $\mbG^\sigma$, and are the fixed points under $\sigma$ and $w\sigma$, where $w$ is an element that induces the graph automorphism on $A_2$. But then this means that the two classes are not the same type, so one must be $\PGL_3(q).2$ and the other $\PGU_3(q).2$. This is true for both $\mbG$-classes, so we obtain four classes in the adjoint group, swapped in pairs by the graph. Thus we obtain $2e$ classes of each in $E_6(q)$ and $2e'$ classes of each in ${}^2\!E_6(q)$. (This tallies with \cite[pp.5--6]{aschbacherE6Vun}.)

\begin{table}
\begin{center}
\begin{tabular}{llll}
\hline Group & Conditions & No. classes & Stabilizer
\\\hline $[q^{15}].\Sp_6(q).(q-1)$ & - & $1$ & $\gen\phi$
\\ $[q^{20}].(\SL_2(q)\times \SL_3(q)).(q-1)$ & - & $1$ & $\gen\phi$
\\ $[q^{20}].(\SL_3(q)\times \SL_2(q)).(q-1)$ & - & $1$ & $\gen\phi$
\\ $[q^{15}].2\cdot\Omega_7(q).(q-1)$ & - & $1$ & $\gen\phi$
\\ $2\cdot \Omega_9(q)$ & - & $1$ & $\gen\phi$
\\ $2^2\cdot \POmega_8^+(q).\Sym(3)$ & -  & $1$ & $\gen\phi$
\\ ${}^3\!D_4(q).3$ & - & $1$ & $\gen\phi$
\\ $(\Sp_6(q)\circ \SL_2(q)).2$ & - & $1$ & $\gen\phi$
\\ $(\SL_3(q)\circ \SL_3(q)).\gcd(3,q-1).2$ & - & $1$ & $\gen\phi$
\\ $(\SU_3(q)\circ \SU_3(q)).\gcd(3,q+1).2$ & - & $1$ & $\gen\phi$
\\ $\PGL_2(q)\times G_2(q)$ & $q\neq 3$ & $1$ & $\gen\phi$
\\ $F_4(q_0)$ & $q=q_0^r$, $r$ prime & $1$ & $\gen\phi$
\\ $\PGL_2(q)$ & $p\geq 13$ & $1$ & $\gen\phi$
\\ $G_2(q)$ & $p=7$ & $1$ & $\gen\phi$
\\ $3^3\rtimes \SL_3(3)$ & $q=p\geq 5$ & $1$ & $1$
\\ \hline
\end{tabular}
\end{center}
\caption{Subgroups $H$ such that $N_{\bar G}(H)$ is a maximal subgroup of an almost simple group $\bar G$ with socle $F_4(q)$ not belonging to $\mathcal{S}$, $q$ odd.}
\label{tab:f4othermaximalsqodd}
\end{table}

\begin{table}
\begin{center}
\begin{tabular}{llll}
\hline Group & Conditions & No. classes & Stabilizer
\\\hline $[q^{15}].\Sp_6(q)\times (q-1)$ & - & $2$ & $\gen\phi$
\\ $[q^{20}].(\SL_2(q)\times \SL_3(q)).(q-1)$ & - & $2$ & $\gen\phi$
\\ $\Sp_8(q)$ & - & $2$ & $\gen\phi$
\\ $\POmega_8^+(q).\Sym(3)$ & - & $2$ & $\gen\phi$
\\ ${}^3\!D_4(q).3$ & - & $2$ & $\gen\phi$
\\ $(\SL_3(q)\circ \SL_3(q)).\gcd(3,q-1).2$ & - & $1$ & $\gen{\gamma,\phi}$
\\ $(\SU_3(q)\circ \SU_3(q)).\gcd(3,q+1).2$ & - & $1$ & $\gen{\gamma,\phi}$
\\ $F_4(q_0)$ & $q=q_0^r$, $r$ prime & $1$ & $\gen{\gamma,\phi}$
\\ ${}^2\!F_4(q)$ & $q$ an odd power of $2$ & $1$ & $\gen{\gamma,\phi}$
\\ $[q^{20}].\Sp_4(q).(q-1)^2$ & {Nov.} & $1$ & $\gen{\gamma,\phi}$
\\ $[q^{22}].(\SL_2(q)\times \SL_2(q)).(q-1)^2$ & {Nov.} & $1$ & $\gen{\gamma,\phi}$
\\ $(\Sp_4(q)\times \Sp_4(q)).2$ & {Nov.} & $1$ & $\gen{\gamma,\phi}$
\\ $\Sp_4(q^2).2$ & {Nov.} & $1$ & $\gen{\gamma,\phi}$
\\ $(q-1)^4.W(F_4)$ & {Nov.}, $q>4$ & $1$ & $\gen{\gamma,\phi}$
\\ $(q+1)^4.W(F_4)$ & {Nov.}, $q>2$ & $1$ & $\gen{\gamma,\phi}$
\\ $(q^2+q+1)^2.(3\times \SL_2(3))$ & {Nov.} & $1$ & $\gen{\gamma,\phi}$
\\ $(q^2+1)^2.(\SL_2(3)\rtimes 4)$ & {Nov.}, $q>2$ & $1$ & $\gen{\gamma,\phi}$
\\ $(q^2-q+1)^2.(3\times \SL_2(3))$ & {Nov.}, $q>2$ & $1$ & $\gen{\gamma,\phi}$
\\ $(q^4-q^2+1).12$ & {Nov.}, $q>2$ & $1$ & $\gen{\gamma,\phi}$
\\ \hline
\end{tabular}
\end{center}
\caption{Subgroups $H$ such that $N_{\bar G}(H)$ is a maximal subgroup of an almost simple group $\bar G$ with socle $F_4(q)$ not belonging to $\mathcal{S}$, $q$ even. Those labelled `Nov.' are novelty maximals that only occur when $\bar G$ induces a graph automorphism on $G$.}
\label{tab:f4othermaximalsqeven}
\end{table}

\begin{table}
\begin{center}
\begin{tabular}{llll}
\hline Group & Conditions & No. classes & Stabilizer
\\\hline $[q^{16}].\Spin_{10}^+(q).(q-1)/e$& - & $2$ & $\gen{\delta,\phi}$
\\ $[q^{21}].d.\PSL_6(q).(q-1)$ & - & $1$ & $\gen{\delta,\gamma,\phi}$
\\ $[q^{25}].(\SL_2(q)\times \SL_5(q)).(q-1)/e$ & - & $2$ & $\gen{\delta,\phi}$
\\ $[q^{29}].(\SL_3(q)\circ \SL_3(q)\times \SL_2(q)).(q-1)$ & - & $1$ & $\gen{\delta,\gamma,\phi}$
\\ $d.(\PSL_2(q)\times \PSL_6(q)).d$ & - & $1$ & $\gen{\delta,\gamma,\phi}$
\\ $e.(\PSL_3(q)^{\times 3}).e.\Sym(3)$ & - & $1$ & $\gen{\delta,\gamma,\phi}$
\\ $e'.(\PSL_3(q^2)\times \PSU_3(q)).e'.2$ & - & $1$ & $\gen{\delta,\gamma,\phi}$
\\ $\PSL_3(q^3).3$ & - & $1$ & $\gen{\delta,\gamma,\phi}$
\\ $d^2.(\POmega_8^+(q)\times ((q-1)/d)^2/e).d^2.\Sym(3)$ & $q>2$ & $1$ & $\gen{\delta,\gamma,\phi}$
\\ $({}^3\!D_4(q)\times (q^2+q+1)/e).3$ & - & $1$ & $\gen{\delta,\gamma,\phi}$
\\ $((q-1)^6/e).W(E_6)$ & $q>4$ & $1$ & $\gen{\delta,\gamma,\phi}$
\\ $(q^2+q+1)^3/e.(3^{1+2}.\SL_2(3))$ & - & $1$ & $\gen{\delta,\gamma,\phi}$
\\ $F_4(q)$ & - & $e$ & $\gen{\gamma,\phi}$
\\ $\PSp_8(q).2$ & $p\neq 2$ & $e$ & $\gen{\gamma,\phi}$
\\ $\PSL_3(q)\times G_2(q)$ & - & $1$ & $\gen{\delta,\gamma,\phi}$ 
\\ $G_2(q)$ & $q\equiv 1,2,4\bmod 7$ & $2e$ &  $\gen{\phi\gamma}$ or $\gen\phi$
\\ $\PGL_3(q).2$ & $p\geq 5$, $q\equiv 1\bmod 4$ & $2e$ & $\gen{\phi\gamma}$ or $\gen\phi$
\\ $\PGU_3(q).2$ & $p\geq 5$, $q\equiv 1\bmod 4$ & $2e$ & $\gen{\phi\gamma}$ or $\gen\phi$
\\ 
$E_6(q_0).\gcd(e,r)$ & $q=q_0^r$, $r$ prime & $\gcd(e,r)$ & $\gen{\delta^{\gcd(e,r)},\gamma,\phi}$
\\ ${}^2\!E_6(q_0)$ & $q=q_0^2$ & $\gcd(q_0-1,3)$ & $\gen{\delta^{\gcd(q_0-1,3)},\gamma,\phi}$
\\ $3^{3+3}\rtimes \SL_3(3)$ & $q=p\geq 5$, $3\mid p-1$ & $3$ & $\gen\gamma$
\\ $[q^{24}]. \Spin_8^+(q).(q-1)^2/e$ & {Nov.} & $1$ & $\gen{\delta,\gamma,\phi}$
\\ $[q^{31}].(\SL_2(q)\times \SL_2(q) \times \PSL_3(q)).(q-1)^2/e$ & {Nov.} & $1$ & $\gen{\delta,\gamma,\phi}$
\\ $f.(\POmega^+_{10}(q)\times (q-1)/ef).f$ & {Nov.} & $1$ & $\gen{\delta,\gamma,\phi}$
\\ \hline
\end{tabular}
\end{center}
\caption{Subgroups $H$ such that $N_{\bar G}(H)$ is a maximal subgroup of an almost simple group $\bar G$ with socle $E_6(q)$ not belonging to $\mathcal{S}$. Those labelled `Nov.' are novelty maximals that only occur when $\bar G$ induces a graph automorphism on $G$.}
\label{tab:e6othermaximals}
\end{table}

\begin{table}
\begin{center}
\begin{tabular}{llll}
\hline Group & Conditions & No. classes & Stabilizer
\\\hline $[q^{21}].d.\PSU_6(q).(q-1)$ & - & $1$ & $\gen{\delta,\phi}$
\\ $[q^{24}].\Spin_8^-(q).(q^2-1)/e'$ & - & $1$ & $\gen{\delta,\phi}$
\\ $[q^{29}].(\PSL_3(q^2)\times \SL_2(q)).(q-1)$ & - & $1$ & $\gen{\delta,\phi}$
\\ $[q^{31}].(\SL_3(q)\times \SL_2(q^2)).(q^2-1)/e'$ & - & $1$ & $\gen{\delta,\phi}$
\\ $f'.(P\Omega_{10}^-(q)\times (q+1)/e'f').f'$ & - & $1$ & $\gen{\delta,\phi}$
\\  $d.(\PSL_2(q)\times \PSU_6(q)).d$ & - & $1$ & $\gen{\delta,\phi}$
\\ $e'.(\PSU_3(q)^{\times 3}).e'.\Sym(3)$ & - & $1$ & $\gen{\delta,\phi}$
\\ $e.(\PSL_3(q^2)\times \PSL_3(q)).e.2$ & - & $1$ & $\gen{\delta,\phi}$
\\ $\PSU_3(q^3).3$ & - & $1$ &  $\gen{\delta,\phi}$
\\ $d^2.(\POmega_8^+(q)\times ((q+1)/d)^2/e').d^2.\Sym(3)$ & - & $1$ & $\gen{\delta,\phi}$
\\ $({}^3\!D_4(q)\times (q^2-q+1)/e').3$ & $q>2$ & $1$ & $\gen{\delta,\phi}$
\\ $(q+1)^6/e'.W(E_6)$ & $q>2$ & $1$ & $\gen{\delta,\phi}$
\\ $(q^2-q+1)^3/e'.(3^{1+2}.\SL_2(3))$ & $q>2$ & $1$ & $\gen{\delta,\phi}$
\\ $F_4(q)$ & - & $e'$ & $\gen\phi$
\\ $\PSp_8(q).2$ & $p\neq 2$ & $e'$ & $\gen\phi$ 
\\ $\PSU_3(q)\times G_2(q)$ & - & $1$ & $\gen{\delta,\phi}$
\\ $G_2(q)$ & $q\equiv 3,5,6\bmod 7$ & $2e'$ & $\gen{\phi^2}$
\\ $\PGL_3(q).2$ & $p\geq 5$, $q\equiv 3\bmod 4$ & $2e'$ & $\gen{\phi^2}$
\\ $\PGU_3(q).2$ & $p\geq 5$, $q\equiv 3\bmod 4$ & $2e'$ & $\gen{\phi^2}$
\\ ${}^2\!E_6(q_0).\gcd(e',r)$ & $q=q_0^r$, $r$ odd prime & $\gcd(e',r)$ & $\gen{\delta^{\gcd(e',r)},\phi}$
\\ $3^{3+3}\rtimes \SL_3(3)$ & $q=p\geq 5$, $3\mid p+1$ & $3$ & $\gen\phi$
\\ $({}^3\!D_4(q)\times (q^2-q+1)/e').3$ & {Nov.}, $q=2$ & $1$ & $\gen{\delta,\phi}$
\\ $(q+1)^6/e'.W(E_6)$ & {Nov.}, $q=2$ & $1$ & $\gen{\delta,\phi}$
\\ \hline 
\end{tabular}
\end{center}
\caption{Subgroups $H$ such that $N_{\bar G}(H)$ is a maximal subgroup of an almost simple group $\bar G$ with socle ${}^2\!E_6(q)$ not belonging to $\mathcal{S}$. Those below the line are novelty maximals that only occur when $\bar G$ induces a diagonal automorphism on $G$.}
\label{tab:2e6othermaximals}
\end{table}

\section{Checking overgroups}
\label{sec:proofmaximal}

What we have not yet done is to check that there are no overgroups of the claimed maximal subgroups, apart from those mentioned in the text when they are not maximal, or form novelty maximal subgroups. The quasisimple subgroups $H$ we need to consider are those in Tables \ref{tab:f4curlyS}, \ref{tab:e6curlyS} and \ref{tab:2e6curlyS}, and the possible overgroups $X$ are from those tables, and Tables \ref{tab:f4othermaximalsqodd} to \ref{tab:2e6othermaximals}.

Suppose that both $H$ and $X$ are members of $\mathcal{S}$. Notice that no group $\PSL_2(r)$ contains a subgroup $\PSL_2(s)$ unless $r$ is a power of $s$ or $\PSL_2(s)\cong \Alt(5)$, which does not appear in our tables. Thus $X$ must not be of type $\PSL_2(r)$. 

Let $G=F_4(q)$ first. Thus $X$ is one of $\PSL_4(3)$ and ${}^3\!D_4(2)$. The group $\PSL_4(3)$ cannot contain any of the other groups in Table \ref{tab:f4curlyS}. The group ${}^3\!D_4(2)$ only contains $\PSL_2(8)$, and this is unique up to conjugacy. The restriction of $26$ to it has factors $8^2,7,1^3$ in characteristic $0$, and so this cannot be the one in Table \ref{tab:f4curlyS}.

Now let $G=E_6(q)$, so that $X$ is one of $M_{12}$, $J_3$ and ${}^2\!F_4(2)'$. The group $M_{12}$ contains $\PSL_2(11)$, but this does not appear in Table \ref{tab:e6curlyS} in characteristic $5$. The group $J_3$ contains $\PSL_2(19)$, and this is why $\PSL_2(19)$ does not appear in Table \ref{tab:e6curlyS} in characteristic $2$. Finally, ${}^2\!F_4(2)$ contains no other group in Table \ref{tab:e6curlyS}.

Finally, for $G={}^2\!E_6(q)$, $X$ is one of $Fi_{22}$ and ${}^2\!F_4(2)'$. The group $Fi_{22}$ contains $\Omega_7(3)$ (which is why it is a novelty maximal), and ${}^2\!F_4(2)'$ and $\PSL_2(8)$, but $p=2$ and neither of these groups is in Table \ref{tab:2e6curlyS}. If $X$ is ${}^2\!F_4(2)'$ then again there are no possibilities for $H$.

\medskip

Thus we may assume that $X$ is one of the groups from Tables \ref{tab:f4othermaximalsqodd} to \ref{tab:2e6othermaximals}. Notice that none of our possibilities for $H$ has a projective representation of dimension less than $5$, so subgroups of type $A_1$ and $A_2$ can be ignored. We can also of course exclude subfield subgroups like $F_4(q_0)\leq F_4(q)$ and ${}^\varepsilon\!E_6(q_0)\leq {}^\pm E_6(q)$. We will deal with the possibility $H\leq G_2(q)$ later, so ignore this case for now.

If $G=F_4(q)$ and $X=B_4(q)$ then $X$ stabilizes a line on $M(F_4)^\circ$ unless $p=3$. The same holds for the $D_4$ and ${}^3\!D_4$ subgroups, so these can be eliminated. If $p=3$ on the other hand, $B_4$ acts with factors $9,16$, and this is incompatible with the factors of the groups in Table \ref{tab:f4curlyS} for $p=3$. The $\Sp_6$-parabolic stabilizes a $6$-space on $M(F_4)$ and none of the members of $\mathcal{S}$ does. The $B_3$-parabolic has a Levi subgroup contained in $B_4$, so this cannot occur. The group $\SL_3(3)$ is a minimal simple group, so $3^3\rtimes \SL_3(3)$ cannot contain $H$. Apart from $G_2$ subgroups, all other options for $X$ have been eliminated.

If $G={}^\varepsilon\!E_6(q)$, then parabolic subgroups stabilize spaces of (co)dimension $1$, $2$, $3$ or $6$, and $H$ does not, and $X$ cannot be $F_4(q)$ as then $H$ stabilizes a line on $M(E_6)$. All of the $D_4$ type subgroups have connected component inside $D_5$, so can be ignored. The $3$-local subgroup $3^{3+3}\rtimes \SL_3(3)$ again cannot be $X$. Apart from $G_2$ subgroups, we have the Weyl group $W(E_6)$ and the subgroup $\PSp_8(q)$. The Weyl group has simple factor $\PSp_4(3)\cong \PSU_4(2)$, and does not contain any of our subgroups $H$.

The group $X=\PSp_8(q)$ when $p$ is odd acts irreducibly on $M(E_6)$, so could contain some of our groups. The groups $H$ that can embed in $X$ are $\PSL_2(r)$ for $r=8,11,13$. Since the minimal degree for a non-trivial projective representation for $H$ is $5$, the restriction of $M(C_4)$ to (a central extension of) $H$ has a single non-trivial factor. If it has dimension less than $7$, then $H$ stabilizes a line and a hyperplane on $M(C_4)$, and we see that $H$ must stabilize a line on $M(E_6)$ (note that $S^2(M(C_4))=1\oplus M(E_6){\downarrow_{C_4}}$ as $p$ is odd). If it has dimension $7$, then $H$ acts on $M(C_4)$ as $1\oplus 7$, but this cannot stabilize an alternating form. Thus $H$ acts irreducibly on $M(C_4)$. But the three choices for $r$ cannot yield an irreducible embedding $H\to X$. (The group $\PSL_2(17)$ \emph{can} be irreducibly embedded in $\Sp_8(q)$, and this is why it does not appear in Tables \ref{tab:e6curlyS} or \ref{tab:2e6curlyS}.)

\medskip

The last remaining case is when $X$ has a $G_2$ factor. If $X$ is of type $A_1G_2\leq F_4$ or $A_2G_2\leq E_6$ then $H$ embeds in the $G_2$ factor alone. But this stabilizes lines on $M(F_4)$ and $M(E_6)$ respectively, so this cannot be the case. Thus $H$ embeds in a $G_2$ that acts irreducibly on the minimal module. The possibilities for $H$ are tabulated in \cite[Tables 8.41--8.43]{bhrd}, and are $\PSL_2(8)$ and $\PSL_2(13)$. But the irreducible $G_2$ has already been noted to contain these options for $H$ in their respective sections, so no new cases can emerge.

This completes the proof that there are no overgroups of the members of $\mathcal{S}$ listed in the tables.

\bigskip

\bigskip\textbf{Acknowledgments}: The author gratefully acknowledges Michael Aschbacher for allowing the inclusion of his work on the subgroup $\SL_2(8)$ from \cite{aschbacherE6Vun}, without which it seems possible that the normalizer of that subgroup would not have been computed at all. The author would also like to thank the referee of this manuscript, who provided a large number of detailed and helpful comments on the paper, and significantly improved its exposition.

This version also acknowledges the people who have found errors in this paper since its publication. They are listed in the appendix along with the particular error they found.

\appendix

\section{Changes from published version}

\begin{enumerate}
\item On the bottom of p638, I use ${}^2\!F_4(\sqrt q)$ rather than the more standard notation ${}^2\!F_4(q)$, which I also use elsewhere. This has been removed. (See also \ref{error:Tits} below.)

\item In Table 8 there is a torus normalizer with structure $(q^2+1)^2\cdot (4 \circ \mathrm{GL}_2(3))$. The group $4 \circ \mathrm{GL}_2(3)$ is not a subgroup of the Weyl group of type $F_4$, and the correct group should be $\mathrm{SL}_2(3)\rtimes 4$, which is the centralizer of an element of order $4$ in $W(F_4)$. The original paper of Liebeck--Saxl--Seitz also contains this error, and I copied their tables over without double-checking.

Thanks to Mikko Korhonen for noticing this.
\item\label{error:Tits} In Table 8 the large Ree subgroups have a badly explained structure. I have written ${}^2\!F_4(q_0)$ in analogy with the fact that $\mathrm{SU}_n(q)\leq \mathrm{SL}_n(q^2)$. This doesn't really work for large Ree groups, and the more standard notation is that ${}^2\!F_4(2)\leq F_4(2)$, for example. So this line should read that ${}^2\!F_4(q)\leq F_4(q)$.

Thanks to Tim Burness for mentioning this.

More importantly, when writing this erratum I noticed that the condition that $q$ is an odd power of $2$ is missing from this, but since there is no such group as ${}^2\!F_4(4)$ this is heavily implied, even if it is an error that it is not present.

\item In Table 9 the novelty subgroup should be $[q^{31}].(\SL_2(q) \times \SL_2(q) \times SL_3(q)).(q-1)^2/e$. There's a missing $(q-1)$ factor in the torus (obvious) and the action when removing the centre is to quotient by a diagonal $e$ in the $\SL_3(q) \times (q-1)^2$ subgroup. This yields a central product, so we can write the structure above.

Thanks to Tim Burness for highlighting this.

\item In Table 9 the subgroup $d^2.(\mathrm{P}\Omega_8^+(q)\times ((q-1)/d)^2/e).d^2.\Sym(3)$ is correct. Comparing this with Table 10, we find the correct version should be $d^2.(\mathrm{P}\Omega_8^+(q)\times ((q + 1)/d)^2/e').d^2.\Sym(3)$, not the $d^2.(\mathrm{P}\Omega_8^+(q)\times ((q + 1)/d)^2).d^2.\Sym(3)$ as given.

Thanks to Melissa Lee for noticing this.

\end{enumerate}

\providecommand{\bysame}{\leavevmode\hbox to3em{\hrulefill}\thinspace}

\end{document}